\input amstex
\input xy
\xyoption{all}
\documentstyle{amsppt}
\document

\magnification 1000

\def\gen{\frak{g}}

\def\den{\frak{d}}
\def\ten{\frak{t}}

\def\een{\frak{e}}

\def\len{\frak{l}}

\def\nen{\frak{n}}

\def\men{\frak{m}}

\def\Sen{\frak{S}}

\def\a{{\alpha}}
\def\g{{\gamma}}

\def\l{{{\lambda}}}

\def\eps{{\varepsilon}}

\def\ab{{\bold a}}
\def\bb{{\bold b}}

\def\eb{{\bold e}}
\def\fb{{\bold f}}

\def\ib{{\bold i}}
\def\ibt{{{\yb}}}

\def\jb{{\bold j}}
\def\kb{{\bold k}}

\def\mb{{\bold m}}

\def\forb{{{\bold{for}}}}

\def\yb{{\bold y}}
\def\Ab{{\bold A}}

\def\Bb{{\bold B}}

\def\Db{{\bold D}}
\def\Eb{{\bold E}}
\def\Fb{{\bold F}}

\def\Modb{{\bold{Mod}}}
\def\modb{{\bold{mod}}}
\def\fmodb{{\bold{fmod}}}
\def\fModb{{\bold{fMod}}}

\def\projb{{\bold{proj}}}
\def\Dcb{{\pmb{\Dc}}}
\def\Pcb{{\pmb{\Pc}}}
\def\Qcb{{\pmb{\Qc}}}

\def\Lbb{{\pmb{{\Lambda}}}}
\def\Gb{{\bold G}}
\def\zb{{\bold z}}

\def\Kb{{\bold K}}

\def\Nb{{\bold N}}
\def\Hb{{\bold H}}

\def\Qb{{\bold Q}}
\def\Rb{{\bold R}}

\def\Sb{{\bold S}}

\def\Lb{{\bold L}}

\def\Vb{{\bold V}}

\def\Wb{{\bold W}}
\def\Yb{{\bold Y}}
\def\Zb{{\bold Z}}

\def\k{{\kb}}

\def\B{{\roman B}}
\def\C{{\roman C}}

\def\A{{\roman A}}
\def\B{{\roman B}}
\def\D{{\roman D}}

\def\R{{\roman R}}
\def\R{{\roman R}}

\def\th{{{}^\theta}}

\def\top{{\text{top}}}
\def\gr{\roman{gr}}

\def\Ext{\roman{Ext}}

\def\Hom{{\roman{Hom}}}
\def\hom{{\roman{hom}}}
\def\gHom{{\roman{hom}}}

\def\sup{{\roman{sup}}}

\def\Eu{{\roman{eu}}}

\def\dim{{\roman{dim}}}
\def\gdim{{\roman{gdim}}}

\def\Ind{{\roman{Ind}}}

\def\End{{\roman{End}}}
\def\Hom{{\roman{Hom}}}

\def\eu{{\roman{eu}}}
\def\ch{{\roman{ch}}}

\def\Sum{\ts\sum}

\def\Prod{\ts\prod}

\def\Lie{\roman{Lie}}

\def\proj{{{\bold{{proj}}}}}

\def\CC{{\Bbb C}}

\def\NN{{\Bbb N}}

\def\PP{{\Bbb P}}
\def\QQ{{\Bbb Q}}

\def\ZZ{{\Bbb Z}}

\def\Ac{{\Cal A}}
\def\Bc{{\Cal B}}
\def\Cc{{\Cal C}}
\def\Dc{{\Cal D}}
\def\Ec{{\Cal E}}
\def\Fc{{\Cal F}}

\def\Hc{{\Cal H}}

\def\Kc{{\Cal K}}
\def\Lc{{\Cal L}}

\def\Oc{{\Cal O}}
\def\Pc{{\Cal P}}
\def\Qc{{\Cal Q}}
\def\Rc{{\Cal R}}

\def\Vc{{\Cal V}}
\def\Vcb{{\pmb{\Vc}}}

\def\Ccb{{\pmb{\Cc}}}
\def\Bcb{{\Bb}}

\def\Zc{{\Cal Z}}

\def\mod{{{\roman{mod}}}}

\def\soc{\roman{soc}}
\def\top{\roman{top}}

\def\and{{\text{and}}}
\def\low{{\text{low}}}
\def\up{{\text{up}}}

\def\max{{\roman{max}}}

\def\ds{\displaystyle}
\def\ts{\textstyle}

\def\qed{\hfill $\sqcap \hskip-6.5pt \sqcup$}        
\overfullrule=0pt                                    

\def\7dag{{{\!\!\!\!\!\!\!\dag}}}
\def\6dag{{{\!\!\!\!\!\!\dag}}}
\def\5dag{{{\!\!\!\!\!\dag}}}
\def\4dag{{{\!\!\!\!\dag}}}
\def\3dag{{{\!\!\!\dag}}}
\def\2dag{{{\!\!\dag}}}
\def\1dag{{{\!\dag}}}

\def\la{{\langle}}
\def\ra{{\rangle}}

\newdimen\Squaresize\Squaresize=14pt
\newdimen\Thickness\Thickness=0.5pt
\def\Square#1{\hbox{\vrule width\Thickness
          \alphaox to \Squaresize{\hrule height \Thickness\vss
          \hbox to \Squaresize{\hss#1\hss}
          \vss\hrule height\Thickness}
          \unskip\vrule width \Thickness}
          \kern-\Thickness}
\def\Vsquare#1{\alphaox{\Square{$#1$}}\kern-\Thickness}

\nologo

\topmatter
\title Canonical bases and affine Hecke algebras of type B 
\endtitle
\rightheadtext{} \abstract
We prove a series of conjectures of Enomoto and Kashiwara 
on canonical bases and branching rules of affine
Hecke algebras of type B. The main ingredient of the proof is
a new graded Ext-algebra associated with 
quiver with involutions that we compute 
explicitly.
\endabstract
\author M. Varagnolo, E. Vasserot\endauthor
\address 
Universit\'e de Cergy-Pontoise, UMR CNRS 8088,
F-95000 Cergy-Pontoise,\endaddress \email
michela.varagnolo\@math.u-cergy.fr\endemail
\address D\'epartement de Math\'ematiques,
Universit\'e Paris 7, 175 rue du Chevaleret, F-75013 Paris,
\endaddress
\email vasserot\@math.jussieu.fr\endemail
\thanks
2000{\it Mathematics Subject Classification.} Primary ??; Secondary
??.
\endthanks
\endtopmatter
\document

\head Introduction\endhead
A new family of graded algebras, called KLR algebras, has been recently introduced in \cite{KL}, \cite{R}.
These algebras yield a categorification of $\fb$, the negative part 
of the quantized enveloping algebra of any type.
In particular, one can obtain a new interpretation of the canonical bases, see \cite{VV}. In type A 
or $\A^{(1)}$ the KLR algebras are Morita equivalent to the affine Hecke 
algebras and their cyclotomic quotients. Hence they give a new way
to understand the categorification of the simple highest weight modules  and
the categorification of $\fb$ 
via some Hecke algebras of type A or $\A^{(1)}$. 
See \cite{BK} for instance.
One of the advantages of KLR algebras is that they are graded, while the affine Hecke algebras are not.
This explain why KLR algebras are better adapted than affine Hecke algebras to describe canonical bases.
Indeed one could view KLR algebras as an intermediate object between the representation theory of
affine Hecke algebras and its Kazhdan-Lusztig geometric counterpart in term of perverse sheaves.
This is central in \cite{VV}, where KLR algebras are proved to be isomorphic to the
Ext-algebras of some complex of constructible sheaves.

In the other hand, the (branching rules
for) affine Hecke algebras of type B have been investigated
quite recently, see \cite{E}, \cite{EK1,2,3}, \cite{Ka}, \cite{M}. 
Lusztig's description of the canonical basis of
$\fb$ in type $\A^{(1)}$ in \cite{L1} implies that this basis can be naturally 
identified with the set of isomorphism classes of 
simple objects of a category of modules of the affine Hecke algebras of type A.
This identification was mentioned in \cite{G}, and was used
in \cite{A}. 
More precisely, there is a linear isomorphism between 
$\fb$ and the
Grothendieck group of finite dimensional modules of the
affine Hecke algebras of type A, and it is proved in \cite{A}
that the induction/restriction functors
for affine Hecke algebras are given by the action of the Chevalley generators  and 
their transposed operators with respect to 
some symmetric bilinear form on $\fb$. 
In \cite{E}, \cite{EK1,2,3} a similar behavior is conjectured and 
studied for affine Hecke algebras of type B. Here $\fb$ is replaced by
an explicit module $\th\Vb(\l)$ over an explicit algebra $\th\Bcb$.
First, it is
conjectured that $\th\Vb(\l)$ admits a canonical basis.
Next, it is conjectured that this basis is naturally identified with
the set of isomorphism classes of
simple objects of a category of modules of the affine Hecke algebras of type B.
Further, in this identification the 
branching rules of the affine Hecke algebras of type B
are given by the $\th\Bcb$-action on
$\th\Vb(\l)$. The first conjecture has been proved in \cite{E} under
the restrictive assumption that $\l=0$. 
Here we prove the whole set of conjectures.
Indeed, our construction is slightly more general, see the appendix.

Roughtly speaking our argument is as follows.
In \cite{E} a geometric description of the canonical basis of
$\th\Vb(0)$ was given. This description is similar to Lusztig's 
description of the canonical basis of $\fb$ via perverse sheaves on the moduli
stack of representations of some quiver. It is given in terms of
perverse sheaves on the moduli
stack of representations of a quiver with involution.
First we give a analogue of this for $\th\Vb(\l)$ for any $\l$.
This yields the existence of a canonical basis $\th\!\Gb^\low(\l)$
for $\th\Vb(\l)$ for arbitrary $\l$.
Then we compute explicitely the Ext-algebras between complexes
of constructible sheaves naturally attached to quivers with involutions.
These complexes enter in a natural way in the definition of $\th\!\Gb^\low(\l)$.
This computation yields a new family of graded algebras $\th\Rb_m$
where $m$ is a nonnegative integer. We prove that the algebras $\th\Rb_m$
are Morita equivalent to the affine Hecke algebras of type B.
Finally we describe $\th\Vb(\l)$ and the basis $\th\!\Gb^\low(\l)$ in terms
of the Grothendieck group of $\th\Rb_m$.

The plan of the paper is the following.
Section 1 contains some basic notation for Lusztig's theory of perverse sheaves on the moduli
stack of representations of quivers. Section 2 yields similar notation for the case of quivers with involutions.
Our setting is more general than in \cite{E}, where only the case $\l=0$ is considered.
In Section 3 we introduce the convolution algebra associated with a quiver with involution.
The main result of Section 4 is Theorem 4.17 
where the polynomial representation of the Ext-algebra
$\Zb^\delta_{\Lbb,\Vb}$
associated with a quiver with involution is computed. 
Here $\Lbb$ is a $I$-graded $\CC$-vector space of 
dimension vector $\l\in\NN I$, while
$\Vb$ is a $I$-graded $\CC$-vector space with a non-degenerate symmetric 
bilinear form of dimension vector $\nu\in\NN I$.
The polynomial representation of 
$\Zb^\delta_{\Lbb,\Vb}$
is faithful.
In Section 5 we give the main properties of the graded algebra $\th\Rb(\Gamma)_{\l,\nu}$.
In Section 6 we introduce the affine Hecke algebra of type B 
and we prove that it is Morita equivalent
to $\th\Rb_m$, a specialization of $\th\Rb(\Gamma)_{\l,\nu}$.
Section 7 is a reminder on KLR algebras and on the main result of \cite{VV}.
In Section 8 we categorify the module $\th\Vb(\l)$ from \cite{EK1} using the graded algebra $\th\Rb_m$.
In Section 9 we prove the isomorphism
$\th\Rb(\Gamma)_{\l,\nu}=\Zb^\delta_{\Lbb,\Vb}$. This is essential to compare the construction from Section 8 with that in Section 10. In Section 10 we give a categorification of $\th\Vb(\l)$ ``\`a la Lusztig'' in terms of perverse sheaves on the moduli stack of representations of quivers with involution. This is essentially the same construction as in \cite{E}. However, since we need a more general setting than in loc.~ cit.~ we have briefly
reproduced the main steps of the construction. One of our initial motivations was to give a completely algebraic proof of the conjectures, without any perverse sheaves at all. We still do not know how to do this.
The main result of the paper is Theorem 10.19.

The same technic yields similar results for 
affine Hecke algebras of any classical type. The case of type D is done in \cite{SVV},
the case of type C is done in the appendix.
Note that the idea to use canonical bases technics 
to study affine Hecke algebras in non A type is not new, 
see \cite{L3}, \cite{L4}.
At the moment we do not know the precise relation between loc.~cit.~
and our approach.

{\it Acknowledgment.}
We are grateful to M. Kashiwara and G. Lusztig for some remark 
on the material of this paper.

\head Contents\endhead
\item {0.}~Notation
\vskip 1mm 
\item {1.}~Reminder on quivers and extensions
\vskip 1mm 
\item {2.}~Quivers with involutions
\vskip 1mm 
\item {3.}~The convolution algebra
\vskip 1mm 
\item {4.}~The polynomial representation of the graded algebra 
$\th\Zb^\delta_{\Lbb,\Vb}$
\vskip 1mm 
\item {5.}~The graded $\kb$-algebra $\th\Rb(\Gamma)_{\l,\nu}$
\vskip 1mm 
\item {6.}~Affine Hecke algebras of type B
\vskip 1mm 
\item {7.}~Global bases of $\fb$ and projective graded modules of KLR algebras
\vskip 1mm 
\item {8.}~Global bases of $\th\Vb(\l)$ and projective graded $\th\Rb$-modules
\vskip 1mm 
\item {9.}~Presentation of the graded algebra
$\th\Zb^\delta_{\Lbb,\Vb}$
\vskip 1mm 
\item {10.}~Perverse sheaves on $\th\Eb_{\Lbb,\Vb}$ 
and the global bases of $\th\Vb(\l)$
\vskip 1mm 
\item {A.}~Appendix
\vskip 1mm 
\item {}~Index of notation

\vskip2cm

\head 0.~Notation\endhead

\subhead 0.1.~Combinatorics\endsubhead Given a positive integer $m$ and a tuple
$\mb=(m_1,m_2,\dots m_r)$ of positive integers we write $\Sen_m$ for
the symmetric group and $\Sen_\mb$ for the group
$\prod_{l=1}^r\Sen_{m_l}$. Set
$$|\mb|=\sum_{l=1}^rm_l,\quad
\ell_\mb=\sum_{l=1}^r\ell_{m_l},\quad\ell_m=m(m-1)/2.$$ We use the
following notation for $v$-numbers
$$\la m\ra=\sum_{l=1}^mv^{m+1-2l},
\quad \la m\ra!=\prod_{l=1}^m\la l\ra, \quad
\left\langle\matrix m+n\cr n\endmatrix\right\rangle={\la m+n\ra!\over\la m\ra!\la n\ra!},\quad
\la\mb\ra!=\prod_{l=1}^r\la m_l\ra!.$$ 
Given two tuples $\mb=(m_1,m_2,\dots m_r)$,
$\mb'=(m'_1,m'_2,\dots m'_{r'})$ we define the tuple
$$\mb\mb'=(m_1,m_2,\dots m_r,m'_1,m'_2,\dots m'_{r'}).$$

\vskip1mm

\subhead 0.2.~Graded modules over graded algebras\endsubhead 
Let $\kb$ be an algebraically closed field of caracteristic 0.
Let $\Rb=\bigoplus_d\Rb_d$ be a graded $\kb$-algebra.
Unless specified otherwise the word graded we'll always mean
$\ZZ$-graded.
Let $\Rb$-$\modb$ be the category of finitely
generated graded $\Rb$-modules, $\Rb$-$\fmodb$ be the full
subcategory of finite-dimensional graded modules and $\Rb$-$\proj$ be the
full subcategory of projective objects. Unless specified otherwise a module
is always a left module. We'll abbreviate
$$K(\Rb)=[\Rb\text{-}\proj],\quad
G(\Rb)=[\Rb\text{-}\fmodb].$$ 
Here $[\Ccb]$ denotes the Grothendieck group
of an exact category $\Ccb$. 
Assume that the $\kb$-vector spaces $\Rb_d$ are finite dimensional for each $d$.
Then $K(\Rb)$ is a
free Abelian group with a basis formed by the isomorphism classes of
the indecomposable objects in $\Rb$-$\proj$,
and $G(\Rb)$ is a free Abelian group with a basis
formed by the isomorphism classes of the simple objects in
$\Rb$-$\fmodb$. Given an object $M$ of $\Rb$-$\proj$ or
$\Rb$-$\fmodb$ let $[M]$ denote its class in $K(\Rb)$, $G(\Rb)$
respectively. 
When there is no risk of confusion we abbreviate $M=[M]$.
We'll write $[M:N]$ for the composition multiplicity of the $\Rb$-module $N$
in the $\Rb$-module $M$.
Consider the ring $\Ac=\ZZ[v,v^{-1}]$.
If the grading of $\Rb$ is bounded below then
the $\Ac$-modules $K(\Rb)$, $G(\Rb)$ are free. Here
$\Ac$ acts on $G(\Rb)$, $K(\Rb)$ as follows
$$vM=M[1],\quad v^{-1}M=M[-1].$$
For any $M,N$ in $\Rb$-$\modb$ let 
$$\gHom_\Rb(M,N)=\bigoplus_d\Hom_\Rb(M,N[d])$$
be the graded $\kb$-vector space of all $\Rb$-module homomorphisms.
If $\Rb=\kb$
we'll omit the subscript $\Rb$ in hom's and in tensor products. 
As much as possible we'll use the following convention : 
graded objects are denoted by minuscules and non-graded ones by majuscules. 
In particular $\Rb$-$\Modb$ will denote the category of finitely generated 
(non-graded) $\Rb$-modules. We 'll abbreviate
$$\Hom=\Hom_\kb,\quad\otimes=\otimes_\kb,\quad etc.$$
For a graded $\kb$-vector space $M=\bigoplus_dM_d$ we'll write
$$\gdim(M)=\sum_dv^d\dim(M_d),$$
where $\dim$ is the dimension over $\kb$.

\vskip3mm

\subhead 0.3.~Constructible sheaves\endsubhead
Given an action of a complex linear
algebraic group $G$ on a quasiprojective algebraic variety $X$ over
$\CC$ we write $\Dcb_G(X)$ for the bounded derived category of
complexes of $G$-equivariant sheaves of $\k$-vector spaces on $X$.
Objects of $\Dcb_G(X)$ are referred to as complexes. If $G=\{e\}$,
the trivial group, we abbreviate $\Dcb(X)=\Dcb_G(X)$. For each
complexes $\Lc$, $\Lc'$ we'll abbreviate
$$\Ext_G^*(\Lc,\Lc')=\Ext^*_{\Dcb_G(X)}(\Lc,\Lc'),\quad
\Ext^*(\Lc,\Lc')=\Ext^*_{\Dcb(X)}(\Lc,\Lc')$$ if no confusion is
possible. The constant sheaf on $X$ with stalk $\kb$ will be denoted $\k$. For
any object $\Lc$ of $\Dcb_G(X)$ let $H_G^*(X,\Lc)$ be the space of
$G$-equivariant cohomology with coefficients in $\Lc$. Let $\Dc\in
\Dcb_G(X)$ be the $G$-equivariant dualizing complex, see \cite{BL,
def.~3.5.1}. For each $\Lc$ let 
$$\Lc^\vee=\Hc om(\Lc,\Dc)$$
be its Verdier dual, where $\Hc om$ is the internal Hom. Recall that
$$(\Lc^\vee)^\vee=\Lc,\quad
\Ext^*_G(\Lc,\Dc)=H^*_G(X,\Lc^\vee),\quad
\Ext^*_{G}(\k,\Lc)=H_G^*(X,\Lc).$$ We define the space of
$G$-equivariant homology by
$$H^G_*(X,\k)=H^*_G(X,\Dc).$$
Note that $\Dc=\k[2d]$ if $X$ is a smooth $G$-variety
of pure dimension $d$. Consider the following graded $\kb$-algebra
$$\Sb_G=H^{*}_G(\bullet,\k).$$
The graded $\k$-vector space $H^G_*(X,\k)$ has a natural structure
of a graded $\Sb_G$-module. We have
$$H_*^G(\bullet,\k)=\Sb_G$$ as graded
$\Sb_G$-module. There is a canonical graded $\k$-algebra isomorphism
$$\Sb_G\simeq\k[\gen]^G.$$ Here the symbol $\gen$ denotes the Lie
algebra of $G$ and a $G$-invariant homogeneous polynomial over
$\gen$ of degree $d$ is given the degree $2d$ in $\Sb_G$.

Fix a morphism of quasi-projective algebraic $G$-varieties $f:X\to
Y$. If $f$ is a proper map there is a direct image homomorphism
$$f_*:H_*^G(X,\k)\to H_*^G(Y,\k).$$ If $f$ is a smooth map of relative
dimension $d$ there is an inverse image homomorphism
$$f^*:H_i^G(Y,\k)\to H_{i-2d}^G(X,\k),\quad\forall i.$$
If $X$ has pure dimension $d$ there is a natural homomorphism
$$H_G^i(X,\k)\to H^G_{i-2d}(X,\k).$$ It is invertible if $X$ is
smooth. The image of the unit is called the fundamental class of $X$
in $H_*^G(X,\k)$. We denote it by $[X]$. If $f:X\to Y$ is the
embedding of a $G$-stable closed subset and $X'\subset X$ is the
union of the irreducible components of maximal dimension then the
image of $[X']$ by the map $f_*$ is the fundamental class of $X$ in
$H_*^G(Y,\k)$. It is again denoted by $[X]$.

\vskip2cm

\head 1. Reminder on quivers and extensions
\endhead

\subhead 1.1.~Representations of quivers\endsubhead We assume given
a nonempty quiver $\Gamma$ such that no arrow may join a
vertex to itself. Recall that $\Gamma$ is a tuple $(I,H,h\mapsto
h',h\mapsto h'')$ where $I$ is the set of vertices, $H$ is
the set of arrows, and for $h\in H$ the vertices $h',h''\in I$
are the origin and the goal of $h$ respectively. Note that the set $I$ may be infinite.
For $i,j\in I$
we write
$$H_{i,j}=\{h\in H;h'=i,h''=j\}.$$ We'll abbreviate $i\to j$ for
$H_{i,j}\neq\emptyset$, $i\not\to j$ for $H_{i,j}=\emptyset$, and
$h:i\to j$ for $h\in H_{i,j}$. Let $h_{i,j}$ be the number of
elements in $H_{i,j}$ and set
$$i\cdot j=-h_{i,j}-h_{j,i},\quad i\cdot i=2,\quad i\neq j.$$

Let $\Vcb$ be the category of finite-dimensional $I$-graded
$\CC$-vector spaces $\Vb=\bigoplus_{i\in I}\Vb_i$ with morphisms
being linear maps respecting the grading. For $\nu=\sum_i\nu_i
i$ in $\NN I$ let $\Vcb_\nu$ be the full subcategory of $\Vcb$ whose
objects are those $\Vb$ such that $\dim(\Vb_i)=\nu_i$ for all $i$.
We call $\nu$ the dimension vector of $\Vb$.
Given an object $\Vb$ of $\Vcb$ let
$$E_\Vb=\bigoplus_{h\in H}\Hom(\Vb_{h'},\Vb_{h''}).$$
The algebraic group $G_\Vb=\Prod_iGL(\Vb_i)$ acts on $E_\Vb$ by
$(g,x)\mapsto gx=y$ where $y_h=g_{h''}x_hg_{h'}^{-1}$, $g=(g_i)$,
$x=(x_h)$, and $y=(y_h)$.

Fix a nonzero element $\nu$ of $\NN I$.
Let $Y^\nu$ be the set of all pairs $\ibt=(\ib,\ab)$ where
$\ib=(i_1,i_2,\dots i_k)$ is a sequence of elements of $I$ and
$\ab=(a_1,a_2,\dots a_k)$ is a sequence of positive integers such
that $\sum_la_l\,i_l=\nu$. Note that the assignment
$$\ibt\mapsto(a_1 i_1,a_2 i_2,\dots a_k i_k)\leqno(1.1)$$
identifies $Y^\nu$  
with a set of sequences 
$$\nu^1,\nu^2,\dots,\nu^k\in\NN I\quad\roman{with}\ 
\nu=\sum_{l=1}^k\nu^l.\leqno(1.2)$$
For each pair $\yb=(\ib,\ab)$ as above we'll call $\ab$ the multiplicity of
$\yb$.
Let $I^\nu\subset Y^\nu$
be the set of all pairs $\yb$ with multiplicity $(1,1,\dots,1)$.
We'll abbreviate $\ib$ for a pair $\yb=(\ib,\ab)$ which lies in $I^\nu$.
Given a positive integer $m$ we have 
$\bigsqcup_{\nu}I^\nu=I^m,$ where $\nu$ runs
over the set  of elements $\nu$ of $\NN I$ with $|\nu|=m$.
Here, we write $\nu=\sum_{i\in I}\nu_i i$ and $|\nu|=\sum_i\nu_i$.
In a similar way, we define $Y^m=\bigsqcup_{\nu}Y^\nu.$

\vskip3mm

\subhead 1.2.~Flags\endsubhead Let $\nu\in\NN I$, $\nu\neq 0$, and assume
that $\Vb$ lies in $\Vcb_\nu$.
For each sequence $\yb=(\nu^1,\nu^2,\dots\nu^k)$ 
as in (1.1), (1.2), a flag of type $\yb$ in $\Vb$ is a sequence
$$\phi=(\Vb=\Vb^0\supset\Vb^{1}\supset\cdots\supset\Vb^k=0)$$ of
$I$-graded subspace of $\Vb$ such that for any $l$ the $I$-graded
subspace $\Vb^{l-1}/\Vb^{l}$ belongs to $\Vcb_{\nu^l}$. Let
$F_{\Vb,\yb}$ be the variety of all flags of type $\yb$ in $\Vb$.
The group $G_\Vb$ acts transitively on $F_{\Vb,\yb}$ in the obvious
way, yielding a smooth projective $G_\Vb$-variety structure on
$F_{\Vb,\yb}$.

If $x\in E_\Vb$ we say that the flag $\phi$ is $x$-stable if
$x_h(\Vb^l_{h'})\subset\Vb^l_{h''}$ for all $h$, $l$. Let
$\widetilde F_{\Vb,\yb}$ be the variety of all pairs $(x,\phi)$ such
that $\phi$ is $x$-stable. Set
$d_{\yb}=\dim(\widetilde F_{\Vb,\yb}).$
The group $G_\Vb$ acts on $\widetilde F_{\Vb,\yb}$ by
$g:(x,\phi)\mapsto(gx,g\phi)$. The first projection
gives a $G_\Vb$-equivariant proper morphism
$$\pi_{\yb}:\widetilde F_{\Vb,\yb}\to E_\Vb.$$

\vskip3mm

\subhead 1.3.~Ext-algebras \endsubhead 
Let $\nu\in\NN I$, $\nu\neq 0$, and assume that
$\Vb\in\Vcb_\nu$. We abbreviate $\Sb_\Vb=\Sb_{G_\Vb}.$ For each
sequence $\yb\in Y^\nu$ we have the following semisimple complexes in
$\Dcb_{G_\Vb}(E_\Vb)$
$$\Lc_{\yb}=(\pi_{\yb})_!(\k),\quad\Lc_\yb^\vee=\Lc_\yb[2d_\yb],\quad
\Lc^\delta_\yb=\Lc_\yb[d_\yb].$$
For $\yb,\yb'$ in $Y^\nu$ we consider the graded
$\Sb_\Vb$-module
$$\Zb_{\Vb,\yb,\yb'}=
\Ext^*_{G_\Vb}(\Lc_\yb^\vee,\Lc_{\yb'}^\vee).$$ For $\yb,\yb',\yb''$
in $Y^\nu$ the Yoneda composition is a homogeneous
$\Sb_\Vb$-bilinear map of degree zero
$$\star:\Zb_{\Vb,\yb,\yb'}\times\Zb_{\Vb,\yb',\yb''}\to\Zb_{\Vb,\yb,\yb''}.$$
The map $\star$ equips the graded $\k$-vector space
$$\Zb_\Vb=\bigoplus_{\ib,\ib'\in
I^\nu}\Zb_{\Vb,\ib,\ib'}$$ with the structure of an
associative graded $\Sb_\Vb$-algebra with 1. If there is no ambiguity
we'll omit the symbol $\star$. We
set
$$\Fb_{\Vb,\yb}=\Ext^*_{G_\Vb}(\Lc_\yb^\vee,\Dc), \quad
\Fb_\Vb=\bigoplus_{\ib\in I^\nu}\Fb_{\Vb,\ib}.$$  For
$\yb,\yb'$ in $Y^\nu$ the Yoneda product gives a graded
$\Sb_\Vb$-bilinear map
$\Zb_{\Vb,\yb,\yb'}\times\Fb_{\Vb,\yb'}\to\Fb_{\Vb,\yb}.$ This
yields a left graded representation of $\Zb_\Vb$ on $\Fb_\Vb$. For
each $\ib\in I^\nu$ let $1_{\Vb,\ib}\in\Zb_{\Vb,\ib,\ib}$ denote the
identity of $\Lc_\ib$. The elements $1_{\Vb,\ib}$ form a complete
set of orthogonal idempotents of $\Zb_\Vb$ such that
$$\Zb_{\Vb,\ib,\ib'}=1_{\Vb,\ib}\star\Zb_\Vb\star 1_{\Vb,\ib'},\quad
\Fb_{\Vb,\ib}=1_{\Vb,\ib}\star\Fb_\Vb.$$
We'll change the grading of $\Zb_{\Vb}$ in the following way. Put
$$
\Zb_{\Vb,\ib,\ib'}^\delta= \Ext^*_{
G_\Vb}(\Lc_{\ib}^\delta,\Lc_{\ib'}^\delta),\quad
\Zb_{\Vb}^\delta= 
\bigoplus_{\ib,\ib'\in
I^\nu}\Zb^\delta_{\Vb,\ib,\ib'}.$$
The graded $\kb$-algebra $\Zb_\Vb^\delta$ depends only on the dimension 
vector of $\Vb$. We'll write $$\Rb(\Gamma)_\nu=\Zb_\Vb^\delta.$$
This graded $\kb$-algebra has been computed explicitely in \cite{VV}.
The same result has also been anounced by R. Rouquier.
See Section 7 for more details.
We set also $I^0=\{\emptyset\}$, 
$\Lc^\delta_\emptyset=\kb$ (the constant sheaf over $\{0\}$) 
$$\Rb(\Gamma)_0=\Zb_{\{0\}}^\delta=\kb.$$

\vskip2cm

\head 2.~Quivers with involutions
\endhead

In this section we introduce an analogue of the Ext-algebra $\Rb(\Gamma)_\nu$.
It is associated with a quiver with an involution.

\subhead 2.1.~Representations of quivers with involution\endsubhead
Fix a  nonempty quiver $\Gamma$ such that no arrow may join a
vertex to itself. 
An involution $\theta$ on $\Gamma$ is a pair of involutions on $I$ and $H$, 
both denoted by $\theta$, 
such that 
the following properties hold 
for $h\in H$
\vskip2mm
\itemitem{$\bullet$}
$\theta(h)'=\theta(h'')$ and $\theta(h)''=\theta(h')$,
\vskip2mm
\itemitem{$\bullet$}
$\theta(h')=h''$ iff $\theta(h)=h$.
\vskip2mm
\noindent {\it We'll always assume that $\theta$ has no fixed points in $I$}, i.e.,
there is no $i\in I$ such that $\theta(i)=i$. To simplify we'll say
that $\theta$ has no fixed points.

Let $\th\Vcb$ be the category of finite-dimensional $I$-graded
$\CC$-vector spaces $\Vb$ with a non-degenerate symmetric
bilinear form $\varpi$ such that 
$$(\Vb_i)^\perp=\bigoplus_{j\neq\theta(i)}\Vb_j.$$ 
To simplify we'll say that 
$\Vb$ belongs to $\th\Vcb$ if there is a bilinear form $\varpi$ such that
the pair $(\Vb,\varpi)$ lies in $\th\Vcb$. The morphisms in
$\th\Vcb$ are the linear maps which respect the grading and
the bilinear form.
Let $$\th\NN I=\{\nu=\Sum_i\nu_i\,i\in\NN I
;\nu_{\theta(i)}=\nu_i,\forall i\}.$$ For  $\nu\in\th\NN I$ let $\th\Vcb_\nu$ be the full subcategory of
$\th\Vcb$ consisting of the pairs $(\Vb,\varpi)$ such that
$\Vb$ lies in $\Vcb_\nu$. Note that $|\nu|$ is an even integer.
We'll usually write $|\nu|=2m$ with $m\in\NN$.
Given $\Vb$ in $\th\Vcb$ and $\Lbb$ in $\Vcb$ we let
$$\gathered
\th\! E_\Vb=\{x=(x_h)\in E_\Vb;x_{\theta(h)}=-{}^tx_h,\,\forall h\in H\},\vspace{2mm}
\th\! G_\Vb=\{g\in
G_\Vb;g_{\theta(i)}={}^tg_i^{-1},\forall i\in I\},\vspace{2mm}
\th\! E_{\Lbb,\Vb}=\th\! E_\Vb\times L_{\Lbb,\Vb},\quad L_{\Lbb,\Vb}=\Hom_\Vcb(\Lbb,\Vb).
\endgathered$$
The algebraic groups $\th\! G_\Vb$, $G_\Lbb$
act on $\th\! E_\Vb$, $L_{\Lbb,\Vb}$
in the obvious way.

\vskip3mm

\subhead 2.2.~Generalities on isotropic flags\endsubhead 
Given a finite dimensional
$\CC$-vector space $\Wb$ with a non-degenerate symmetric
bilinear form $\varpi$, an {\it isotropic flag in $\Wb$} is a sequence of
subspaces
$$\phi=(\Wb=\Wb^{-k}\supset\Wb^{1-k}\supset\cdots\supset\Wb^{k}=0)$$
such that $(\Wb^l)^\perp=\Wb^{-l}$ for any $l=-k,1-k,\dots,k-1,k$. 
Here the symbol $\perp$ means the orthogonal relative to $\varpi$.
In particular $\Wb^0$ is a Lagrangian subspace of $\Wb$. Let
$F(\Wb)$ be the variety of all complete flags in $\Wb$, and
$F(\Wb,\varpi)$ be the subvariety of all complete isotropic flags,
i.e., we require that $\phi=(\Wb^l)$ is an isotropic flag such
that $\Wb^l$ has the dimension $m-l$ and $k=m$. If $\Wb$ has
dimension $2m$ then $F(\Wb,\varpi)$ has dimension $2\ell_m=m(m-1)$.

\vskip3mm

\subhead 2.3.~Sequences\endsubhead Fix a nonzero
dimension vector $\nu$ in $\th\NN I$. 
Let $\th Y^\nu$ be the set of all pairs $\yb=(\ib,\ab)$ in $Y^\nu$ such that
$$\ib=(i_{1-k},\dots,i_{k-1}, i_{k}),\quad\ab=(a_{1-k},\dots,a_{k-1},a_{k}),
\quad \theta(i_l)=i_{1-l},\quad a_l=a_{1-l}.$$ As in (1.1) we can identify
a pair $\yb$ as above with a sequence
$$\nu^{1-k},\dots,\nu^{k-1},\nu^{k}\in\NN I,\quad
\theta(\nu^l)=\nu^{1-l},\quad\sum_l\nu^l=\nu.$$
Let $\th\! I^\nu\subset \th Y^\nu$ be
the set of all pairs $\yb$ of multiplicity $(1,1,\dots, 1)$.
We'll abbreviate $\ib=(\ib,\ab)$ for each pair in $\th\!
I^\nu$. Note that a sequence in $\th I^\nu$ contains
$|\nu|=2m$ terms. 
Unless specified otherwise the entries of a sequence
$\ib$ in $\th\!I^\nu$ will always denoted by
$$\ib=(i_{1-m},\dots,i_{m-1}, i_{m}).$$ 
Finally, we set 
$$\th\! I^m=\bigcup_\nu\th\! I^\nu,
\quad\nu\in\th\NN I,\quad|\nu|=2m,$$
and we define $\th Y^m$ in the same way.

\vskip3mm

\subhead 2.4.~Definition of the map $\th\pi_{\Lbb,\yb}$\endsubhead
Fix $\nu\in\th\NN I$, $\nu\neq 0$, and $\l\in\NN I$.
Fix an object $\Vb$ in $\th\Vcb_\nu$ and an object
$\Lbb$ in $\Vcb_\l$.
For $\yb$ in $\th Y^\nu$ an {\it isotropic  flag of type $\yb$}
in $\Vb$ is an isotropic flag
$$\phi=(\Vb=\Vb^{-k}\supset\Vb^{1-k}\supset\cdots\supset\Vb^k=0)$$ 
such that $\Vb^{l-1}/\Vb^l$ lies in $\Vcb_{\nu^l}$ for each $l$. 
We define $\th\! F_{\Vb,\yb}$ to be
the variety of all isotropic flags of type $\yb$ in $\Vb$. 
Next, we define $\th\!\widetilde F_{\Lbb,\Vb,\yb}$
to be the variety of all tuples $(x,y,\phi)$ 
satisfying the 
following conditions :  
\vskip1mm
\itemitem{$\bullet$} 
$x\in \th\! E_\Vb$ and
$\phi\in\th\! F_{\Vb,\yb}$ is {\it stable by $x$}, i.e.,
$x(\Vb^l)\subset\Vb^l$ for each $l$, 
\vskip1mm
\itemitem{$\bullet$} 
$y\in L_{\Lbb,\Vb}$ and 
$y(\Lbb)\subset\Vb^0$. 
\vskip1mm
\noindent
We set $$d_{\l,\yb}=\dim(\th\!\widetilde F_{\Lbb,\Vb,\yb}).$$ 
We have the following formulas.

\proclaim{2.5.~Proposition} 
For $\ib\in\th\!I^\nu$ we have 

\vskip1mm

(a) $\dim(\th\! F_{\Vb,\ib})=\ell_{\nu}/2,$
\vskip2mm

(b) 
$d_{\l,\ib}=
\ell_\nu/2+\sum_{k<l;\,k+l\neq 1}h_{i_k,i_l}/2+\sum_{1\leqslant l}\l_{i_l}.$
\endproclaim

\noindent{\sl Proof :} Fix a subset $J\subset I$ such that 
$I=J\sqcup\theta(J)$.
Set $\Vb_J=\bigoplus_{j\in J}\Vb_j$.
The assignment
$(\Vb^k)\mapsto(\Vb^k\cap\Vb_J)$
takes $\th\! F_{\Vb,\ib}$ isomorphically onto 
$$\prod_{j\in J}F(\Vb_j).$$
Thus we have
$$\dim(\th\! F_{\Vb,\ib})=\sum_{j\in J}\ell_{\nu_j}=
\sum_{i\in I}\ell_{\nu_i}/2.$$
Next, fix a sequence $\ib$ as above and fix
a flag $\phi=(\Vb^k)$ in $\th\! F_{\Vb,\ib}$.
Then we have
$$d_{\l,\ib}=\ell_\nu/2+\dim\{x\in\th\! E_{\Vb};\,x(\Vb^k)\subset\Vb^k,\,
\forall k\}+\dim\{y\in L_{\Lbb,\Vb};\,
y(\Lbb)\subset\Vb^0\}.$$
Finally we have (see the discussion in Section 4.9)
$$\gathered
\dim\{x\in\th\! E_{\Vb};\,x(\Vb^k)\subset\Vb^k,\,\forall k\}=
\sum_{k<l;\,k+l\neq 1}h_{i_k,i_l}/2,
\vspace{1mm}\dim\{y\in L_{\Lbb,\Vb};\,
y(\Lbb)\subset\Vb^0\}=\sum_{1\leqslant l\leqslant m}\l_{i_l}.\endgathered$$

\qed

\vskip3mm

The group
$\th\! G_\Vb$ acts transitively on $\th\! F_{\Vb,\yb}$.
It acts also on $\th\! \widetilde F_{\Lbb,\Vb,\yb}$.
The first projection gives a $\th\!
G_\Vb$-equivariant proper morphism
$$\th\pi_{\Lbb,\yb}:
\th\!\widetilde F_{\Lbb,\Vb,\yb}\to \th\! E_{\Lbb,\Vb}.$$
For a future use we introduce also the obvious projection
$$p:\th\!\widetilde F_{\Lbb,\Vb}\to \th\! F_\Vb,\quad
\th\!\widetilde
F_{\Lbb,\Vb}=\coprod_{\ib\in \th\! I^\nu}
\th\!\widetilde F_{\Lbb,\Vb,\ib},\quad \th
F_\Vb=\coprod_{\ib\in\th\! I^\nu} \th\! F_{\Vb,\ib}.$$

\vskip3mm

\subhead 2.6.~Ext-algebras\endsubhead 
Let  $\l$, $\nu$, $\Lbb$, $\Vb$ be as above.
We abbreviate
$\th\Sb_\Vb=\Sb_{\th\! G_\Vb}.$ For 
$\yb\in\th\! Y^\nu$ we define the following semisimple complexes in
$\Dcb_{\th\! G_\Vb}(\th\! E_{\Lbb,\Vb})$
$$\th\! \Lc_{\yb}=(\th\! \pi_{\Lbb,\yb})_!(\k),
\quad\th\! \Lc_{\yb}^\vee=
\th\! \Lc_{\yb}[2d_{\l,\yb}],\quad
\th\!\Lc^\delta_{\yb}=\th\!\Lc_{\yb}[d_{\l,\yb}].$$ For
$\ib,\ib'$ in $\th\!I^\nu$ we consider the graded $\th\!
\Sb_\Vb$-module
$$\th \Zb_{\Lbb,\Vb,\ib,\ib'}=
\Ext^*_{\th\!
G_\Vb}(\th\!\Lc_{\ib}^\vee,\th\!\Lc_{\ib'}^\vee).$$ The Yoneda composition is a
homogeneous $\th \Sb_\Vb$-bilinear map of degree zero
$$\th \Zb_{\Lbb,\Vb,\ib,\ib'}\times\th
\Zb_{\Lbb,\Vb,\ib',\ib''}\to\th\Zb_{\Lbb,\Vb,\ib,\ib''},\quad\ib,\ib',\ib''\in\th\!I^\nu.$$ 
It equips the $\k$-vector space
$$\th\Zb_{\Lbb,\Vb}=\bigoplus_{\ib,\ib'\in
\th\! I^\nu}\th\Zb_{\Lbb,\Vb,\ib,\ib'}$$ with the structure of
a unital associative graded $\th\Sb_\Vb$-algebra. For
$\ib\in \th\!I^\nu$ we have the graded
$\th\Sb_\Vb$-modules
$$\th\Fb_{\Lbb,\Vb,\ib}=\Ext^*_{\th\! G_\Vb}(\th\!\Lc_{\Vb,\ib}^\vee,\Dc), \quad
\th\Fb_{\Lbb,\Vb}=\bigoplus_{\ib\in \th\! I^\nu}\th\Fb_{\Lbb,\Vb,\ib}.$$ For each $\ib,\ib'$ in $\th\!I^\nu$ the Yoneda
product gives a graded $\th\Sb_\Vb$-bilinear map
$\th\Zb_{\Lbb,\Vb,\ib,\ib'}\times\th\Fb_{\Lbb,\Vb,\ib'}\to\th\Fb_{\Lbb,\Vb,\ib}.$
This yields a left graded representation of $\th\Zb_{\Lbb,\Vb}$ on
$\th\Fb_{\Lbb,\Vb}$. Our first goal is to compute the graded algebra
$\th\Zb_{\Lbb,\Vb}$ and the graded representation $\th\Fb_{\Lbb,\Vb}$. 
For $\ib\in\th\! I^\nu$ let
$1_{\Lbb,\Vb,\ib}$ be the identity of $\th\! \Lc_{\ib}$, 
regarded as an element of
$\th\Zb_{\Lbb,\Vb,\ib,\ib}$. 
The elements $1_{\Lbb,\Vb,\ib}$ form a complete set
of orthogonal idempotents of $\th\Zb_{\Lbb,\Vb}$ such that
$$\th\Zb_{\Lbb,\Vb,\ib,\ib'}=
1_{\Lbb,\Vb,\ib}\,\th\Zb_{\Lbb,\Vb}\, 1_{\Lbb,\Vb,\ib'},\quad
\th\Fb_{\Lbb,\Vb,\ib}=1_{\Lbb,\Vb,\ib}\,\th\Fb_{\Lbb,\Vb}.$$

\vskip3mm

\subhead 2.7.~Remark\endsubhead
Fix a pair $\yb$ in $\th Y^\nu$.
Let $\ib$ be the sequence of $\th\! I^\nu$ obtained by expanding $\yb$. We have
an isomorphism of complexes in the derived category
$$\th\!\Lc_{\ib}^\delta=\bigoplus_{w\in\Sen_\bb}\th\!\Lc^\delta_\yb[\ell_\bb-2\ell(w)].$$ 
Here $\bb=(b_1,\dots, b_m)$ is a sequence such that  the multiplicity of $\yb$
is 
$$\theta(\bb)\bb:=(b_m,\dots,b_2,b_1,b_1,b_2,\dots b_m).$$
We'll abbreviate $\th\!\Lc_{\ib}^\delta=\la\bb\ra!\,\th\!\Lc^\delta_\yb.$

\vskip3mm
\subhead 2.8.~Shift of the grading\endsubhead
Let  $\l$, $\nu$, $\Lbb$, $\Vb$ be as above.
We define a new grading on $\th\Zb_{\Lbb,\Vb}$ and
$\th\Fb_{\Lbb,\Vb}$ by
$$\gathered
\th\Zb_{\Lbb,\Vb,\ib,\ib'}^\delta= \Ext^*_{\th
G_\Vb}(\th\!\Lc_{\ib}^\delta,\th\!\Lc_{\ib'}^\delta)=
\th\Zb_{\Lbb,\Vb,\ib,\ib'}[d_{\l,\ib}-d_{\l,\ib'}],\vspace{2mm}
\th\Zb_{\Lbb,\Vb}^\delta= 
\bigoplus_{\ib,\ib'\in
\th\! I^\nu}\th\Zb^\delta_{\Lbb,\Vb,\ib,\ib'},\vspace{1mm}
\th\!\Lc_{\Vb}=\bigoplus_{\ib\in
\th \!I^\nu}\th\!\Lc_{\ib},\quad
\th\!\Lc_{\Vb}^\delta=\bigoplus_{\ib\in\th \!I^\nu}\th\!\Lc_{\ib}^\delta.
\endgathered$$
We set also $\th\!I^0=\{\emptyset\}$, 
$\th\!\Lc^\delta_\emptyset=\kb$,
and
$\th\Zb_{\Lbb,\{0\}}^\delta=\kb$ as a graded $\kb$-algebra.
Here $\kb$ is regarded as the constant sheaf over $\{0\}$.

\vskip2cm

\head 3.~The convolution algebra
\endhead

Fix a quiver $\Gamma$ with set of vertices $I$ and set of arrows $H$.
Fix an involution $\theta$ on $\Gamma$. 
Assume that $\Gamma$ has no 1-loops and that
$\theta$ has no fixed points. Fix a
dimension vector $\nu\neq 0$ in $\th\NN I$ and a
dimension vector $\l$ in $\NN I$. 
Fix an object $(\Vb,\varpi)$ in $\th\Vcb_\nu$ 
and an object $\Lbb$ in $\Vcb_\l$.
For each sequences $\ib$, $\ib'$ in $\th\! I^\nu$ we
set
$$\th\! Z_{\Lbb,\Vb,\ib,\ib'}=\th\!\widetilde F_{\Lbb,\Vb,\ib}
\times_{\th\! E_{\Lbb,\Vb}}\th\!\widetilde F_{\Lbb,\Vb,\ib'},
\quad\th\! Z_{\Lbb,\Vb}=\coprod_{\ib,\ib'\in
\th\! I^\nu}\th\! Z_{\Lbb,\Vb,\ib,\ib'}.$$
the reduced fiber product relative to the maps $\th\!\pi_{\Lbb,\ib}$,
$\th\!\pi_{\Lbb,\ib'}$.
Next we set $$\th\!\Zc_{\Lbb,\Vb}=\bigoplus_{\ib,\ib'\in
\th\! I^\nu}\th\!\Zc_{\Lbb,\Vb,\ib,\ib'},\quad
\th\Fc_{\Lbb,\Vb}=\bigoplus_{\ib\in \th\! I^\nu}\th\Fc_{\Lbb,\Vb,\ib},$$ where
$$\th\!\Zc_{\Lbb,\Vb,\ib,\ib'}=H_{*}^{\th\! G_\Vb}
(\th\! Z_{\Lbb,\Vb,\ib,\ib'},\k),
\quad
\th\Fc_{\Lbb,\Vb,\ib}=H_{*}^{\th\! G_\Vb}(\th\widetilde F_{\Lbb,\Vb,\ib},\k).$$
We have
$$\th\Fb_{\Lbb,\Vb,\ib}=\Ext^*_{\th\! G_\Vb}(\kb,\th\!\Lc_{\ib})=
H^*_{\th\! G_\Vb}(\th\! E_\Vb,\th\!\Lc_{\ib})=
H^*_{\th\! G_\Vb}(\th\widetilde F_{\Lbb,\Vb,\ib},\kb).$$
We have also
$$H^*_{\th\! G_\Vb}(\th\widetilde F_{\Lbb,\Vb,\ib},\k)=
H^{*}_{\th\! G_\Vb}(\th\widetilde F_{\Lbb,\Vb,\ib},\Dc)[-2d_{\l,\ib}]=
\th\Fc_{\Lbb,\Vb,\ib}[-2d_{\l,\ib}].
\leqno(3.1)$$ This yields a graded
$\th\Sb_\Vb$-module isomorphism
$$\th\Fb_{\Lbb,\Vb,\ib}=\th\Fc_{\Lbb,\Vb,\ib}[-2d_{\l,\ib}].\leqno(3.2)$$
We equip the
$\th\Sb_\Vb$-module $\th\!\Zc_{\Lbb,\Vb}$ with the
convolution product relative to the closed embedding of $\th
Z_{\Lbb,\Vb}$ into $\th\widetilde F_{\Lbb,\Vb}
\times\th\widetilde F_{\Lbb,\Vb}$. See
\cite{CG, sec.~8.6} for details. We obtain an associative 
graded $\th \Sb_\Vb$-algebra $\th\!\Zc_{\Lbb,\Vb}$ with 1 which acts on the graded $\th
\Sb_\Vb$-module $\th\Fc_{\Lbb,\Vb}$. The unit  is the fundamental
class of the closed subvariety $\th\! Z^e_{\Lbb,\Vb}$ of $\th\! Z_{\Lbb,\Vb}$. 
See Section 4.6 below for the notation.

\proclaim{3.1.~Proposition} (a) The left $\th\!\Zc_{\Lbb,\Vb}$-module 
$\th\Fc_{\Lbb,\Vb}$ is
faithful.

\vskip1mm

(b) There is a canonical $\th\Sb_\Vb$-algebra isomorphism
$\th\Zb_{\Lbb,\Vb}=\th\!\Zc_{\Lbb,\Vb}$ such that 
$(3.2)$ identifies the $\th\Zb_{\Lbb,\Vb}$-action on
$\th\Fb_{\Lbb,\Vb}$ and the $\th\!\Zc_{\Lbb,\Vb}$-action on $\th\Fc_{\Lbb,\Vb}$.
\endproclaim

\noindent{\sl Proof :}
This is standard material, see e.g., \cite{VV}.
Let us give one proof of $(a)$. It is a consequence
of the following general fact.
Let $G$ be a linear algebraic group over $\CC$ and let $M$ be a smooth 
quasi-projective $G$-variety over $\CC$. Let $T\subset G$ be a maximal torus.
Let $\Qb$ be the fraction field of $\Sb=\Sb_T$.
Let $Z\subset M\times M$ be a closed $G$-stable subset
(for the diagonal action on $M\times M$) such that $p_{1,3}$
restricts to a proper map
$$p_{1,2}^{-1}(Z)\cap p_{2,3}^{-1}(Z)\to Z,$$
where $p_{i,j}:M\times M\times M\to M\times M$ is the projection 
along the factor not named. The convolution product equips
$H_*^G(Z,\kb)$ with a $\Sb_G$-algebra structure and $H_*^G(M,\kb)$ with
a $H_*^G(Z,\kb)$-module structure, see e.g., \cite{CG}.
Assume now that the $T$-spaces $M$, $Z$ are {\it equivariantly formal},
see e.g., \cite{GKM, Sec.~1.2},
and assume that we have the following equality of $T$-fixed points subsets
$$Z^T=M^T\times M^T.\leqno(3.3)$$
Consider the following commutative diagram 
of algebra homomorphisms 
$$\xymatrix{
H_*^T(Z,\kb)\otimes_\Sb\Qb\ar[r]^-c&\End_\Sb(H_*^T(M,\kb))\otimes_\Sb\Qb\cr
H_*^T(Z,\kb)\ar[u]^b\ar[r]&\End_\Sb(H_*^T(M,\kb))\ar[u]\cr
H_*^G(Z,\kb)\ar[u]^a\ar[r]&\End_{\Sb_G}(H_*^G(M,\kb)).\ar[u]
}$$
The map $c$ is invertible by (3.3) and
the localization theorem in equivariant homology.
The map $b$ is injective because $Z$ is equivariantly formal.
The map $a$ is injective, compare Section 4.10 below.
Thus the lower map is injective, i.e., the
$H_*^G(Z,\kb)$-module $H_*^G(M,\kb)$ is faithful.

\qed

\vskip2cm

\head 4.~ The polynomial representation
of the graded algebra $\th\Zb_{\Lbb,\Vb}^\delta$\endhead

Fix a quiver $\Gamma$ with set of vertices $I$ and set of arrows $H$.
Fix an involution $\theta$ on $\Gamma$. 
Assume that $\Gamma$ has no 1-loops and that
$\theta$ has no fixed points. Fix a
dimension vector $\nu\neq 0$ in $\th\NN I$ and a
dimension vector $\l$ in $\NN I$. Set $|\nu|=2m$.
Fix an object $(\Vb,\varpi)$ in $\th\Vcb_\nu$ 
and an object $\Lbb$ in $\Vcb_\l$.
The main result of this section is Theorem 4.17 which yields an explicit
faithful representation of the graded $\kb$-algebra
$\th\Zb_{\Lbb,\Vb}^\delta$.

\vskip3mm

\subhead 4.1.~Notations\endsubhead  Let $G=O(\Vb,\varpi)$ be the
orthogonal group, and $F=F(\Vb,\varpi)$ be the isotropic flag
manifold. We can regard $F$  as the (non connected) flag manifold
of the (non connected) group $G$. 
Next, the group $\th\! G_\Vb$ is canonically identified with a 
L\'evi subgroup of $G$, i.e., with the subgroup of elements which 
preserve the decomposition $\Vb=\bigoplus_i\Vb_i$. Then 
$\th\! F_\Vb$ is canonically identified
with the closed subvariety of $F$ consisting
of all flags which are fixed under the action of the center of $\th\!
G_\Vb$. Fix once for all a maximal torus $T$ of $\th\! G_\Vb$. Let
$W_\Vb$ and $W$ be the Weyl groups of the pairs $(\th\!
G_\Vb,T)$ and $(G,T)$. The canonical inclusion
$\th\! G_\Vb\subset G$ yields a canonical inclusion
$W_\Vb\subset W$. 

\vskip3mm

\subhead 4.2.~The root systems\endsubhead 
Fix once for all a $T$-fixed flag
$\phi_\Vb$ in $\th\! F_\Vb$. 
We fix once for all
one-dimensional $T$-submodules $\Db_{1-m},\dots,\Db_{m-1},\Db_{m}$ of $\Vb$
such that
$$\phi_\Vb=(\Vb^l),\quad 
\Vb^l=\Db_{l+1}\oplus\dots\oplus\Db_{m-1}\oplus\Db_{m}
.$$ 
Let $\chi_l\in\ten^*$ be the weight of  $\Db_l$.
Note that $\Db_l\simeq \Vb^{l-1}/\Vb^{l}$ and that
the bilinear form $\varpi$ yields a non-degenerate pairing
$(\Vb^{l-1}/\Vb^{l})\times(\Vb^{-l}/\Vb^{1-l})\to\CC$,
because $(\Vb^l)^\perp=\Vb^{-l}$. 
Thus we have $$\chi_{1-l}=-\chi_l.$$
Let $B$ be the stabilizer of the flag
$\phi_\Vb$ in $G$.  Let $\Delta$ be the set of roots of $(G,T)$ and
let $\Delta^+$ be the subset of positive roots relative to
the Borel subgroup $B$. We abbreviate $\Delta^-=-\Delta^+$. Let
$\Pi$ be the set of simple roots in $\Delta^+$. 
We have
$$\gathered
\Delta^+=\{\chi_k\pm \chi_l;\,1\leqslant l<k\leqslant m\},\vspace{2mm}
\Pi=\{\chi_{l+1}-\chi_{l},\,\chi_{2}+\chi_1;\,l=1,2,\dots,m-1\}.\endgathered$$
Let $\leqslant$ and $\ell$ denote the Bruhat order and the length
function on $W$.
Note that $W$ is an extended Weyl group of type
$\D_m$. In particular we have $$\ell(w)=0\iff w=e,\eps_1,$$
where $\eps_1$is as below, and the set $S$ of simple reflections is given by
$$S=\{s_0,s_1,\dots,s_{m-1}\},$$ with 
$s_k$, $k=0,1,\dots, m-1$ the reflection
with respect to 
$$\a_0=\chi_2+\chi_1,\quad\a_1=\chi_2-\chi_1,\quad \dots
\quad\a_{m-1}=\chi_{m}-\chi_{m-1}.$$
Note that $u(\Delta^+)=\Delta^+$ if $\ell(u)=0$.
Next, let $\th\!\Delta_\Vb\subset\Delta$
be the set of roots of $(\th\! G_\Vb,T)$.
Note that $\th\!G_\Vb$ is a product of general linear groups
(this is due to the fact that $\theta$ has no fixed points).
Indeed, we can (and we will) assume that
$$\th\!\Delta_\Vb\subset\{\chi_l-\chi_k;\,l\neq k,\,l,k=1,2,\dots,m\}.$$
More precisely, given a subset $J\subset I$ such that 
$I=J\sqcup\theta(J)$ it is enough
to choose the flag $\phi_\Vb$ such that $\Vb^0=\bigoplus_{j\in J}\Vb_j$.
Finally, let
$\th\!\Delta^+_\Vb$ be the subset of positive roots relative to
the Borel subgroup $\th\! B_\Vb=B\cap\th\! G_\Vb$.
We have $$\th\!\Delta_\Vb^+=\Delta^+\cap\th\!\Delta_\Vb.$$

\vskip3mm

\subhead 4.3.~The wreath product\endsubhead 
Let $\Sen_m$ be the symmetric group, and $\ZZ_2=\{-1,1\}$.
Consider the wreath product 
$W_m=\Sen_m\wr\ZZ_2.$ For $l=1,2,\dots m$ let $\eps_l\in(\ZZ_2)^m$
be $-1$ placed at the $l$-th position. We'll regard $\eps_l$ as in element of $W_m$ in the obvious way.
There is a
unique action of $W_m$ on the set $\{1-m,\dots,m-1,m\}$ such that
$\Sen_m$ permutes $1,2,\dots m$ and such that
$\eps_l$ fixes $k$ if $k\neq l, 1-l$ and switches $l$ and
$1-l$. 
The group $W_m$ acts also on $\th\!I^\nu$. Indeed,
view a sequence $\ib$ as the map 
$$\{1-m,\dots, m-1,m\}\to I,\quad l\mapsto i_l.$$ 
Then we set $w(\ib)=\ib\circ w^{-1}$ for $w\in W_m$.

\vskip3mm

\subhead 4.4.~The $W$-action on the set of $T$-fixed flags \endsubhead 
The sets $F^T$ and $(\th\! F_\Vb)^{T}$
consisting of the flags which are fixed by the $T$-action are equal.
The group $W$ acts freely transitively on both.
We'll write $e$ for the unit in $W$.
Put $$\phi_{\Vb,w}=w(\phi_\Vb),\quad\forall w\in W. $$
Thus we have $F^T=\{\phi_{\Vb,w};w\in W\}$. 
There is a unique group
isomorphism $W=W_m$ such that
$$\phi_{\Vb,w}=(\Vb^l_w),\quad 
\Vb^l_w=\Db_{w(l+1)}\oplus\dots\oplus\Db_{w(m-1)}\oplus\Db_{w(m)}
.$$ 
We'll use this identification whenever it is convenient without recalling
it explicitly.
We set also $$w(\chi_l)=\chi_{w(l)},\quad \forall w,l.$$
Let $\th\! B_{\Vb,w}$ be the stabilizer of the flag
$\phi_{\Vb,w}$ under the $\th\! G_\Vb$-action. It is the Borel subgroup
of $\th\! G_\Vb$ containing $T$ associated with the set of positive roots
$$w(\Delta^+)\cap\th\!\Delta_\Vb.$$ 
Let $\th\! N_{\Vb,w}$ be the unipotent
radical of $\th\! B_{\Vb,w}$.
Finally, let $\ib_w$ be the unique sequence in
$\th\!I^\nu$ such that $\phi_{\Vb,w}$ lies in $\th\! F_{\Vb,\ib_w}$. 
Write
$$\ib_e=(i_{1-m},\dots,i_{m-1},i_{m}).\leqno(4.1)$$ 
Since $\phi_\Vb$ is a flag of type $\ib_e$, we have
$$\Db_l\subset\Vb_{i_l},\quad w^{-1}(\ib_e)=\ib_{w}=(i_{w(1-m)},\dots,i_{w(m-1)},i_{w(m)}).$$
Let
$W_\nu$ be the image of the group
$W_\Vb$ by the isomorphism $W\to W_m$. 
It is the parabolic subgroup given by
$$W_\nu=\{w\in W_m;w(\ib_e)=\ib_e\}.$$
Note that the choices made in Section 4.2 imply that
$$W_\nu\subset\Sen_m.\leqno(4.2)$$
There is a bijection
$$W_\nu\setminus W_m\to \th\!I^\nu,\quad W_\nu w\mapsto
\ib_w.$$
For each  $\ib$ in $\th\! I^\nu$ we have $$(\th\widetilde
F_{\Vb,\ib})^{T}\simeq (\th
F_{\Vb,\ib})^{T}=\{\phi_{\Vb,w};w\in W_{\ib}\},
\quad W_{\ib}=\{w\in W;\ib_w=\ib\}.$$
 We'll abbreviate
$$\th\! F_{\Vb,w}=\th\! F_{\Vb,\ib_w},\quad W_w=W_{\ib_w},\quad
\th\pi_{\Lbb,w}=\th\pi_{\Lbb,\ib_w}.$$ 
We'll also omit
the symbol $w$ if $w=e$. For instance we write $\th\! B_\Vb=\th
B_{\Vb,e}$ and $\th\! N_\Vb=\th\! N_{\Vb,e}.$ Note that
$W_w=W_\Vb\, w$ and that we have an isomorphism of
$\th\!G_\Vb$-varieties
$$\th\! G_\Vb/\th\! B_{\Vb,w}\to \th\! F_{\Vb,w},\quad g\mapsto g\phi_{\Vb,w}.$$

\vskip3mm

\subhead 4.5.~The stratification of $\th\! F_\Vb\times\th\! F_\Vb$\endsubhead
The group $G$ acts diagonally on $F\times
F$. The action of the subgroup $\th\! G_\Vb$ preserves the subset $\th\!
F_\Vb\times\th\! F_\Vb$. For  $w\in W$ let 
$\th\!O^{w}_\Vb$ be the set of
all pairs of flags in $\th\! F_\Vb\times \th\! F_\Vb$ which are in
relative position $w$. More precisely, we write
$$\th\!O^{w}_\Vb=(\th\! F_\Vb\times\th\! F_\Vb)\cap (G\phi_{\Vb,e,w}),\quad
\phi_{\Vb,x,y}=(\phi_{\Vb,x},\phi_{\Vb,y}),\quad\forall
x,y\in W.$$
Let $\th\bar O_\Vb^{w}$ be the Zariski
closure of $\th\!O_\Vb^{w}$. 
For any $w,x,y$ in $W$ we write also
$$\th\!O^{w}_{\Vb,x,y}=\th\!O^{w}_\Vb\cap(
\th\! F_{\Vb,x}\times\th\! F_{\Vb,y}),\quad \th\bar
O^{w}_{\Vb,x,y}=\th\bar O^{w}_\Vb\cap( \th\! F_{\Vb,x}\times\th
\!F_{\Vb,y}).$$ 
We define $\th P_{\Vb,w,ws}$, $s\in S$, as the smallest parabolic
subgroup of $\th\! G_\Vb$ containing $\th\! B_{\Vb,w}$ and $\th
B_{\Vb,ws}$.

\proclaim{4.6. Lemma} Let $w,x,y,s,u\in W$.

\vskip1mm

(a) The set of $T$-fixed elements in $\th\!O^x_\Vb$ is 
$\{\phi_{\Vb,w,wx};\,w\in W\}$.
\vskip1mm

(b) Assume that $\ell(u)=0$.
We have $\th\bar O_\Vb^u=\th\!O_\Vb^u$. It is a smooth 
$\th\! G_\Vb$-variety isomorphic to $\th\! F_\Vb$.
We have $\th\!O_{\Vb,x,y}^u=\emptyset$ unless $y=xu$.
\vskip1mm

(c) Assume that $\ell(s)=1$. 
Set $s=s'u$ with $s'\in S$ and $\ell(u)=0$.
We have $\th\bar O_\Vb^{s}=\th\!O_\Vb^{s}\cup\th\!O_\Vb^{u}$. 
It is a smooth variety. 
We have
$\th\bar O^{s}_{\Vb,x,y}=\emptyset$ if $y\neq xs, xu$. 

\vskip1mm

\itemitem\itemitem{$\bullet$} If $xs\notin W_{xu}$ then 
$$\th\! F_{\Vb,xs}\neq\th\! F_{\Vb,xu},\quad
\th\! B_{\Vb,xs}=\th\! B_{\Vb,x},\quad \th\!O^u_{\Vb,x,xs}=\th
O^s_{\Vb,x,xu}=\emptyset.$$

\vskip1mm

\itemitem\itemitem{$\bullet$} If $xs\in W_{xu}$ then 
$$\gathered
\th\! F_{\Vb,xs}=\th\! F_{\Vb,xu},\quad \th\! B_{\Vb,xs}\neq
\th\! B_{\Vb,x},
\vspace{2mm}
\th\! G_\Vb\times_{\th\! B_{\Vb,x}}(\th P_{\Vb,x,xs}/\th\! B_{\Vb,x})
{\buildrel\sim\over\to}
\th\bar O^s_{\Vb,x,xu}=\th\bar O^s_{\Vb,x,xs},\ 
(g,h)\mapsto(g\phi_{\Vb,x},gh\phi_{\Vb,xu}).
\endgathered$$

\endproclaim

\noindent{\sl Proof :}
The proof is standard and is left to the reader.
Note that $\th\! B_{\Vb,x}=\th\! B_{\Vb,xu}$ and
$\th\! B_{\Vb,xs}=\th\! B_{\Vb,xs'}$ because $\ell(u)=0$.
Note also that 
$$\th\! B_{\Vb,x}=\th\! B_{\Vb,xs}\iff x(\a)\notin\th\!\Delta_\Vb
\iff xs'\in W_x\iff xs\in W_{xu},$$
where $\a$ is the simple root associated with $s'$.

\qed

\vskip3mm

For a future use let us introduce the following
notation. Let $q$ be the obvious projection $\th\! Z_{\Lbb,\Vb}\to \th\! F_\Vb\times \th\! F_\Vb,$
and, for each $x\in W$, let
$\th\! Z_{\Lbb,\Vb}^{x}$ 
be the Zariski closure in $\th\! Z_{\Lbb,\Vb}$ of
the locally closed subset $q^{-1}(\th\!O^{x}_\Vb)$.

\vskip3mm

\subhead 4.7.~Euler classes in $\Sb$\endsubhead
Consider the graded $\kb$-algebra  $\Sb=\Sb_T$. The weights
$\chi_1,\chi_2,\dots\chi_{m}$
are algebraically independent generators of  $\Sb$ and they are homogeneous of degree 2.
The reflection representation  on $\ten$ yields a $W$-action on $\Sb$.
Recall that we have
$$w(\chi_l)=\chi_{w(l)},\quad\forall l,w.$$
Now, let $M$ be a finite dimensional
representation of $\ten$ and fix a linear form $\l\in\ten^*$.
Let $M[\l]\subset M$ be the weight subspace associated with $\l$.
The {\it character} of $M$ is the linear form 
$\ch(M)=\sum_\l \dim (M[\l])\,\l$.
Let $\Eu(M)$ be the
determinant of $M$, viewed as  an element of degree $2\dim(M)$ of
$\Sb$. We'll call $\Eu(M)$ the {\it Euler class of  $M$}.
If $M$ is a finite dimensional representation of $T$
let $\Eu(M)$ be the Euler class of the
differential of $M$, a module over $\ten$.
Now, assume that $X$ is a quasi-projective $T$-variety and that
$x\in X^{T}$ is a smooth point of $X$. The cotangent space
$T^*_xX$ at $x$ is equipped with a natural representation of
$T$. We'll abbreviate $\Eu(X,x)=\Eu(T_x^*X)$.
We'll be particularly interested in the following elements 
$$\Lambda_{w}=\Eu(\th\widetilde F_{\Lbb,\Vb},\phi_{\Vb,w}),\quad
\Lambda^{x}_{w,w'}=\Eu(\th\! Z_{\Lbb,\Vb}^{x},\phi_{\Vb,w,w'})^{-1}$$ where
$\ell(x)=0,1$. Note that $\Lambda_w$ lies in $\Sb$ and has the degree
$2d_{\l,w}$.

\vskip3mm

\subhead{4.8.~ Description of the $\th\! G_\Vb$-varieties
$\th\widetilde F_{\Lbb,\Vb,w}$}\endsubhead
Let $\th\!\gen_\Vb$, $\ten$,
$\th\!\nen_{\Vb,w}$, $w\in W$, be the Lie algebras of $\th\! G_\Vb$,
$T$, $\th\! N_{\Vb,w}$ respectively.
Consider the flag 
$$\phi_{\Vb,w}=(\Vb=\Vb^{-m}_w\supset\dots\supset\Vb^{m-1}_w\supset\Vb^m_w=0).$$ 
The $\th\! G_\Vb$-action on
$\th\! E_{\Lbb,\Vb}$ yields a representation of $\th\! B_{\Vb,w}$ on the space
$$\th\!\een_{\Lbb,\Vb,w}=\{(x,y)\in \th\! E_{\Lbb,\Vb};\,x(\Vb^l_w)\subset\Vb^l_w,\,y(\Lbb)\subset\Vb^m_w\}.$$
There is an isomorphism of $\th\! G_\Vb$-varieties
$$\th\! G_\Vb\times_{\th\! B_{\Vb,w}}\th\!\een_{\Lbb,\Vb,w}\to
\th\widetilde F_{\Lbb,\Vb,w},
\quad(g,x,y)\mapsto (g\phi_{\Vb,w},gx,gy).$$ Under this isomorphism the
map $\th \pi_{\Lbb,w}$ is identified with the map
$$\th\! G_\Vb\times_{\th\! B_{\Vb,w}}\th \een_{\Lbb,\Vb,w}\to\th\! E_{\Lbb,\Vb},\quad
(g,x,y)\mapsto (gx,gy).$$

\vskip3mm

\subhead{4.9. ~Character formulas}\endsubhead 
In this section we gather some character formula for
a later use. 
For $w,w'\in W$ we write
$$\gathered
\th\!\een_{\Lbb,\Vb,w,w'}=\th\!\een_{\Lbb,\Vb,w}\cap\th\!\een_{\Lbb,\Vb,w'},
\vspace{2mm}
\th\den_{\Lbb,\Vb,w,w'}=\th\!\een_{\Lbb,\Vb,w}/\th\!\een_{\Lbb,\Vb,w,w'},
\vspace{2mm}
\th\!\nen_{\Vb,w,w'}=\th\!\nen_{\Vb,w}\cap\th\!\nen_{\Vb,w'},
\vspace{2mm}
\th\!\men_{\Vb,w,w'}=\th\!\nen_{\Vb,w}/\th\!\nen_{\Vb,w,w'}.
\endgathered$$
We have the following $T$-module isomorphisms
$$\gathered
\th\!\nen_{\Vb,w}=\th\!\nen_{\Vb,w,w'}\oplus\th\!\men_{\Vb,w,w'},
\vspace{2mm}
\th\!\een_{\Lbb,\Vb,w}=\th\!\een_{\Lbb,\Vb,w,w'}\oplus\th\den_{\Lbb,\Vb,w,w'},
\vspace{2mm}
\th\!\men_{\Vb,w,w'}=(\th\!\men_{\Vb,w',w})^*.
\endgathered$$
Write 
$\th\!\een_{\Vb,w}=\th\!\een_{\{0\},\Vb,w}.$
As $T$-modules we have
$$\gathered
\th\!\nen_{\Vb,w}=\bigoplus_\a\gen[\a],\quad
\a\in w(\Delta^+)\cap\th\!\Delta_\Vb,\vspace{2mm}
\th\!\een_{\Vb,w}=\bigoplus_\a\th\!E_\Vb[\a],\quad
\a\in w(\Delta^+).
\endgathered\leqno{(4.3)}$$
Recall that
$\Vb=\bigoplus_l\Db_l$ as $I$-graded $T$-modules, where
$l=1-m,\dots,m-1,m$. 
We'll use the notation in (4.1). Thus  $i_l$, $\chi_l$ are the dimension vector and the character of $\Db_l$.
Note that 
$$\gathered
h_{i_k,i_l}=h_{i_{1-l},i_{1-k}},\quad
\chi_l=-\chi_{1-l},\vspace{2mm}
\nu_{\theta(i)}=\nu_i,\quad
\nu=\sum_i\nu_i\,i=\sum_{l=1}^m(i_l+i_{1-l}).
\endgathered$$
Set $H^0=\{h\in H;\,h'=\theta(h'')\}$, $H^1=H\setminus H^0$, 
and $\l=\sum_i\l_i\,i$.
Note that $H^0=\{h\in H;\,h=\theta(h)\}$.
Decomposing a tuple $x\in\th\!E_\Vb$ as the sum of
$(x_h)_{h\in H^1}$ and $(x_h)_{h\in H^0}$ we get the following formula 
$$\dim(\th\! E_\Vb)=
\sum_{h\in H^1}\nu_{h'}\nu_{h''}/2+\sum_{h\in H^0}\nu_{h'}(\nu_{h'}-1)/2.$$
Next, the decomposition (4.3) yields the following formula
$$\ch(\th\!\een_{\Vb,w})=
\sum_{\chi_l-\chi_k\in w(\Delta^+)}h_{i_k,i_l}(\chi_l-\chi_k).$$
Here the sum runs over $\a\in w(\Delta^+)$, and for each $\a$ we choose
one pair $(l,k)$ such that $\a=\chi_l-\chi_k$.
In a similar way we have also
$$\gathered
\ch(\th\! E_\Vb)=\sum_{\chi_l-\chi_k\in\Delta}h_{i_k,i_l}(\chi_l-\chi_k),
\vspace{2mm} 
\ch(L_{\Lbb,\Vb})=\sum_l\l_{i_l}\chi_l.
\endgathered$$
Here the first sum runs over $\Delta$.
Since $\Vb^0=\bigoplus_{l\geqslant 1}\Db_l$
we have also
$$\ch(\th\!\een_{\Lbb,\Vb,w})=\sum_{\chi_l-\chi_k}h_{i_k,i_l}(\chi_l-\chi_k)+
\sum_l\l_{i_l}\chi_l,\leqno(4.4)$$ 
where the first sum runs over all roots in $w(\Delta^+)$ 
and the second one over all $l$ in
$\{w(1),w(2),\dots, w(m)\}$.
Note that (4.4) can be rewriten in the following way
$$\ch(\th\!\een_{\Lbb,\Vb,w})=
\sum_{\chi_l-\chi_k\in\Delta^+}h_{i_{w(k)},i_{w(l)}}w(\chi_l-\chi_k)+
\sum_{1\leqslant l\leqslant m}\l_{i_{w(l)}}w(\chi_l).$$ 
By (4.3) the Euler class
$\Eu(\th\!\nen_{\Vb,w})$ is the product of all roots in $\th\!\Delta_\Vb\cap
w(\Delta^+)$. Therefore, for $s\in S$ the following formulas hold

\vskip1mm

\itemitem{$\bullet$} either $ws\notin W_w$ and we have
$$\gathered
\Eu(\th\!\nen_{\Vb,ws})=\Eu(\th\!\nen_{\Vb,w}),\vspace{2mm}
\Eu(\th\!\men_{\Vb,w,ws})=\Eu(\th\!\men_{\Vb,ws,w})=0,\endgathered$$

\itemitem{$\bullet$} 
or $ws\in W_w$ and we have
$$\gathered\Eu(\th\!\nen_{\Vb,ws})=-\Eu(\th\!\nen_{\Vb,w}),\vspace{2mm}
\Eu(\th\!\men_{\Vb,w,ws})=-\Eu(\th\!\men_{\Vb,ws,w})=w(\a),
\endgathered$$
where $\a$ is the simple root associated with $s$.

\vskip1mm

\noindent Finally, let $s=s_l$ with $l=0,1,\dots,m-1$.
Formula (4.4) yields the following. 
\vskip1mm
\itemitem{$\bullet$} 
We have 
$$\gathered
\eu(\th\!\den_{\Lbb,\Vb,w,w\eps_1})=w(\chi_1)^{\l_{i_{w(1)}}}.\endgathered$$
\vskip1mm
\itemitem{$\bullet$} 
If $l\neq 0$ we have
$$\gathered 
\eu(\th\!\den_{\Lbb,\Vb,w,ws_l})=w(\a_l)^{h_{i_{w(l)},i_{w(l+1)}}}.
\endgathered$$

\vskip1mm
\itemitem{$\bullet$} 
We have
$$\gathered
\eu(\th\!\den_{\Lbb,\Vb,w,ws_0})=
w(\chi_{1})^{\l_{i_{w(1)}}}w(\chi_{2})^{\l_{i_{w(2)}}}
w(\a_{0})^{h_{i_{w(0)},i_{w(2)}}}.\endgathered$$

\vskip3mm

\subhead{4.10.~Reduction to the torus}\endsubhead
The restriction of functions from $\th\!\gen_\Vb$ to $\ten$ gives an
isomorphism of graded $\k$-algebras
$$\th\Sb_\Vb=\kb[\chi_1,\chi_2,\dots,\chi_m]^{W_\nu}.$$
The group $\th\! G_\Vb$ is a product of several copies
of the general linear group. Hence it is connected with a simply connected
derived subgroup. It is a general fact that if $X$ is a
$\th\! G_\Vb$-variety then the $\Sb$-module
$H_*^{T}(X,\k)$ is equipped with a 
$\Sb$-skewlinear representation of the group $W_\Vb$ such that
the forgetful map gives a $\th\Sb_\Vb$-module
isomorphism
$$H_*^{\th\! G_\Vb}(X,\k)\to H_*^{T}(X,\k)^{W_\Vb},$$
see e.g., \cite{HS, thm.~2.10}.
We'll call this action on $H_*^{T}(X,\k)$ the canonical $W_\Vb$-action.

\vskip3mm

\subhead{4.11.~The $W$-action and the $\th\Sb_\Vb$-action
on $\th\Fb_{\Lbb,\Vb}$}\endsubhead Fix
a tuple $\ib$ in $\th I^\nu$ and an integer $l=1,2,\dots, m$. We
define $\Oc_{\Lbb,\Vb,\ib}(l)$ to be the $\th
G_\Vb$-equivariant line bundle over $\th\!\widetilde F_{\Lbb,\Vb,\ib}$
whose fiber at the triple $(x,y,\phi)$ with
$$\phi=(\Vb=\Vb^{-m}\supset\Vb^{1-m}\supset\cdots\supset\Vb^{m}=0)$$ 
is equal to $\Vb^{l-1}/\Vb^{l}$.
Assigning to a formal variable $x_\ib(l)$ of degree 2 the first
equivariant Chern class of $\Oc_{\Lbb,\Vb,\ib}(l)^{-1}$ we 
get a graded $\k$-algebra isomorphism
$$\k[x_\ib(1),x_\ib(2),\dots x_\ib(m)]=
H^*_{\th\!G_\Vb}(\th\!\widetilde F_{\Lbb,\Vb,\ib},\k).$$
So (3.1), (3.2) yield canonical isomorphisms of graded
$\k$-vector spaces
$$\k[x_\ib(1),x_\ib(2),\dots x_\ib(m)]
=\th\Fc_{\Lbb,\Vb,\ib}[-2d_{\l,\ib}]=
\th\Fb_{\Lbb,\Vb,\ib}.\leqno(4.5)$$ 
For a future use we set also
$$x_\ib(l)=-x_\ib(1-l),\quad l=1-m,2-m,\dots,0.$$
For $w\in W_m$ we set
$$w f(x_\ib(1),\dots,
x_\ib(m))=f(x_{w(\ib)}(w(1)),\dots,x_{w(\ib)}(w(m))).$$
This yields a $W_m$-action on $\th\Fb_{\Lbb,\Vb}$ such that
$w(\th\Fb_{\Lbb,\Vb,\ib})=\th\Fb_{\Lbb,\Vb,w(\ib)}$.

The multiplication of polynomials equip both 
$\th\Fb_{\Lbb,\Vb,\ib}$ and $\th\Fb_{\Lbb,\Vb}$ with an obvious
structure of graded $\k$-algebras. For $w\in W_{\ib}$ the pull-back by the 
inclusion $\{\phi_{\Vb,w}\}\subset \th\widetilde F_{\Lbb,\Vb,\ib}$ yields
a graded $\k$-algebra isomorphism
$$\th\Fb_{\Lbb,\Vb,\ib}\to\Sb,\quad
f(-x_\ib(1),\dots,-x_\ib(m))\mapsto
f(\chi_{w(1)},\dots\chi_{w(m)}).\leqno(4.6)$$ We'll abbreviate 
$$w(f)=f(\chi_{w(1)},\dots\chi_{w(m)}).$$ 
The isomorphism (4.6) is not
canonical, because it depends on the choice of $w$. 

Now, consider the
canonical $\th\Sb_\Vb$-action on $\th\Fb_{\Lbb,\Vb}$ coming from the
$\th\!G_\Vb$-equivariant cohomology. It can be regarded as a
$\th\Sb_\Vb$-action on $\bigoplus_\ib\k[x_\ib(1),x_\ib(2),\dots
x_\ib(m)]$ which is described in the following
way.  The composition of the obvious
projection $\th\Fb_{\Lbb,\Vb}\to\th\Fb_{\Lbb,\Vb,\ib}$ with the map (4.6)
identifies the
graded $\k$-algebra
of the $W_m$-invariant polynomials in the $x_\ib(l)$'s,
with 
$$\th\Sb_\Vb=\Sb^{W_\nu}=\kb[\chi_1,\chi_2,\dots,\chi_m]^{W_\nu}.$$
This isomorphism
does not depend on the choice of $\ib,w$.
The 
$\th\Sb_\Vb$-action on $\th\Fb_{\Lbb,\Vb}$ 
is the composition of this isomorphism
and of the multiplication by $W_m$-invariant polynomials.

\vskip3mm

\subhead 4.12.~Localization and the convolution product\endsubhead
Let $\Qb$ be the fraction field of $\Sb$.
Write
$$\gathered
\th\Fc'_{\Lbb,\Vb}=H_*^T(\th\widetilde F_{\Lbb,\Vb},\k),
\vspace{2mm}
\th\Fc''_{\Lbb,\Vb}=\th\Fc'_{\Lbb,\Vb}\otimes_{\Sb}\Qb,
\vspace{2mm}
\th\!\Zc'_{\Lbb,\Vb}=H_*^T(\th\! Z_{\Lbb,\Vb},\k),
\vspace{2mm}
\th\!\Zc''_{\Lbb,\Vb}=\th\!\Zc'_{\Lbb,\Vb}\otimes_{\Sb}\Qb.
\endgathered$$
Let
$\psi_{w}$ be the fundamental class of the singleton $\{\phi_{\Vb,w}\}$
in $\th\Fc'_{\Lbb,\Vb}$, and let
$\psi_{w,w'}$ be the fundamental class of $\{\phi_{\Vb,w,w'}\}$ in 
$\th\!\Zc'_{\Lbb,\Vb}.$
Let
$\psi_w$, $\psi_{w,w'}$ denote also the corresponding  elements in the 
$\Qb$-vector spaces $\th\Fc''_{\Lbb,\Vb}$,
$\th\!\Zc''_{\Lbb,\Vb}.$
Now, we consider the convolution products
$$\th\!\Zc'_{\Lbb,\Vb}\times \th\!\Zc'_{\Lbb,\Vb}\to \th\!\Zc'_{\Lbb,\Vb},\quad
\th\!\Zc'_{\Lbb,\Vb}\times \th\Fc'_{\Lbb,\Vb}\to\th\Fc'_{\Lbb,\Vb}$$
relative to the inclusion of
$\th\! Z_{\Lbb,\Vb}$ in the smooth scheme $\th\widetilde F_{\Lbb,\Vb}\times
\th\widetilde F_{\Lbb,\Vb}$.
Both may be denoted by the symbol $\star$. We'll use the notation in (4.1).

\proclaim{4.13.~Proposition} 
(a) The $\Sb$-modules
$\th\Fc'_{\Lbb,\Vb}$ and $\th\!\Zc'_{\Lbb,\Vb}$ are free.
The canonical $W_\Vb$-action on the $T$-equivariant homology spaces
$\th\Fc'_{\Lbb,\Vb}$ and $\th\!\Zc'_{\Lbb,\Vb}$ is given by $w(\psi_{x})=\psi_{wx}$ and $w(\psi_{x,y})=\psi_{wx,wy}$.
The  inclusions
$\th\!\Zc_{\Lbb,\Vb}\subset\th\!\Zc'_{\Lbb,\Vb}$ and 
$\th\Fc_{\Lbb,\Vb}\subset\th\Fc'_{\Lbb,\Vb}$
commute with the convolution products.

\vskip1mm

(b) The elements $\psi_w$, $\psi_{w,w'}$ yield $\Qb$-bases of
$\th\Fc''_{\Lbb,\Vb}$,
$\th\!\Zc''_{\Lbb,\Vb}$ respectively.
For each $\ib$ the map (4.5) 
yields an inclusion of
$\k[x_\ib(1),\dots, x_\ib(m)]$ into $\th\Fc''_{\Lbb,\Vb,\ib}$ such that
$$f(-x_\ib(1), \dots,- x_\ib(m))\mapsto\sum_{w\in W_{\ib}}w(f)\Lambda_w^{-1}\psi_w.$$

\vskip-1mm

(c) We have
$\psi_{w',w}\star\psi_{w}=\Lambda_{w}\,\psi_{w'}$ and
$\psi_{w'',w'}\star\psi_{w',w}=\Lambda_{w'}\,\psi_{w'',w}.$

\vskip2mm

(d) If $\ell(s)=0,1$ then
$[\th\! Z_{\Lbb,\Vb}^{s}]=\sum_{w,w'}\Lambda_{w,w'}^{s}\,\psi_{w,w'}$
in $\th\!\Zc''_{\Lbb,\Vb}$.
\vskip2mm

(e) We have $\Lambda_w=
\Eu(\th\!\een^*_{\Lbb,\Vb,w}\oplus\th\!\nen_{\Vb,w})$. 

\vskip2mm

(f) If $\ell(u)=0$ then
$\Lambda^u_{w,w'}=0$ if $w'\neq wu$, and 
$$\Lambda^e_{w,w}=\Lambda_w^{-1},\quad\Lambda^{\eps_1}_{w,w\eps_1}=
(\chi_{w(0)})^{\lambda_{i_{w(1)}}}\Lambda_w^{-1}=
(\chi_{w(1)})^{\lambda_{i_{w(0)}}}\Lambda_{w\eps_1}^{-1}.$$

\vskip1mm

(g) If $l=0,1,\dots,m-1$ then 
\vskip1mm
\itemitem\itemitem{$\bullet$} 
either $ws_l\notin W_w$ and 
$$\Lambda^{s_l}_{w,ws_l}=
\Lambda^{s_l}_{w,w}=
\Eu(\th\!\een_{\Lbb,\Vb,w,ws_l}^*\oplus\th\!\nen_{\Vb,w})^{-1},$$
\vskip1mm
\itemitem\itemitem{$\bullet$} or $ws_l\in W_w$ and 
$$\aligned
&\Lambda^{s_l}_{w,w}=
\Eu(\th\!\een_{\Lbb,\Vb,w,ws_l}^*\oplus
\th\!\nen_{\Vb,w}\oplus\th\!\men_{\Vb,w,ws_l})^{-1},
\hfill\vspace{2mm} &\Lambda^{s_l}_{w,ws_l}=
\Eu(\th\!\een_{\Lbb,\Vb,w,ws_l}^*\oplus\th\!\nen_{\Vb,w}
\oplus\th\!\men_{\Vb,ws_l,w})^{-1}.
\endaligned$$

\endproclaim

\noindent{\sl Proof :} Parts $(a)$ to $(d)$ are left to the reader. 
The fiber at $\phi_{\Vb,w}$ of the vector
bundle $$p:\th\widetilde F_{\Lbb,\Vb}\to\th\! F_\Vb$$ is isomorphic to 
$\th\!\een_{\Lbb,\Vb,w}$
as a $T$-module. Thus the cotangent space to $\th\widetilde
F_{\Lbb,\Vb}$ at the point $\phi_{\Vb,w}$ is isomorphic to
$\th\!\een_{\Lbb,\Vb,w}^*\oplus\th \nen_{\Vb,w}$ as a $T$-module. 
This yields $(e)$.
Next,
observe that the variety $\th\! Z^s_{\Lbb,\Vb}$ is smooth if 
$\ell(s)\leqslant 1$.
First, assume that $\ell(u)=0$. 
The fiber at $\phi_{\Vb,w,w'}$ of the vector bundle $$q:\th\! Z_{\Lbb,\Vb}^u\to
\th\! F_\Vb\times \th\! F_\Vb$$ is isomorphic to $\th\!\een_{\Lbb,\Vb,w,wu}$
as a $T$-module if $w'=wu$ and it is zero else. 
Thus we have
$$\aligned
\Lambda^{u}_{w,wu}
&=
\Eu(\th\den^*_{\Lbb,\Vb,w,wu})
\Eu(\th\!\een_{\Lbb,\Vb,w}^*)^{-1}
\Eu(\th\! F_\Vb,\phi_{\Vb,w})^{-1},\vspace{2mm}
&=
\Eu(\th\den^*_{\Lbb,\Vb,w,wu})
\Lambda_w^{-1}.
\endaligned$$ 
Therefore, Section 4.9 yields 
$$
\Lambda^{u}_{w,wu}=\cases
\Lambda_w^{-1}&\roman{if}\ u=e,\cr
(-\chi_{w(1)})^{\lambda_{i_{w(1)}}}\Lambda_w^{-1}
&\roman{if}\ u=\eps_1.
\endcases
$$ 
Note that
$$\chi_{w(1-l)}=-\chi_{w(l)},\quad\forall w,l.$$
This yields
$(f)$. 
Finally, let us concentrate on  part $(g)$.
The fiber at $\phi_{\Vb,w,w'}$ of the vector bundle 
$$q:\th\! Z_{\Lbb,\Vb}^{s_l}\to
\th\! F_\Vb\times \th\! F_\Vb$$ is isomorphic to $\th\!\een_{\Lbb,\Vb,w,ws_l}$
as a $T$-module if $w'=w,ws_l$ and it is zero else. Therefore, we have
$$\Lambda^{s_l}_{w,w'}=\cases
\Eu(\th\!\een_{\Lbb,\Vb,w,ws_l}^*)^{-1}
\Eu(\th\bar O_\Vb^{s_l},\phi_{\Vb,w,w'})^{-1}&\roman{if}\ w'=w,ws_l\cr
0&\roman{else}.\endcases$$ 
Next, by Lemma 4.6$(c)$, if $ws_l\notin W_w$ the cotangent spaces to the variety $\th\bar
O_\Vb^{s_l}$ at the points $\phi_{\Vb,w,ws_l}$ and $\phi_{\Vb,w,w}$ are given by
$$\gathered
T^*_{{w,ws_l}}\th\bar O_\Vb^{s_l}=T^*_{{w,ws_l}}\th O_\Vb^{s_l}=
\th\!\nen_{\Vb,w,ws_l}=\th\!\nen_{\Vb,w},
\vspace{2mm}
T^*_{{w,w}}\th\bar O_\Vb^{s_l}=T^*_{{w,w}}\th O_\Vb^{e}=
\th\!\nen_{\Vb,w,w}=\th\!\nen_{\Vb,w}.
\endgathered$$
Thus they are both isomorphic to $\th\!\nen_{\Vb,w}$ as $T$-modules.
Similarly, if $ws_l\in W_w$ the cotangent spaces to the variety $\th\bar
O_\Vb^{s_l}$ at the points $\phi_{\Vb,w,ws_l}$, $\phi_{\Vb,w,w}$ are given by
$$\gathered
T^*_{{w,ws_l}}\th\bar O_\Vb^{s_l}=
\th\!\nen_{\Vb,w}\oplus\th\!\men_{\Vb,ws_l,w},
\vspace{2mm}
T^*_{{w,w}}\th\bar O_\Vb^{s_l}=
\th\!\nen_{\Vb,w}\oplus\th\!\men_{\Vb,w,ws_l},
\endgathered$$
because 
$\Lie(\th\!P_{\Vb,w,ws_l})/\Lie(\th\!B_{\Vb,w})$
is dual to
$\th\!\men_{\Vb,w,ws_l}=
\th\!\nen_{\Vb,w}/\th\!\nen_{\Vb,w,ws_l}$
as a $T$-module.

\qed

\vskip3mm

\subhead 4.14.~Description of the
$\th\!\Zc_{\Lbb,\Vb}$-action on $\th\Fc_{\Lbb,\Vb}$\endsubhead
Using the computations in the previous proposition we can now
describe explicitly the representation of $\th\!\Zc_{\Lbb,\Vb}$ in $\th\Fc_{\Lbb,\Vb}$.
For $k=0,1,\dots,m-1$ let $\sigma_{\Lbb,\Vb}(k)$ be the fundamental class of
$\th\! Z^{s_k}_{\Lbb,\Vb}$ in $\th\!\Zc_{\Lbb,\Vb}^{\leqslant s_k}$. Next, let
$\pi_{\Lbb,\Vb}(1)$ be the fundamental class of $\th\! Z_{\Lbb,\Vb}^{\eps_1}$ in
$\th\!\Zc_{\Lbb,\Vb}^{\eps_1}$.
Finally, for $l=1,2,\dots, m$ the pull-back of the first equivariant Chern
class of the line bundle $\bigoplus_\ib\Oc_{\Lbb,\Vb,\ib}(l)^{-1}$ 
by the obvious map
$$\th\! Z_{\Lbb,\Vb}^e
\to\th\widetilde F_{\Lbb,\Vb}$$ belongs to $H^*_{\th\! G_\Vb}(\th
Z_{\Lbb,\Vb}^e,\k)$. So it yields an element $\varkappa_{\Lbb,\Vb}(l)$ in
$\th\!\Zc_{\Lbb,\Vb}^e$. 
Now, recall that $\th\!\Zc_{\Lbb,\Vb}^{\leqslant w}$
embeds into $\th\!\Zc_{\Lbb,\Vb}$. Thus the classes
$\sigma_{\Lbb,\Vb,\ib',\ib}(k)$,
$\pi_{\Lbb,\Vb,\ib',\ib}(1)$
and
$\varkappa_{\Lbb,\Vb,\ib',\ib}(l)$ can all be regarded as elements
of $\th\!\Zc_{\Lbb,\Vb}$.
We write
$$\gathered
\sigma_{\Lbb,\Vb,\ib',\ib}(k)=1_{\Lbb,\Vb,\ib'}\star\sigma_{\Lbb,\Vb}(k)\star
1_{\Lbb,\Vb,\ib},\vspace{2mm}
\pi_{\Lbb,\Vb,\ib',\ib}(1)=1_{\Lbb,\Vb,\ib'}\star\pi_{\Lbb,\Vb}(1)\star
1_{\Lbb,\Vb,\ib},
\vspace{2mm}
\varkappa_{\Lbb,\Vb,\ib',\ib}(l)=1_{\Lbb,\Vb,\ib'}\star 
\varkappa_{\Lbb,\Vb}(l)\star
1_{\Lbb,\Vb,\ib}.
\endgathered$$  
For a sequence $\ib=(i_{1-m},\dots,i_{m-1},i_{m})$ 
and integers $l=1-m,\dots,m-1,m$ and $k=1,\dots,m-1,m$, we set
$$\gathered
\l_\ib(l)=\l_{i_l},
\quad
h_\ib(k)=\cases
-1& \roman{if}\ s_k\ib=\ib,\vspace{2mm}
h_{i_k,i_{k+1}}&\roman{if}\ s_k\ib\neq\ib,\,k\neq 0,\vspace{2mm}
h_{i_{0},i_{2}}&\roman{if}\ s_0\ib\neq\ib,\,k=0.
\endcases\endgathered$$ 
Finally, recall that $\th\!\Zc_{\Lbb,\Vb}$ acts on
$\th\Fc_{\Lbb,\Vb}=\bigoplus_\ib\th\Fc_{\Lbb,\Vb,\ib}$ and that we identify
$\th\Fc_{\Lbb,\Vb,\ib}$ with 
$$\k[x_\ib(1),x_\ib(2),\dots x_\ib(m)]
=\th\Fb_{\Lbb,\Vb,\ib}$$ 
via (4.5). The later is given the obvious $\kb$-algebra structure.

\proclaim{4.15. Proposition} For $\ib,\ib',\ib''$ in $\th I^\nu$
and $f$ in $\th\Fc_{\Lbb,\Vb,\ib}$ the following hold :

\vskip1mm

(a)  $1_{\Lbb,\Vb,\ib'}\star f=f$ if $\ib=\ib'$ and 
$1_{\Lbb,\Vb,\ib'}\star f=0$ else.
\vskip1mm

(b)  $\varkappa_{\Lbb,\Vb,\ib'',\ib'}(l)\star f=0$ unless
$\ib''=\ib'=\ib$ and $\varkappa_{\Lbb,\Vb}(l)\star f=x_\ib(l)f$.
\vskip1mm

(c)  $\pi_{\Lbb,\Vb,\ib'',\ib'}(1)\star f=0$ unless
$\ib'=\ib$, $\ib''=\eps_1\ib$ and 
$$\pi_{\Lbb,\Vb}(1)\star f=x_{\eps_1\ib}(0)^{\l_{\eps_1\ib}(0)}\eps_1(f).$$
\vskip1mm

(d)  $\sigma_{\Vb,\ib'',\ib'}(k)\star f=0$ unless $\ib'=\ib$ and
$\ib''=s_k\ib$ or $\ib$, and we have 
\vskip1mm\itemitem\itemitem{$\bullet$} if $s_k\ib=\ib$ and $k\neq 0$ then
$$\sigma_{\Vb}(k)\star f=
(x_\ib(k+1)-x_\ib(k))^{h_\ib(k)}(s_k(f)-f),
$$ 
\vskip1mm
\itemitem\itemitem{$\bullet$} if $s_0\ib=\ib$ then
$$\sigma_{\Vb}(0)\star f=
(x_\ib(2)-x_\ib(0))^{h_\ib(0)}x_\ib(1)^{\lambda_\ib(1)}
x_\ib(2)^{\lambda_\ib(2)}(s_0(f)-f),
$$ 
\vskip1mm
\itemitem\itemitem{$\bullet$} if $s_k\ib\neq\ib$
and $k\neq 0$ then
$$\aligned
&\sigma_{\Vb,s_k\ib,\ib}(k)\star f=
(x_{s_k\ib}(k+1)-x_{s_k\ib}(k))^{h_{s_k\ib}(k)} s_k(f),
\vspace{2mm}
&\sigma_{\Vb,\ib,\ib}(k)\star f=
(x_\ib(k+1)-x_\ib(k))^{h_{\ib}(k)} f
\endaligned$$
\vskip1mm
\itemitem\itemitem{$\bullet$} if $s_0\ib\neq\ib$ then
$$\aligned
&\sigma_{\Vb,s_0\ib,\ib}(0)\star f=
(x_{s_0\ib}(2)-x_{s_0\ib}(0))^{h_{s_0\ib}(0)}
x_{s_0\ib}(1)^{\lambda_{s_0\ib}(1)}x_{s_0\ib}(2)^{\lambda_{s_0\ib}(2)}s_0(f),\vspace{2mm}
&\sigma_{\Vb,\ib,\ib}(0)\star f=
(x_\ib(2)-x_\ib(0))^{h_{\ib}(0)}x_\ib(1)^{\lambda_{\ib}(1)}x_{\ib}(2)^{\lambda_\ib(2)}f.
\endaligned$$
\endproclaim

\noindent{\sl Proof :} Parts $(a)$, $(b)$ are  left
to the reader. Let $w\in W_\ib$. Recall that
$\th\Fc_{\Lbb,\Vb,\ib}\subset\th\Fc'_{\Lbb,\Vb,\ib}$ 
and that $\psi_w$ lies in $\th\Fc'_{\Lbb,\Vb,\ib}$.
Under the map (4.6) the multiplication in
$\th\Fc_{\Lbb,\Vb,\ib}$ and the $\Sb$-action on
$\th\Fc'_{\Lbb,\Vb,\ib}$ are related by the following formula
$$f(-x_\ib(1),\dots,-x_\ib(m))\,\psi_w=
w(f)\,\psi_w.$$
Further we have $\eps_1(\ib_w)=\ib_{w\eps_1}$.
Therefore, part $(c)$ follows from the following computation,
see Proposition 4.13$(c),(d),(f)$,
$$\aligned
[\th\! Z_{\Lbb,\Vb}^{\eps_1}]\star\psi_w
=\Lambda^{\eps_1}_{w\eps_1,w}\Lambda_w\psi_{w\eps_1}
=(\chi_{w\eps_1(1)})^{\l_{i_{w\eps_1(0)}}}\psi_{w\eps_1}
=(-\chi_{w\eps_1(0)})^{\l_{i_{w\eps_1(0)}}}\psi_{w\eps_1},
\endaligned$$
where $i_{w\eps_1(l)}$ is the $l$-th component of the sequence $\ib_{w\eps_1}$.
Let us concentrate on $(d)$. The first claim is
obvious because $\th\! Z^{s_k}_{w,w'}=\emptyset$ unless $w'=w,ws_k$ 
by Lemma 4.6, and $\ib_{ws_k}=s_k\ib_w$. Now, given $\ib'=\ib$ or
$s_k\ib$ we must compute the linear operator
$$\th\Fc_{\Vb,\ib}\to\th\Fc_{\Vb,\ib'},\quad
f\mapsto\sigma_{\Vb,\ib',\ib}(k)\star f.\leqno(4.7)$$
Proposition 4.13$(b)$ yields an embedding
$$\th\Fc_{\Vb,\ib}\to \bigoplus_{w\in W_{\ib}}\Qb\psi_w,
\quad f(-x_\ib(1),\dots,-x_\ib(m))\mapsto
\sum_{w\in W_{\ib}}w(f)\Lambda_w^{-1}\psi_w.$$
Under this inclusion the map (4.7) is of the following form
$$\sum_{w\in W_{\ib}}w(f)\Lambda_w^{-1}\psi_w
\mapsto\sum_{w'\in W_{\ib'}} g_{w'}\psi_{w'},\quad
g_{w'}=\sum_{w\in W_{\ib}}w(f)\Lambda^{s_k}_{w',w}$$
by Proposition 4.13$(c),(d)$.
We claim that the right hand side is the image of a polynomial
$g$ in $\th\Fc_{\Vb,\ib'}$ that we'll compute explicitly. 
The polynomial $g$ is completely determined
by the following relations $$
g_{w'}=w'(g)\Lambda_{w'}^{-1},\quad\forall
w'\in W_{\ib'}.\leqno(4.8)$$ In the rest of the proof we'll
fix $w,w'$ in the following way
$$w\in W_{\ib},\quad w'\in W_{\ib'},\quad
w'=w\ \roman{or}\ ws_k.$$ In particular we have $\ib=\ib_w$,
$\ib'=\ib_{w'}$, and $\ib'=\ib$ or $s_k\ib$.

\vskip1mm

{$(i)$} First, assume that $s_k\ib=\ib$. Then $\ib'=\ib$,
$w's_k\in W_{w'}$,  and we have
$$g_{w'}=w'(f)\Lambda^{s_k}_{w',w'}+
w's_k(f)\Lambda^{s_k}_{w',w's_k}.$$ 
Section 4.9 and Proposition 4.13 yield 
$$
\gathered \Lambda_{w'} ={\Eu}(\th\!\een_{\Lbb,\Vb,w'}^*\oplus\th\!\nen_{\Vb,w'}),\vspace{2mm}
\Lambda^{s_k}_{w',w'}
=\Eu(\th\!\een_{\Lbb,\Vb,w',w's_k}^*\oplus
\th\!\nen_{\Vb,w'}\oplus\th\!\men_{\Vb,w',w's_k})^{-1}
,\vspace{2mm}
\Lambda^{s_k}_{w',w's_k}=
\Eu(\th\!\een_{\Lbb,\Vb,w',w's_k}^*\oplus
\th\!\nen_{\Vb,w'}\oplus\th\!\men_{\Vb,w's_k,w'})^{-1},
\vspace{2mm}
\Eu(\th\!\men_{\Vb,w',w's_k})=-\Eu(\th\!\men_{\Vb,w's_k,w'})=w'(\a_k).
\endgathered$$
So we have
$$
\aligned 
\Lambda^{s_k}_{w',w'}
=\Eu(\th\!\den_{\Lbb,\Vb,w',w's_k}^*)
w'(\a_k)^{-1}\Lambda_{w'}^{-1}=-
\Lambda^{s_k}_{w',w's_k}.
\endaligned
$$
Therefore we obtain
$$g_{w'}=w'(f-s_k(f))
\,\Eu(\th\!\den_{\Lbb,\Vb,w',w's_k}^*)
w'(\a_k)^{-1}\Lambda_{w'}^{-1}.$$ 
Now, assume that $k\neq 0$. 
There is no arrow joining $i_{w'(k)}$ and $i_{w'(k+1)}$, because $i_{w'(k)}= i_{w'(k+1)}$.
Thus Section 4.9 yields $$\Eu(\th\!\den_{\Lbb,\Vb,w',w's_k}^*)=1.$$
Hence $$
\aligned
g_{w'}
&=w'(f-s_k(f))
w'(\a_k)^{-1}\Lambda_{w'}^{-1}\vspace{2mm}
&=w'(g)\Lambda_{w'}^{-1},\vspace{2mm}
g&=(f-s_k(f))\a_k^{-1}.\endaligned$$
Next, assume that $k=0$. 
There is no arrow joining $i_{w'(0)}$ and $i_{w'(2)}$.
Thus Section 4.9 yields 
$$\Eu(\th\!\den_{\Lbb,\Vb,w',w's_0}^*)=
(-\chi_{w'(1)})^{\l_{i_{w'(1)}}}(-\chi_{w'(2)})^{\l_{i_{w'(2)}}}.$$
Therefore we have
$$\aligned
g_{w'}&=w'(f-s_0(f))
(-\chi_{w'(1)})^{\l_{i_{w'(1)}}}(-\chi_{w'(2)})^{\l_{i_{w'(2)}}}
w'(\a_0)^{-1}
\Lambda_{w'}^{-1}\vspace{2mm}
&=w'(g)\Lambda_{w'}^{-1},
\vspace{2mm} 
g&=(f-s_0(f))(-\chi_{1})^{\l_{i_{w'(1)}}}(-\chi_{2})^{\l_{i_{w'(2)}}}\a_0^{-1}.\endaligned$$

\vskip2mm

{$(ii)$} Finally, assume that $s_k\ib\neq\ib$, i.e., that $ws_k\notin W_w$.
Section 4.9 and Proposition 4.13 yield
$$\gathered
\Eu(\th\!\nen_{\Vb,ws_k})=\Eu(\th\!\nen_{\Vb,w}),\vspace{2mm} \Lambda_{w}
={\Eu}(\th\!\een_{\Lbb,\Vb,w}^*\oplus\th\!\nen_{\Vb,w}),\vspace{2mm}
\Lambda_{ws_k}
={\Eu}(\th\!\een_{\Lbb,\Vb,ws_k}^*\oplus\th\!\nen_{\Vb,w}),\vspace{2mm}
\Lambda^{s_k}_{w,w} 
=\Eu(\th\!\een_{\Lbb,\Vb,w,ws_k}^*\oplus\th\!\nen_{\Vb,w})^{-1}=
\Lambda^{s_k}_{ws_k,w}.
\endgathered$$
So we have
$$\gathered
\Lambda^{s_k}_{ws_k,w}\Lambda_{ws_k}=
\Eu(\th\!\den_{\Lbb,\Vb,ws_k,w}^*),
\vspace{2mm}
\Lambda^{s_k}_{w,w}\Lambda_{w}=
\Eu(\th\!\den_{\Lbb,\Vb,w,ws_k}^*).
\endgathered$$
Next, one of the two following alternatives holds :

\vskip1mm

\item{$\bullet$}
either $\ib'=s_k\ib$, $w'=ws_k$ and 
$$\aligned
g_{w'}
&=w's_k(f)\Lambda^{s_k}_{w',w's_k}
\vspace{2mm}
&=w's_k(f)(\Lambda^{s_k}_{ws_k,w}\Lambda_{ws_k})\Lambda_{w'}^{-1}
\vspace{2mm}
&=w's_k(f)\,\Eu(\th\!\den_{\Lbb,\Vb,w',w's_k}^*)\,\Lambda_{w'}^{-1}.
\endaligned$$ 
\vskip1mm

\item{$\bullet$}
or $\ib'=\ib$, $w'=w$ and $$\aligned
g_{w'}
&=w'(f)\Lambda^{s_k}_{w',w'}
\vspace{2mm}
&=w'(f)(\Lambda^{s_k}_{w,w}\Lambda_{w})\Lambda_{w'}^{-1}
\vspace{2mm}
&=
w'(f)\,\Eu(\th\!\den_{\Lbb,\Vb,w',w's_k}^*)\,\Lambda_{w'}^{-1}.
\endaligned$$ 
\vskip1mm
\noindent
Now we consider the cases $k\neq 0$ and $k=0$.
First, assume that $k\neq 0$. By Section 4.9 we have 
$$\gathered
\Eu(\th\!\den_{\Lbb,\Vb,w',w's_k})=w'(\a_k)^{h_{i_{w'(k)},i_{w'(k+1)}}}.
\cr\endgathered$$
Thus (4.8) holds with $$g=s_k(f)(-\a_k)^{h_{i_{w'(k)},i_{w'(k+1)}}}$$
in the first case and
with $$g=f(-\a_{k})^{h_{i_{w'(k)},i_{w'(k+1)}}}$$ in the second one.
Next, assume that $k=0$. By Section 4.9 we have
$$\gathered
\Eu(\th\!\den_{\Lbb,\Vb,w',w's_0})
=w'(\chi_{1})^{\l_{i_{w'(1)}}}w'(\chi_{2})^{\l_{i_{w'(2)}}}w'(\a_0)^{h_{i_{w'(0)},i_{w'(2)}}}.
\cr\endgathered$$
Thus (4.8) holds with 
$$g=s_0(f)(-\chi_{1})^{\l_{i_{w'(1)}}}(-\chi_{2})^{\l_{i_{w'(2)}}}
(-\a_0)^{h_{i_{w'(0)},i_{w'(2)}}}$$ in the first case and
with  
$$g=f(-\chi_{1})^{\l_{i_{w'(1)}}}(-\chi_{2})^{\l_{i_{w'(2)}}}(-\a_{0})^{h_{i_{w'(0)},i_{w'(2)}}}$$
in the second one.

\qed


\vskip3mm

\subhead 4.16.~Description of the graded $\kb$-algebra
$\th\Zb_{\Lbb,\Vb}$\endsubhead 
We can use the previous computations concerning
$\th\Zc_{\Lbb,\Vb}$ to get informations on $\th\Zb_{\Lbb,\Vb}$.
The action of $$\gathered
1_{\Lbb,\Vb,\ib},\quad\varkappa_{\Lbb,\Vb}(l),\quad
\sigma_{\Lbb,\Vb}(k),\quad\pi_{\Lbb,\Vb}(1),\endgathered$$
yields linear operators in
$\End(\th\Fc_{\Lbb,\Vb})$. Recall that $\th\Fc_{\Lbb,\Vb}$ is a faithful left
$\th\!\Zc_{\Lbb,\Vb}$-module and that there are canonical isomorphisms 
$$\th\Zb_{\Lbb,\Vb}=\th\!\Zc_{\Lbb,\Vb},\quad
\th\Fb_{\Lbb,\Vb}=\th\Fc_{\Lbb,\Vb}.$$ 
Thus the graded left $\th\Zb_{\Lbb,\Vb}$-module
$\th\Fb_{\Lbb,\Vb}$ is also faithful. Recall also that
$$\th\Fb_{\Lbb,\Vb}=\bigoplus_{\ib\in\th\! I^\nu}\th\Fb_{\Lbb,\Vb,\ib},\quad
\th\Fb_{\Lbb,\Vb,\ib}=\kb[x_\ib(1),x_\ib(2),\dots,x_\ib(m)].$$ We obtain
the following.

\proclaim{4.17.~Theorem}
The graded $\kb$-algebra $\th\Zb_{\Lbb,\Vb}$ 
is isomorphic to a graded $\kb$-subalgebra
of $\End(\th\Fb_{\Lbb,\Vb})$ which contains the linear operators
$$\gathered 1_{\Lbb,\Vb,\ib},\quad\varkappa_{\Lbb,\Vb,\ib}(l),\quad
\sigma_{\Lbb,\Vb,\ib}(k),\quad\pi_{\Lbb,\Vb,\ib}(1),\vspace{2mm}
\ib\in\th\! I^\nu,\quad k=0,1,\dots,m-1,\quad l=1,2,\dots,m,
\endgathered$$
defined as follows :

\vskip1mm \itemitem{$(a)$} $1_{\Lbb,\Vb,\ib}$ is the projection to
$\th\Fb_{\Lbb,\Vb,\ib}$ relatively to
$\bigoplus_{\ib'\neq\ib}\th\Fb_{\Lbb,\Vb,\ib'}$,

\vskip1mm \itemitem{$(b)$} $\varkappa_{\Lbb,\Vb,\ib}(l)=0$ on
$\th\Fb_{\Lbb,\Vb,\ib'}$ if $\ib'\neq\ib$, and it acts by multiplication by
$x_\ib(l)$ on $\th\Fb_{\Lbb,\Vb,\ib}$,

\vskip1mm \itemitem{$(c)$} $\sigma_{\Lbb,\Vb,\ib}(k)=0$  on
$\th\Fb_{\Lbb,\Vb,\ib'}$ if  $\ib'\neq\ib$, and it takes 
a polynomial $f$ in $\th\Fb_{\Lbb,\Vb,\ib}$ to \vskip-2mm
$$\matrix
&(x_\ib(k+1)-x_\ib(k))^{h_\ib(k)}(s_k(f)-f)
\hfill& \text{if}\
s_k\ib=\ib\,,k\neq 0,\vspace{2mm}
&(x_\ib(2)-x_\ib(0))^{h_\ib(0)}x_\ib(1)^{\lambda_\ib(1)}
x_\ib(2)^{\lambda_\ib(2)}(s_0(f)-f)
\hfill& \text{if}\ s_k\ib=\ib\,,k=0,\vspace{2mm}
&(x_{s_k\ib}(k+1)-x_{s_k\ib}(k))^{h_{s_k\ib}(k)}s_k(f)\hfill&
\text{if}\ s_k\ib\neq \ib\,,k\neq 0,\vspace{2mm}
&(x_{s_0\ib}(2)-x_{s_0\ib}(0))^{h_{s_0\ib}(0)}
x_{s_0\ib}(1)^{\lambda_{s_0\ib}(1)}x_{s_0\ib}(2)^{\lambda_{s_0\ib}(2)}
s_0(f)&
\text{if}\ s_k\ib\neq \ib,\,k=0,
\endmatrix$$

\vskip1mm \itemitem{$(d)$} $\pi_{\Lbb,\Vb,\ib}(1)=0$  on
$\th\Fb_{\Lbb,\Vb,\ib'}$ if  $\ib'\neq\ib$, and it takes 
a polynomial $f$ in $\th\Fb_{\Lbb,\Vb,\ib}$ to 
$$x_{\eps_1\ib}(0)^{\l_{\eps_1\ib}(0)}\eps_1(f).$$ 

\vskip1mm

\noindent The degrees of these operators
are given by the following formulas
$$\aligned
&\deg(1_{\Lbb,\Vb,\ib})=0,\vspace{2mm}
&\deg(\varkappa_{\Lbb,\Vb,\ib}(l))=2,\vspace{2mm}
&\deg(\pi_{\Lbb,\Vb,\ib}(1))=2\lambda_{\eps_1\ib}(0),\vspace{2mm}
&\deg(\sigma_{\Lbb,\Vb,\ib}(0))=
2h_{s_0\ib}(0)+2\lambda_{s_0\ib}(1)+2\lambda_{s_0\ib}(2),\vspace{2mm}
&\deg(\sigma_{\Lbb,\Vb,\ib}(k))=
2h_{s_k\ib}(k)\quad\roman{if}\ k\neq 0.
\endaligned$$

\endproclaim
\vskip3mm

\subhead 4.18.~Shift of the grading\endsubhead 
We are mostly interested by the 
graded $\kb$-algebra $\th\Zb^\delta_{\Lbb,\Vb}$,
whose grading differs from the grading of
$\th\Zb^\delta_{\Lbb,\Vb}$. Let us compute the degree of the generators
of $\th\Zb^\delta_{\Lbb,\Vb}$.
We have
$$\th\Zb^\delta_{\Lbb,\Vb,s_k\ib,\ib}=
\th\Zb_{\Lbb,\Vb,s_k\ib,\ib}[d_{\l,\ib}-d_{\l,s_k\ib}].$$
Recall that $h_{\theta(i),j}=h_{\theta(j),i}$ for each $i,j$. 
Hence, an easy computation using Proposition 2.5 yields
$$\aligned
d_{\l,\ib}-d_{\l,s_k\ib}&=\cases
h_{\ib}(k)-h_{s_k\ib}(k)&\roman{if}\ k\neq0,\vspace{2mm}
h_{\ib}(0)-h_{s_0\ib}(0)+
\lambda_{\ib}(2)+\lambda_{\ib}(1)-\lambda_{s_0\ib}(2)-
\lambda_{s_0\ib}(1)&\roman{if}\ k=0,\endcases
\vspace{2mm}d_{\l,\ib}-d_{\l,\eps_1\ib}&=\lambda_{\eps_1\ib}(0)-\lambda_\ib(0).\endaligned$$
Therefore the grading of $\th\Zb^\delta_{\Lbb,\Vb}$ is given by the
following rules : 
$$\aligned
&\deg(1_{\Lbb,\Vb,\ib})=0,\vspace{2mm}
&\deg(\varkappa_{\Lbb,\Vb,\ib}(l))=2,\vspace{2mm}
&\deg(\pi_{\Lbb,\Vb,\ib}(1))= \lambda_{\ib}(0)+\lambda_\ib(1),\vspace{2mm}
&\deg(\sigma_{\Lbb,\Vb,\ib}(0))=-i_0\cdot
i_{2}+\lambda_{\ib}(-1)+\lambda_{\ib}(0)+\lambda_{\ib}(1)+\lambda_{\ib}(2),\vspace{2mm}
&\deg(\sigma_{\Lbb,\Vb,\ib}(k))=-i_k\cdot i_{k+1}
\quad\roman{if}\ k\neq 0.
\endaligned$$

\vskip2cm

\head 5.~The graded $\kb$-algebra $\th\Rb(\Gamma)_{\l,\nu}$\endhead

Fix a quiver $\Gamma$ with set of vertices $I$ and set of arrows $H$.
Fix an involution $\theta$ on $\Gamma$. 
Assume that $\Gamma$ has no 1-loops and that
$\theta$ has no fixed points. Fix a
dimension vector $\nu\neq 0$ in $\th\NN I$ and a
dimension vector $\l$ in $\NN I$. Set $|\nu|=2m$.

\vskip3mm

\subhead 5.1.~Definition of the graded $\kb$-algebra
$\th\Rb(\Gamma)_{\l,\nu}$\endsubhead 
Assume that $m>0$. We define a graded 
$\kb$-algebra $\th\Rb(\Gamma)_{\l,\nu}$ with 1 generated by $1_\ib$,
$\varkappa_{l}$, $\sigma_{k}$, $\pi_{1}$ 
with $\ib=(i_{1-m},\dots,i_m)$ in $\th\! I^\nu$, $k=1,\dots,m-1$, 
$l=1,2,\dots,m$,
modulo the following defining relations
\footnote"${}^{\dag}$"{We thank M. Kashiwara who indicate us an error
in a previous version of the relations}

\vskip2mm

\itemitem{$(a)$}
$1_\ib\,1_{\ib'}=\delta_{\ib,\ib'}1_\ib$,
\quad 
$\sigma_{k}1_\ib=
1_{s_k\ib}\sigma_{k}$,
\quad
$\varkappa_{l}1_\ib=
1_{\ib}\varkappa_{l}$,
\quad
$\pi_{1}1_\ib=
1_{\eps_1\ib}\pi_{1}$,

\vskip2mm

\itemitem{$(b)$}
$\varkappa_{l}\varkappa_{l'}=\varkappa_{l'}\varkappa_{l}$,\quad
$\pi_1\varkappa_{l}=\varkappa_{\eps_1(l)}\pi_1$,

\vskip2mm

\itemitem{$(c)$}
$\sigma_{k}^21_\ib=
Q_{i_k,i_{k+1}}(\varkappa_{k+1},\varkappa_{k})1_\ib$,\quad
$\pi_1^21_\ib=
\varkappa_0^{\lambda_{i_0}}
\varkappa_1^{\lambda_{i_1}}1_\ib$,

\vskip2mm

\itemitem{$(d)$}
$\sigma_{k}\sigma_{k'}=\sigma_{k'}\sigma_{k}$ if $k\neq k'\pm 1$,\quad
$\pi_1\sigma_{k}=\sigma_{k}\pi_1$ if $k\neq 1$,

\vskip2mm

\itemitem{$(e)$}
${\ds(\sigma_1\pi_1)^21_\ib=(\pi_1\sigma_1)^21_\ib+
\delta_{i_0,i_2}(-1)^{\l_{i_2}}
{\varkappa_0^{\l_{i_1}+\l_{i_2}}-\varkappa_2^{\l_{i_1}+\l_{i_2}}\over
\varkappa_0-\varkappa_2}\sigma_11_\ib,}$

\vskip2mm

\itemitem{$(f)$}
$(\sigma_{k+1}\sigma_{k}\sigma_{k+1}-
\sigma_{k}\sigma_{k+1}\sigma_{k})1_\ib=$
$$=
\delta_{i_k,i_{k+2}}{Q_{i_k,i_{k+1}}(\varkappa_{k+1},\varkappa_{k})
-Q_{i_k,i_{k+1}}(\varkappa_{k+1},\varkappa_{k+2})\over
\varkappa_{k}-\varkappa_{k+2}}1_\ib,
$$

\vskip2mm

\itemitem{$(g)$}
$(\sigma_{k}\varkappa_{l}-\varkappa_{s_k(l)}\sigma_{k})1_\ib
=\cases -1_\ib&\ \roman{if}\  l=k,\, i_k=i_{k+1}, \cr 1_\ib&\
\roman{if}\  l=k+1,\, i_k=i_{k+1},\cr 0&\ \roman{else}.
\endcases$

\vskip2mm

\noindent 
Here $\delta_{i,j}$ is the Kronecker symbol,
$\varkappa_{1-l}=-\varkappa_l$, and 
$$Q_{i,j}(u,v)=\cases(-1)^{h_{i,j}}(u-v)^{-i\cdot j}&\ \roman{if}\
i\neq j, \cr 0&\ \roman{else}.\endcases\leqno(5.1)$$ 
We'll abbreviate
$\sigma_{\ib,k}=\sigma_{k}1_\ib$,
$\varkappa_{\ib,l}=\varkappa_{l}1_\ib$,
and $\pi_{\ib,1}=\pi_{1}1_\ib$.
The grading on
$\th\Rb(\Gamma)_{\l,\nu}$ is given by the following rules :
$$\aligned
&\deg(1_{\ib})=0,\vspace{2mm}
&\deg(\varkappa_{\ib,l})=2,\vspace{2mm}
&\deg(\pi_{\ib,1})= \lambda_{i_0}+\lambda_{i_1},\vspace{2mm}
&\deg(\sigma_{\ib,k})=-i_k\cdot i_{k+1}.
\endaligned$$
If $\nu=0$ we set $\th\Rb(\Gamma)_{\l,\nu}=\kb$ as a graded
$\kb$-algebra.
Let $\omega$ be the unique anti-involution of the graded
$\kb$-algebra $\th\Rb(\Gamma)_{\l,\nu}$
which fixes $1_\ib$, $\varkappa_l$, $\sigma_k$, $\pi_1$. 

\vskip3mm

\subhead 5.2.~Remarks\endsubhead
$(a)$
We may set $\sigma_0=\pi_1\sigma_1\pi_1$. We have
$$\deg(\sigma_{0}1_\ib)=-i_0\cdot
i_{2}+\lambda_{i_{-1}}+\lambda_{i_0}+
\lambda_{i_1}+\lambda_{i_2}.$$

\vskip1mm

$(b)$
We may also set $\pi_l=\sigma_{l-1}\dots\sigma_2\sigma_1\pi_1\sigma_1\sigma_2\dots\sigma_{l-1}$.
We have
$$\deg(\pi_{l}1_\ib)=-(i_1+i_2+\dots+i_{l-1})\cdot(i_l+i_{1-l})+
\lambda_{i_l}+\lambda_{i_{1-l}}.$$

\vskip3mm

\subhead 5.3.~The polynomial representation and the PBW theorem\endsubhead 
Given any objects $\Vb$ in
$\th\Vcb_\nu$ and $\Lbb$ in $\Vcb_\l$ we abbreviate 
$$\th\Fb_\nu=\th\Fb_{\Lbb,\Vb},
\quad
\th\Fb_{\ib}=\th\Fb_{\Lbb,\Vb,\ib}, 
\quad
\th\Sb_\nu=\th\Sb_{\Vb}.$$

\proclaim{5.4.~Proposition} There is an unique graded
$\kb$-algebra morphism
$\th\Rb(\Gamma)_{\l,\nu}\to\End(\th\Fb_{\nu})$ such that,
for each $\ib\in\th\! I^\nu$, $k=0,1,\dots,m-1$, $l=1,2,\dots,m$, we have
$$\gathered 1_\ib\mapsto 1_{\Lbb,\Vb,\ib},\quad
\varkappa_{\ib,l}\mapsto\varkappa_{\Lbb,\Vb,\ib}(l),\quad
\sigma_{\ib,k}\mapsto\sigma_{\Lbb,\Vb,\ib}(k),\quad\pi_{\ib,1}\mapsto\pi_{\Lbb,\Vb,\ib}(1).
\endgathered$$
\endproclaim

\noindent{\sl Proof :}
The defining relations  of $\th\Rb(\Gamma)_{\l,\nu}$
are checked by a direct computation.
Let us (only) give a few indications concerning the relation 5.1$(e)$.
We have 
$$\sigma_1 1_\ib=(\varkappa_1-\varkappa_2)^{h_{s_1\ib}(1)}
(s_1-\delta_{i_1,i_2})1_\ib,\quad
\pi_1 1_\ib=\varkappa_0^{\l_{i_1}}\eps_11_\ib.$$
This yields
$$\gathered\sigma_1\pi_1 1_\ib=(\varkappa_1-\varkappa_2)^{h_{s_0\ib}(0)}
(s_1-\delta_{i_0,i_2})\varkappa_0^{\l_{i_1}}\eps_11_\ib,
\vspace{2mm}
\pi_1\sigma_1 1_\ib=\varkappa_0^{\l_{i_{2}}}\eps_1
(\varkappa_1-\varkappa_2)^{h_{s_1\ib}(1)}(s_1-\delta_{i_1,i_2})1_\ib.
\endgathered$$
Therefore we have
$$\gathered
(\sigma_1\pi_1)^2 1_\ib=
(\varkappa_1-\varkappa_2)^{h_{s_1\ib}(1)}
(s_1-\delta_{i_1,i_2})
\varkappa_0^{\l_{i_{2}}}
\eps_1
(\varkappa_1-\varkappa_2)^{h_{s_0\ib}(0)}
(s_1-\delta_{i_0,i_2})
\varkappa_0^{\l_{i_1}}
\eps_11_\ib,
\vspace{2mm}
(\pi_1\sigma_1)^2 1_\ib=
\varkappa_0^{\l_{i_{1}}}\eps_1
(\varkappa_1-\varkappa_2)^{h_{s_0\ib}(0)}
(s_1-\delta_{i_0,i_2})
\varkappa_0^{\l_{i_{2}}}\eps_1
(\varkappa_1-\varkappa_2)^{h_{s_1\ib}(1)}
(s_1-\delta_{i_1,i_2})1_\ib.
\endgathered$$
Hence we have
$$\gathered
(\sigma_1\pi_1)^2 1_\ib=
(\varkappa_1-\varkappa_2)^{h_{s_1\ib}(1)}
(\varkappa_0-\varkappa_2)^{h_{s_0\ib}(0)}A,
\vspace{2mm}
(\pi_1\sigma_1)^2 1_\ib=
(\varkappa_1-\varkappa_2)^{h_{s_1\ib}(1)}
(\varkappa_0-\varkappa_2)^{h_{s_0\ib}(0)}B,
\endgathered$$
where
$$\gathered
A=
(s_1-\delta_{i_1,i_2})
\varkappa_0^{\l_{i_{2}}}
\eps_1
(s_1-\delta_{i_0,i_2})
\varkappa_0^{\l_{i_1}}
\eps_11_\ib,
\vspace{2mm}
B=
\varkappa_0^{\l_{i_{1}}}\eps_1
(s_1-\delta_{i_0,i_2})
\varkappa_0^{\l_{i_{2}}}\eps_1
(s_1-\delta_{i_1,i_2})1_\ib.
\endgathered$$
If $i_0\neq i_2$ it is easy to see that $A=B$.
If $i_0=i_2$ a direct computation yields
$$B-A=\bigl(\varkappa_2^{\l_{i_1}}(-\varkappa_2)^{\l_{i_2}}-
\varkappa_1^{\l_{i_2}}(-\varkappa_1)^{\l_{i_1}}\bigr)s_1.$$
The rest of the computation is left to the reader.

\qed

\vskip3mm

The $\kb$-algebra $\th\Rb(\Gamma)_{\l,\nu}$
is a left graded $\th\Fb_\nu$-module such that
$x_\ib(l)$ acts by the
left multiplication with the element
$\varkappa_{\ib,l}$ for each $l=1,2,\dots, m$. 
To unburden the notation we may write $\varkappa_{\ib,l}=x_\ib(l)$.
The following convention is important.

\vskip2mm

\item{}
{\it From now on we'll regard $W_m$ as a Weyl group of type $B_m$,
with the set of simple reflections $\{s_1,s_2,\dots,s_m\}$ where $s_m=\eps_1$,
rather than an extended Weyl group of type $D_m$ as in Section $4.2$}.
\vskip2mm

\noindent
For $w\in W_m$ 
we choose a reduced decomposition $\dot w$ of $w$.
By the observation above $\dot w$ is a minimal
decomposition of the following form
$$w=s_{k_1}s_{k_2}\cdots s_{k_r},\quad 0<k_1,k_2,\dots,k_r\leqslant m,\quad 
s_m=\eps_1.$$
We define an element
$\sigma_{\dot w}$ in $\th\Rb(\Gamma)_{\l,\nu}$
by the following formula
$$\quad\sigma_{\dot w}=\sum_\ib 1_\ib\sigma_{\dot w} ,\quad
1_\ib\sigma_{\dot w}=\cases
1_\ib &\roman{if}\ r=0\vspace{2mm}
1_\ib\sigma_{k_1} \sigma_{k_2}\cdots\sigma_{k_r}
&\roman{else},
\endcases\leqno(5.2)$$ 
where we have set $\sigma_m=\pi_1$. Observe that $\sigma_{\dot w}$ may depend on the choice of the
reduced decomposition $\dot w$.

\proclaim{5.5.~Proposition} The $\kb$-algebra
$\th\Rb(\Gamma)_{\l,\nu}$ is a free (left or right)
$\th\Fb_\nu$-module with basis $\{\sigma_{\dot w};\,w\in W_m\}$.
Its rank is $2^mm!$. The operator $1_\ib\sigma_{\dot w}$ is homogeneous and its degree is independent of 
the choice of the reduced decomposition $\dot w$.
\endproclaim

\noindent{\sl Proof :}
The $\kb$-space $\th\Rb(\Gamma)_{\l,\nu}$ is filtered with
$1_\ib$, $\varkappa_{\ib,l}$ in degree 0 and
$\sigma_{\ib,k}$, $\pi_{\ib,1}$ in degree 1. 
This filtration is a nonnegative increasing $\kb$-algebra filtration.
Each term of the filtration is a graded subspace
of $\th\Rb(\Gamma)_{\l,\nu}$. Therefore the associated graded
$\kb$-algebra $\gr\,\th\Rb(\Gamma)_{\l,\nu}$ is bigraded and the symbol map
preserves the grading.

Now, the {\it Nil Hecke algebra} of type $\B_m$ is the $\kb$-algebra
$\th\Nb\Hb_m$ generated by the elements 
$\bar\pi_1,\bar\sigma_1,\bar\sigma_2,\dots,\bar\sigma_{m-1}$
with the relations
$$\gathered
\bar\sigma_k\bar\sigma_{k'}=\bar\sigma_{k'}\bar\sigma_k\ \roman{if}\ |k-k'|>1,
\quad
\bar\pi_1\bar\sigma_k=\bar\sigma_{k}\bar\pi_1\ \roman{if}\ k\neq 1,
\quad
(\bar\pi_1\bar\sigma_1)^2=(\bar\sigma_1\bar\pi_1)^2,
\vspace{2mm}
\bar\sigma_{k+1}\bar\sigma_{k}\bar\sigma_{k+1}
=\bar\sigma_{k}\bar\sigma_{k+1}\bar\sigma_{k},\quad
\bar\pi_1^2=\bar\sigma_k^2=0.
\endgathered
$$
We can form the semidirect product
$\th\Fb_\nu\rtimes\th\Nb\Hb_m$, which is generated by 
$1_\ib$, $\bar\varkappa_l$, $\bar\pi_1,\bar\sigma_k$
with the  relations above and
$$\gathered
\bar\sigma_k\bar\varkappa_{l}=\bar\varkappa_{s_k(l)}\bar\sigma_k,
\quad
\bar\pi_1\bar\varkappa_l=\bar\varkappa_{\eps_1(l)}\bar\pi_1\ ,
\quad
\bar\varkappa_{l}\bar\varkappa_{l'}
=\bar\varkappa_{l'}\bar\varkappa_{l'}.
\endgathered
$$
We have a surjective $\kb$-algebra morphism
$$\th\Fb_\nu\rtimes\th\Nb\Hb_m\to\gr\,\th\Rb(\Gamma)_{\l,\nu},\quad
1_\ib\mapsto 1_{\ib},\quad
\bar\varkappa_l\mapsto\varkappa_{l},\quad
\bar\pi_1\mapsto\pi_{1},\quad
\bar\sigma_{k}\mapsto\sigma_{k}
.\leqno(5.3)$$
Thus the elements $\sigma_{\dot w}$ with $w\in W_m$ generate 
$\th\Rb(\Gamma)_{\l,\nu}$
as a $\th\Fb_\nu$-module. We must prove that they yield indeed a basis of
$\th\Rb(\Gamma)_{\l,\nu}$. 
This is rather clear, since the images  of these elements in
$\End(\th\Fb_\nu)$ under the polynomial representation 
are independent over $\th\Fb_\nu$ (by Galois theory).
Therefore the map (5.3) is invertible.
The last claim is now clear, because the element $\sigma_{\dot w}$ has
the same degree as its symbol and if $\dot w$, $\ddot w$ are two reduced 
decomposition of $w$ then $\sigma_{\dot w}$ and $\sigma_{\ddot w}$ 
have the same symbol.

\qed

\vskip3mm

Let $\th\Fb'_\nu=\bigoplus_\ib\th\Fb'_\ib,$ where
$\th\Fb'_\ib$ is the localization of the ring $\th\Fb_\ib$ with respect
to the multiplicative system generated by 
$$\{\varkappa_{\ib,l}\pm \varkappa_{\ib,l'};\,
1\leqslant l\neq l'\leqslant m\}\cup
\{\varkappa_{\ib,l};\,l=1,2,\dots,m\}.$$ 

\proclaim{5.6.~Corollary}
The polynomial representation of $\th\Rb(\Gamma)_{\l,\nu}$ 
on $\th\Fb_{\nu}$ is faithful.
The inclusion of
$\th\Rb(\Gamma)_{\l,\nu}$ into $\End(\th\Fb_{\nu})$ 
yields an isomorphism of $\th\Fb'_\nu$-algebras
from
$\th\Fb'_\nu\otimes_{\th\Fb_\nu}\th\Rb(\Gamma)_{\l,\nu}$ to 
$\th\Fb'_\nu\rtimes W_m,$
such that
for each $\ib$ and each $l=1,2,\dots,m$, $k=1,2,\dots, m-1$ we have
$$\aligned
&1_{\ib}\mapsto 1_{\ib},\vspace{2mm}
&\varkappa_{\ib,l}\mapsto \varkappa_l1_{\ib},\vspace{2mm}
&\pi_{\ib,1}\mapsto\varkappa_{0}^{\lambda_{i_1}}\eps_11_\ib,
\vspace{2mm}
&\sigma_{\ib,k}\mapsto
\cases
(\varkappa_k-\varkappa_{k+1})^{-1}(s_k-1)1_\ib
& \text{if}\
i_k=i_{k+1},\vspace{2mm}
(\varkappa_k-\varkappa_{k+1})^{h_{i_{k+1},i_k}}s_k1_\ib\hfill&
\text{if}\ i_k\neq i_{k+1}.
\endcases
\endaligned\leqno(5.4)$$
\endproclaim

\vskip3mm

Restricting the $\th\Fb_\nu$-action on $\th\Rb(\Gamma)_{\l,\nu}$ to the
subalgebra $\th\Sb_\nu$ of $\th\Fb_\nu$ we get a structure of graded
$\th\Sb_\nu$-algebra on $\th\Rb(\Gamma)_{\l,\nu}$. 

\proclaim{5.7.~Proposition} 
(a) $\th\Sb_\nu$ is isomorphic to the center of $\th\Rb(\Gamma)_{\l,\nu}$.
\vskip1mm

(b) $\th\Rb(\Gamma)_{\l,\nu}$ is a free graded module over $\th\Sb_\nu$ of rank
$(2^mm!)^2$.
\endproclaim

\noindent{\sl Proof :}
First we prove $(a)$. Recall that
$$\th\Sb_\nu=
\kb[\chi_1,\chi_2,\dots,\chi_m]^{W_\nu}
=\Bigl(\bigoplus_\ib\kb[\varkappa_{1},\varkappa_{2},\dots,
\varkappa_{m}]1_\ib\Bigr)^{W_m}.$$
Given a sequence $\ib$ in $\th\! I^\nu$ the assignment
$x\mapsto x 1_\ib$ embeds $\th\Sb_\nu$
as a central subalgebra of
$\th\Rb(\Gamma)_{\l,\nu}$.
We must check that this map surjects onto the center of
$\th\Rb(\Gamma)_{\l,\nu}$.
This follows from Corollary 5.6.
Part $(b)$ follows from $(a)$ and Proposition 5.5.

\qed

\vskip3mm

In Section 9 we'll prove the following theorem.

\proclaim{5.8.~Theorem} For any $\nu\in\th\NN I$, $\l\in\NN I$
there is an unique graded
$\th\Sb_\Vb$-algebra isomorphism
$$\Psi:\th\Rb(\Gamma)_{\l,\nu}\to\th\Zb^\delta_{\Lbb,\Vb}$$
which intertwines the representations of 
$\th\Rb(\Gamma)_{\l,\nu}$ and $\th\Zb^\delta_{\Lbb,\Vb}$ on $\th\Fb_\nu$.
\endproclaim

\vskip3mm

\subhead 5.9.~Examples\endsubhead
$(a)$ If $m=0$ then
$\th\Rb(\Gamma)_{\l,\nu}=\kb$ by definition.

\vskip1mm

$(b)$
Assume that $m=1$. Fix a vertex $i$ in $I$ and set $\nu=i+\theta(i)$.
We have $\th\!I^\nu=\{\ib,\theta(\ib)\}$ with $\ib=i\theta(i)$ 
and $\theta(\ib)=\theta(i)i$. 
We have
$$\gathered
\th\Rb(\Gamma)_{\l,\nu}=
(\kb[\varkappa_1]\oplus\pi_1\kb[\varkappa_1])1_\ib\oplus
(\kb[\varkappa_{1}]\oplus\pi_{1}\kb[\varkappa_{1}])1_{\theta(\ib)},
\vspace{2mm}
\pi_1\varkappa_11_\ib=-\varkappa_{1}\pi_11_\ib,
\quad\pi_1\varkappa_11_{\theta(\ib)}=
-\varkappa_1\pi_11_{\theta(\ib)},\vspace{2mm}
\pi_1^21_{\theta(\ib)}=
(-1)^{\l_{\theta(i)}}\varkappa_1^{\l_i+\l_{\theta(i)}}1_{\theta(\ib)},\quad
\pi_1^21_\ib=(-1)^{\l_{i}}\varkappa_1^{\l_i+\l_{\theta(i)}}1_\ib.
\endgathered$$
The inclusion $\th\Sb_\nu\subset\th\Rb(\Gamma)_{\l,\nu}$
is given by 
$$\kb[\chi]\to\th\Rb(\Gamma)_{\l,\nu},\quad
\chi\mapsto(\varkappa_11_\ib,0,-\varkappa_11_{\theta(\ib)},0).$$

\vskip3mm

\vskip2cm

\head 6.~Affine Hecke algebras of type B\endhead

\subhead 6.1.~Affine Hecke algebras of type B\endsubhead 
Given a connected reductive group $G$ we call {\it affine Hecke algebra of $G$}
the Hecke algebra of the extended affine Weyl group
$W\ltimes P$ where $W$ is the Weyl group of $(G,T)$,
$P$ is the group of characters of $T$, and $T$ is a maximal torus of $G$.
Fix $p,q$ in $\kb^\times$. 
For any integer $m\geqslant 0$ we define
the affine Hecke algebra $\Hb_m$ of type $\B_m$ to be
the affine Hecke algebra of $SO(2m+1)$.
It admits the following presentation, see e.g., \cite{Mc}. 
If $m>0$ then $\Hb_m$ is the $\kb$-algebra generated by
$$T_k,\quad X_l^{\pm 1},\quad 
k=0,1,\dots,m-1,\quad l=1,2,\dots,m$$ satisfying the
following defining relations :

\vskip2mm
\itemitem{$(a)$} $X_lX_{l'}=X_{l'}X_l$,

\vskip2mm

\itemitem{$(b)$} $(T_0T_1)^2=(T_1T_0)^2$,
$T_kT_{k-1}T_k=T_{k-1}T_kT_{k-1}$ if $k\neq 0,1$, and
$T_kT_{k'}=T_{k'}T_k$ if $|k-k'|\neq 1$,

\vskip2mm

\itemitem{$(c)$} $T_0X_1^{-1}T_0=X_1$,
$T_kX_kT_k=X_{k+1}$ if $k\neq 0$, and $T_kX_l=X_lT_k$ if $l\neq
k,k+1$,

\vskip2mm

\itemitem{$(d)$}
$(T_k-p)(T_k+p^{-1})=0$ if $k\neq 0$, and $(T_0-q)(T_0+q^{-1})=0$.

\vskip3mm

\noindent If $m=0$ then $\Hb_0=\kb$, the trivial $\kb$-algebra.
Note that $\Hb_1$ is the $\kb$-algebra generated by $T_0$, $X_1^{\pm 1}$
with the defining relations
$$T_0X_1^{-1}T_0=X_1,\quad
(T_0-q)(T_0+q^{-1})=0.$$

\vskip3mm

\subhead 6.2.~Intertwiners and blocks of $\Hb_m$\endsubhead 
We define
$$\gathered
\Ab=\kb[X_1^{\pm 1},X_2^{\pm 1},\dots,X_m^{\pm 1}],\quad
\Ab'=\Ab[\Sigma^{-1}],\quad
\Hb'_m=\Ab'\otimes_\Ab\Hb_m,
\endgathered$$ 
where $\Sigma$ is the multiplicative set generated by
$$1-X_lX_{l'}^{\pm 1},\quad 1-p^2X_lX_{l'}^{\pm 1},\quad
1-X_l^{2},\quad 1-q^2X_l^{\pm 2},\quad l\neq l'.$$
For $k=0,\dots,m-1$ the intertwiner
$\varphi_k$ in $\Hb'_m$ is given by the
following formulas
$$
\matrix
&\displaystyle{\varphi_k-1={X_k-X_{k+1}\over pX_k-p^{-1}X_{k+1}}\,(T_k-p)}
\hfill
&\roman{if}\ k\neq 0, \vspace{2mm}
&\displaystyle{\varphi_0-1={X_{1}^{-2}-1\over
qX_1^{-2}-q^{-1}}\,(T_0-q).}\hfill&
\endmatrix\leqno(6.1)$$
The group $W_m$ acts on $\Ab'$ as follows
$$\gathered(s_ka)(X_1,\dots,X_m)=a(X_1,\dots,X_{k+1},X_k,\dots,X_m),
\vspace{2mm}
(\eps_1a)(X_1,\dots,X_m)=a(X_1^{-1},X_2,\dots,X_m).\endgathered$$
There is an isomorphism of $\Ab'$-algebras
$$\Ab'\rtimes W_m\to\Hb'_m,\quad
s_k\mapsto\varphi_k,\quad \eps_1\mapsto\varphi_0,\quad k\neq 0.$$
The semi-direct product
group $\ZZ\rtimes\ZZ_2=\ZZ\rtimes\{-1,1\}$ acts on $\kb^\times$ by
$(n,\eps):z\mapsto z^\eps p^{2n}$. 
Given a $\ZZ\rtimes\ZZ_2$-invariant subset $I$ of $\kb^\times$ we denote by
$\Hb_m\text{-}\Modb_I$ the category of all finitely generated $\Hb_m$-modules
such that the action of $X_1,X_2,\dots,X_m$  is locally finite and
all the eigenvalues belong to $I$. 
We associate to the set $I$ the quiver $\Gamma$ with set of
vertices $I$ and with one arrow $p^2i\to i$ whenever $i$ lies in $I$.
We equip $\Gamma$ with an involution $\theta$ such that
$\theta(i)=i^{-1}$ for each vertex $i$ and such that $\theta$ takes
the arrow $p^2i\to i$ to the arrow $i^{-1}\to p^{-2}i^{-1}$. 
{\it We'll assume that the set $I$
does not contain $1$ nor $-1$ and that $p\neq 1,-1$}. 
Thus the involution $\theta$ has no
fixed points and no arrow may join a vertex of $\Gamma$ to itself.

\subhead 6.3.~Remark\endsubhead 
We may assume that
either $I$ is a $\ZZ$-orbit or $I$ contains at least one
of $\pm q$, see the discussion in \cite{EK1}. 
Thus, we can assume that one of the following two cases holds : 

\vskip2mm

\itemitem{$(a)$} $I$ is a $\ZZ$-orbit which does not contain $1$, $-1$, $q$, $-q$.
So either $I=\{p^n;\,n\in\ZZ_{\roman{odd}}\}$ or
$I=\{-p^n;\,n\in\ZZ_{\roman{odd}}\}$.
Then $\Gamma$ is of type $\A_\infty$ if $p$ has infinite order and
$\Gamma$ is of type $\A^{(1)}_r$ if $p^2$ is a primitive $r$-th root
of unity.

\vskip2mm

\itemitem{$(b)$} 
$q\in I$ (the case $-q\in I$ is similar)
and $-1,1\notin I$. Then we have
$I=\{qp^{2n};\,n\in\ZZ\}\cup\{q^{-1}p^{2n};\,n\in\ZZ\}$
with $q^2\neq p^{4n}$ for all $n\in\ZZ$. Thus
$\Gamma$ is of type $\A_\infty$, $\A_\infty\times \A_\infty$,
$\A^{(1)}_r$, or $\A^{(1)}_r\times \A^{(1)}_r$.

\vskip3mm

\subhead 6.4.~$\Hb_m$-modules versus $\th\Rb_{m}$-modules\endsubhead
Given an element $\l$ of $\NN I$ we define the graded $\kb$-algebra
$$\th\Rb_{I,\l,m}=\bigoplus_\nu\th\Rb_{I,\l,\nu},\quad
\th\Rb_{I,\l,\nu}=\th\Rb(\Gamma)_{\l,\nu},\quad
\th\! I^m=\coprod_\nu\th\! I^\nu,$$
where $\nu$ runs over the set of all dimension vectors in $\th\NN I$ such that
$|\nu|=2m$.  
When there is no risk of confusion we abbreviate
$\th\Rb_{m}=\th\Rb_{I,\l,m}$ and $\th\Rb_{\nu}=\th\Rb_{I,\l,\nu}$.
Note that the $\kb$-algebra $\th\Rb_{m}$ may not have 1, because the set
$I$ may be infinite, and that
$\th\Rb_0=\kb$ as a graded $\kb$-algebra.
From now on, unless specified otherwise we'll set
$$\l=\sum_{i}i,\quad i\in I\cap\{q,-q\}. 
\leqno(6.2)$$ Given sequences
$$\ib=(i_{1-m},\dots,i_{m-1},i_{m}),
\quad
\ib'=(i'_1,\dots,i'_{m'-1},i'_{m'}),$$
we define a sequence $\theta(\ib')\ib\ib'$
as follows
$$
\theta(\ib')\ib\ib'=(\theta(i'_{m'}),\dots,\theta(i'_1),i_{1-m},
\dots,i_{m},i'_1,\dots,i'_{m'}).
$$
Let $\nu,$ $\nu'$ be dimension vectors in $\th\NN I$ and $\NN I$ 
respectively such that
$|\nu|=2m$, $|\nu'|=m'$, and  $m+m'=m''$. 
We define an idempotent in $\th\Rb_{m''}$ by
$$1_{\nu,\nu'}=\sum_{\ib,\ib'}1_{\theta(\ib')\ib\ib'},\quad
\ib\in\th\!I^\nu,\ \ib'\in I^{\nu'}.$$
For $\nu'_1,\nu'_2,\dots,\nu'_r$ in $\NN I$ 
we define $1_{\nu,\nu'_1,\dots,\nu'_r}$
in the same way.
Finally, for any graded $\th\Rb_{m''}$-module $M$ we set
$$1_{m,\nu'}M=\bigoplus_{\ib,\ib'}
1_{\theta(\ib')\ib\ib'}M,\quad\ib\in\th\!I^{m},\ \ib'\in I^{\nu'}.
\leqno(6.3)$$
If $M$ is a right graded $\th\Rb_{m''}$-module 
we define $M1_{m,\nu'}$ in the same way.

Next, let $\th\Rb_{m}\text{-}\Modb_0$ be the category
of all finitely generated (non-graded)
$\th\Rb_{m}$-modules such that the elements 
$\varkappa_1,\varkappa_2,\dots,\varkappa_m$ 
act locally nilpotently.
Let 
$$\th\Rb_{m}\text{-}\fModb_0\subset\th\Rb_{m}\text{-}\Modb_0,\quad
\Hb_m\text{-}\fModb_I\subset\Hb_m\text{-}\Modb_I$$
be the full subcategories of finite dimensional modules.

Fix a formal series $f(\varkappa)$ in $\kb[[\varkappa]]$ such that
$f(\varkappa)=1+\varkappa$ modulo $(\varkappa^2)$.

\proclaim{6.5.~Theorem}We have an equivalence of categories 
$$\th\Rb_{m}\text{-}\Modb_0\to\Hb_m\text{-}\Modb_I,\quad M\mapsto M$$
which is given by
\vskip1mm
\itemitem{$(a)$}
$X_l$ acts on $1_\ib M$ by $i_l^{-1}f(\varkappa_l)$ for $l=1,2,\dots,m$, 
\vskip1mm
\itemitem{$(b)$}$T_k$ acts on $1_\ib M$ 
as follows for $k=1,2,\dots, m-1,$
$$\matrix
&\displaystyle
{{(pf(\varkappa_k)-p^{-1}f(\varkappa_{k+1}))(\varkappa_k-\varkappa_{k+1})
\over f(\varkappa_k)-f(\varkappa_{k+1})}\sigma_k+p}
\hfill&
\text{if}\ i_{k+1}=i_{k},\hfill\vspace{2mm}
&\displaystyle{
{f(\varkappa_k)-f(\varkappa_{k+1})\over
(p^{-1}f(\varkappa_k)-pf(\varkappa_{k+1}))
({\varkappa_k}-{\varkappa_{k+1})}}\sigma_k+{(p^{-2}-1)f(\varkappa_{k+1})
\over pf(\varkappa_k)-p^{-1}f(\varkappa_{k+1})}}
&\text{if}\  i_{k+1}=p^2i_{k},\hfill\vspace{2mm}
&\displaystyle{
{pi_{k}f(\varkappa_k)-p^{-1}i_{k+1}f(\varkappa_{k+1})\over
i_{k}f(\varkappa_k)-i_{k+1}f(\varkappa_{k+1})}\sigma_k
+{(p^{-1}-p)i_{k}f(\varkappa_{k+1})\over
i_{k+1}f(\varkappa_k)-i_{k}f(\varkappa_{k+1})}}
\hfill&\text{if}\  i_{k+1}\neq i_{k},p^2i_{k} ,\hfill
\endmatrix$$
\vskip1mm
\itemitem{$(c)$}$T_0$ acts on $1_\ib M$ 
as follows 
$$\matrix
\displaystyle{{f(\varkappa_1)^2-1\over
(q^{-1}-qf(\varkappa_1)^2)\varkappa_1}\pi_1+{(q^{-2}-1)f(\varkappa_1)^2\over
q-q^{-1}f(\varkappa_1)^2}}
\hfill&\text{if}\ i_1=\pm q,\vspace{2mm}
\displaystyle{{q-q^{-1}i_1^{2}f(\varkappa_1)^2\over
1-i_1^{2}f(\varkappa_1)^2}\pi_1+{q-q^{-1}\over
1-i_1^{2}f(\varkappa_1)^{-2}}}\hfill&\text{if}\   i_1\neq\pm q.
\endmatrix$$
\endproclaim

\subhead 6.6.~Remark\endsubhead
The first case in $(c)$ does not occur if 
$q=q^{-1}$ because $\theta$ has no fixed points in $I$.
In the second case we have $i_1^2\neq 1$ for the same reason.
Note also that $(f(\varkappa)-1)/\varkappa$ is a formal series in
$\kb[[\varkappa]]$, and that 
$(f(\varkappa_1)-f(\varkappa_2))/(\varkappa_1-\varkappa_2)$ 
is an invertible formal series in $\kb[[\varkappa_1-\varkappa_2]]$.
Finally, recall that $p^2\neq 1$.

\vskip3mm

\noindent{\sl Proof :}
First, recall that $\pm 1\notin I$ and that $p\neq \pm 1$.
Observe also that (6.2) yields $$
\gathered
i_1=\pm q\iff
\l_{i_1}=1,\vspace{2mm} 
i_1\neq\pm q\iff
\l_{i_1}=0.
\endgathered$$
The functor above is well defined by formulas (5.4) and (6.1).
Let $g$ be the inverse of $f$, i.e., 
$g(X)$ is the unique formal series in $\kb[[X-1]]$
such that $gf(\varkappa)=\varkappa$.
For instance, we may choose
$$f(\varkappa)=1+\varkappa,\quad g(X)=X-1.$$
A quasi-inverse functor $\Hb_m\text{-}\Modb_I\to\th\Rb_{m}\text{-}\Modb_0$
such that $M\mapsto M$ is given by the following rules
\vskip1mm
\itemitem{$(a)$} $1_\ib M=\{m\in M;\,(i_lX_l-1)^r m=0,\,r\gg 0\}$, 
\vskip1mm
\itemitem{$(b)$} $\varkappa_l$ acts on $1_\ib M$ by $g(i_lX_l)$ 
for $l=1,2,\dots,m$,
\vskip1mm
\itemitem{$(c)$} $\sigma_k$ acts on $1_\ib M$ 
as follows for $k=1,2,\dots, m-1,$
$$\matrix
&\displaystyle
{{X_k-X_{k+1}
\over (pX_k-p^{-1}X_{k+1})(g(i_kX_k)-g(i_kX_{k+1}))}(T_k-p)}
\hfill&
\roman{if}\ i_{k+1}=i_{k},\hfill\vspace{2mm}
&\displaystyle{
{g(i_{k+1}X_k)-g(i_{k}X_{k+1})\over
pX_k-p^{-1}X_{k+1}}\Bigl((X_k-X_{k+1})T_k+(p-p^{-1})X_{k+1}\Bigr)
}
&\roman{if}\  i_{k+1}=p^2i_{k},\hfill\vspace{2mm}
&\displaystyle{
T_k{X_k-X_{k+1}\over p^{-1}X_k-pX_{k+1}}+(p-p^{-1}){X_{k+1}\over p^{-1}X_k-pX_{k+1}}}
\hfill&\roman{if}\  i_{k+1}\neq i_{k},p^2i_{k} ,\hfill
\endmatrix$$
\vskip1mm
\itemitem{$(d)$} $\pi_1$ acts on $1_\ib M$ 
as follows 
$$\matrix
\displaystyle{{g(i_0X_1)\over q^{-1}-qX_1^{-2}}\Bigl((X_1^{-2}-1)T_0
+q-q^{-1}\Bigr)}
\hfill&\roman{if}\ i_1=\pm q,\vspace{2mm}
\displaystyle{T_0{X_1^{-2}-1\over
q^{-1}X_1^{-2}-q}+{q-q^{-1}\over
q^{-1}X_1^{-2}-q}}\hfill&\roman{if}\   i_1\neq\pm q.
\endmatrix$$

\noindent 
Note that $g(X)/(X-1)$ is a formal series in $\kb[[X-1]]$,
and that $(X_1-X_2)/(g(X_1)-g(X_2))$ is an invertible
formal series in $\kb[[X_1-X_2]]$.

\qed

\vskip3mm

\proclaim{6.7.~Corollary} There is an equivalence of categories
$$\Psi:\th\Rb_{m}\text{-}\fModb_0\to\Hb_m\text{-}\fModb_I,\quad
M\mapsto M.$$
\endproclaim

\subhead 6.8.~Example\endsubhead
Let $m=1$.
Using  Example 5.9$(b)$ it is easy to check that the 1-dimensional
$\th\Rb_m$-modules are labelled by $\{i\in I\,;\,\l_i+\l_{\theta(i)}\neq 0\}$,
and that the irreducible 2-dimensional $\th\Rb_m$-modules are labelled by 
$\{i\in I\,;\,\l_i+\l_{\theta(i)}=0\}/\theta$.
Further we have
$\l_i+\l_{\theta(i)}\neq 0$ iff $i=\pm q^{-1}$ or $\pm q$.
On the other hand
the 1-dimensional objects in $\Hb_m\text{-}\Modb_I$ are given by
\vskip2mm
\itemitem{$(a)$} $X_1=i^{-1}$, $T_0=q$, $i\in I\cap\{\pm q^{-1}\}$,
\vskip2mm
\itemitem{$(b)$} $X_1=i^{-1}$, $T_0=-q^{-1}$, $i\in I\cap\{\pm q\}$,
\vskip2mm
\noindent
and the irreducible 2-dimensional objects in $\Hb_m\text{-}\Modb_I$ are given by
\vskip2mm
\itemitem{$(c)$} $X_1=\pmatrix i&0\cr 0&i^{-1}\endpmatrix$, 
$T_0=\pmatrix -i^2a/b&a^2-b^2/i^2\cr -i^2/b^2&a/b\endpmatrix$
with $a=q-q^{-1}$, $b=1-i^2$, and $i\neq \pm q^{-1},\pm q.$
\vskip2mm

\noindent
Therefore Theorem 6.5 is obvious in this case.

\vskip3mm

\subhead 6.9.~Induction and restriction of $\Hb_m$-modules\endsubhead
For $i\in I$ we define  functors
$$\gathered
E_i:\Hb_m\text{-}\fModb_I\to\Hb_{m-1}\text{-}\fModb_I,\vspace{2mm}
F_i:\Hb_m\text{-}\fModb_I\to\Hb_{m+1}\text{-}\fModb_I,
\endgathered\leqno(6.4)$$
where $E_iM\subset M$ is the generalized 
$i^{-1}$-eigenspace  of the $X_m$-action,
and where
$$F_iM=\Ind_{\Hb_m\otimes\kb[X_{m+1}^{\pm 1}]}^{\Hb_{m+1}}
(M\otimes\kb_{i}).$$
Here $\kb_{i}$ is the 1-dimensional representation of 
$\kb[X_{m+1}^{\pm 1}]$ defined by
$X_{m+1}\mapsto i^{-1}$. 

\vskip3mm

\subhead 6.10.~Remark\endsubhead
The results in Section 6 hold true if  
$\kb$ is any field of characteristic $\neq 2$. Indeed, set
$f(\varkappa)=1+\varkappa$ and $g(X)=X-1$. Then we must check
that the formulas for $T_k$, $T_0$ and that the formulas for $\sigma_k$, $\pi_1$
still make sense. This is straightforward for all cases, except for the first formula for $\pi_1$.
Here one needs that $i_1X_1+1$ is invertible, which holds true if the characteristic is not 2.

\vskip2cm

\head 7.~Global bases of $\fb$ and projective graded modules of KLR algebras
\endhead

This section is a reminder on KLR algebras. Most of the results here are due 
to \cite{KL}. Although we are essentially concerned by KLR algebras of type 
A, everything here holds true in any type. 

\subhead 7.1.~Definition of the graded $\kb$-algebra $\Rb_m$\endsubhead 
Fix a $\ZZ\rtimes\ZZ_2$-invariant subset $I\subset\kb^\times$ as in Section 6.2. 
Let $\Gamma$ be the corresponding quiver.
For each integer $m\geqslant 0$ we put
$$\Rb(I)_m=\bigoplus_\nu\Rb(I)_\nu,\quad
\Rb(I)_\nu=\Rb(\Gamma)_{\nu},$$
where $\nu$ runs over the set of all dimension vectors in $\NN I$ such that
$|\nu|=m$.  Here $\Rb(\Gamma)_{\nu}$ is the graded $\kb$-algebra introduced in Section 1.3.
When there is no risk of confusion we'll abbreviate $\Rb_m=\Rb(I)_m$.
Let $Q_{i,j}(u,v)$ be as in (5.1). If $m>0$ the graded $\kb$-algebra $\Rb_m$
is generated by elements $1_\ib$, $\varkappa_{\ib,l}$, $\sigma_{\ib,k}$ with
$\ib\in I^m$, $l=1,2,\dots, m$ and $k=1,2,\dots, m-1$ 
satisfying the following defining relations 

\vskip2mm

\itemitem{$(a)$}
$1_\ib\,1_{\ib'}=\delta_{\ib,\ib'}1_\ib$,\quad
$\sigma_{\ib,k}=1_{s_k\ib}\sigma_{\ib,k}1_\ib$,\quad
$\varkappa_{\ib,l}=1_{\ib}\varkappa_{\ib,l}1_\ib$,\quad

\vskip2mm

\itemitem{$(b)$}
$\varkappa_{l}\varkappa_{l'}=\varkappa_{l'}\varkappa_{l}$,

\vskip2mm

\itemitem{$(c)$}
$\sigma_{k}^21_\ib=
Q_{i_k,i_{k+1}}(\varkappa_{k+1},\varkappa_{k})1_\ib$,

\vskip2mm

\itemitem{$(d)$}
$\sigma_{k}\sigma_{k'}=\sigma_{k'}\sigma_{k}$ if $|k-k'|>1$,

\vskip2mm

\itemitem{$(e)$}
$(\sigma_{k+1}\sigma_{k}\sigma_{k+1}-
\sigma_{k}\sigma_{k+1}\sigma_{k})1_\ib=$
$$=\cases
{\ds{Q_{i_k,i_{k+1}}(\varkappa_{k+1},\varkappa_{k})
-Q_{i_k,i_{k+1}}(\varkappa_{k+1},\varkappa_{k+2})\over
\varkappa_{k}-\varkappa_{k+2}}1_\ib}&
\roman{if}\ i_k=i_{k+2},
\vspace{2mm}
0&\roman{else},\endcases$$

\vskip2mm

\itemitem{$(f)$}
$(\sigma_{k}\varkappa_{k'}-\varkappa_{s_k(k')}\sigma_{k})1_\ib
=\cases -1_\ib&\ \roman{if}\  k'=k,\, i_k=i_{k+1}, \cr 1_\ib&\
\roman{if}\  k'=k+1,\, i_k=i_{k+1},\cr 0&\ \roman{else}.
\endcases$

\vskip2mm

The grading on $\Rb_m$ is given by the following rules : $1_\ib$
has the degree 0, $\varkappa_{\ib,l}$ has the degree 2,
and $\sigma_{\ib,k}$ has the degree $-i_k\cdot i_{k+1}$.
Given any element $a$ in $1_\ib\Rb_{m}1_{\ib'}$ we write
$\varkappa_ka=\varkappa_{\ib,k}a$, 
$a\varkappa_k=a\varkappa_{\ib',k}$, 
etc.
Note that the $\kb$-algebra
$\Rb_m$ may not have 1, because the set $I$ may be infinite.
If $m=0$ we have $\Rb_m=\kb$ as a graded $\kb$-algebra.

Let $\omega$ be the unique anti-involution of 
the graded $\kb$-algebra $\Rb_m$ given by 
$$\omega:\ 1_\ib,\,\varkappa_l,\,\sigma_k\mapsto
1_\ib, \,\varkappa_l, \,\sigma_k.$$
Note that $Q_{i,j}(u,v)=Q_{j,i}(v,u)$. 
Hence there is an unique involution $\tau$ of 
the graded $\kb$-algebra $\Rb_m$ such that 
$$\tau:\ 1_\ib, \,\varkappa_l, \,\sigma_k
\mapsto1_{w_m(\ib)}, \,\varkappa_{m+1-l}, \,-\sigma_{m-k},$$
where
$w_m$ is the longest element in $\Sen_m$. Finally, we have
$Q_{i,j}(u,v)=Q_{\theta(i),\theta(j)}(-u,-v)$. 
Hence there is an unique involution 
$$\iota:\ 1_\ib, \,\varkappa_l, \,\sigma_k\mapsto 
1_{\theta(\ib)},\,-\varkappa_{l}, \,-\sigma_{k}.$$ 
We define 
$$\kappa=\iota\circ\tau=\tau\circ\iota.\leqno(7.1)$$

\vskip3mm

\subhead 7.2.~The Grothendieck groups of $\Rb_m$\endsubhead 
The graded $\kb$-algebra
$\Rb_m$ is finite dimensional over its center, a
commutative graded $\kb$-subalgebra.
Therefore any
simple object of $\Rb_m$-$\modb$ is finite-dimensional and
there is a finite number of simple modules in $\Rb_m$-$\modb$.
The Abelian group $G(\Rb_m)$ is  free  with a
basis formed by the classes of the simple objects of
$\Rb_m$-$\modb$, see Section 0.2 for the notation. 
The Abelian group $K(\Rb_m)$ is also free, with a
basis formed by the classes of the indecomposable projective
objects. 
Both Abelian groups
are free $\Ac$-modules where $v$ shifts the grading by $1$.
We define
$$\gathered
\Kb_I=\bigoplus_{m\geqslant 0}\Kb_{I,m},\quad \Kb_{I,m}=K(\Rb_m),\vspace{1mm}
\Gb_I=\bigoplus_{m\geqslant 0}\Gb_{I,m},\quad \Gb_{I,m}=G(\Rb_m).
\endgathered$$
Now, fix integers $m,m',m''\geqslant 0$
with $m''=m+m'$.
Given sequences $\ib\in I^m$ and
$\ib'\in I^{m'}$ we write $\ib''=\ib\ib'$.
We'll abbreviate
$$\Rb_{m,m'}=\Rb_{m}\otimes\Rb_{m'}.$$
There is an unique
inclusion of graded $\k$-algebras
$$\gathered
\phi:\Rb_{m,m'}\to\Rb_{m''},\vspace{2mm}
1_\ib\otimes 1_{\ib'}\mapsto1_{\ib''}, \vspace{2mm}
\varkappa_{\ib,l}\otimes 1_{\ib'}\mapsto \varkappa_{\ib'',l},
\vspace{2mm}
1_\ib\otimes \varkappa_{\ib',l}\mapsto \varkappa_{\ib'',m+l},
\vspace{2mm}
\sigma_{\ib,k}\otimes 1_{\ib'}\mapsto \sigma_{\ib'',k},
\vspace{2mm}
1_\ib\otimes \sigma_{\ib',k}\mapsto \sigma_{\ib'',m+k}.
\endgathered\leqno(7.2)$$
This yields a triple of adjoint functors $(\phi_!,\phi^*,\phi_*)$ where
$$\phi^*\,:\,
\Rb_{m''}\text{-}\modb\to
\Rb_m\text{-}\modb\times\Rb_{m'}\text{-}\modb$$ is the restriction and
$\phi_!$, $\phi_*$ are given by
$$\aligned
&\phi_!\,:\,\cases
\Rb_m\text{-}\modb\times\Rb_{m'}\text{-}\modb\to
\Rb_{m''}\text{-}\modb,\vspace{2mm}
(M,M')\mapsto
\Rb_{m''}\otimes_{\Rb_{m,m'}}(M\otimes
M'),\endcases
\vspace{2mm}
&\phi_*\,:\,\cases
\Rb_m\text{-}\modb\times\Rb_{m'}\text{-}\modb\to
\Rb_{m''}\text{-}\modb,\vspace{2mm}
(M,M')\mapsto
\gHom_{\Rb_{m,m'}}(
\Rb_{m''},M\otimes M').\endcases
\endgathered
$$
First, note that the functors $\phi_!$, $\phi^*$, $\phi_*$ commute with the shift of the grading. 
Next, the functor $\phi^*$ is exact
and it takes finite dimensional graded modules to finite dimensional ones.
By \cite{KL, prop.~2.16} the right graded $\Rb_{m,m'}$-module
$\Rb_{m''}$ is free of finite rank.
Thus $\phi_!$ is exact and it 
takes finite dimensional graded modules to finite dimensional ones.
For the same reason the left graded $\Rb_{m,m'}$-module
$\Rb_{m''}$ is free of finite rank.
Thus $\phi_*$ is exact and
it takes finite dimensional graded modules to finite dimensional ones.
Further $\phi_!$ and $\phi^*$ take  projective graded modules to projective ones,
because they are left adjoint to the exact functors $\phi^*$, $\phi_*$ respectively.
To summarize, the functors $\phi_!$, $\phi^*$, $\phi_*$ are exact
and take finite dimensional graded modules to finite dimensional ones,
and the functors $\phi_!$, $\phi^*$ take projective graded modules to projective ones.
Taking the sum over all $m$, $m'$ we get an $\Ac$-bilinear map
$$
\phi_!\,:\,\Kb_{I}\times\Kb_{I}\to \Kb_{I}.$$
In the same way we define also an $\Ac$-linear map
$$
\phi^*\,:\,\Kb_I\to\Kb_{I}\otimes_\Ac\Kb_{I}.
$$ 
From now on, to unburden the notation we may abbreviate
$\Rb=\Rb_m$, hoping it will not create any confusion.
Recall the anti-automorphism $\omega$
from the previous section.
Consider the duality 
$$\Rb\text{-}\proj\to\Rb\text{-}\proj,\quad
P\mapsto P^\sharp=\gHom_{\Rb}(P,\Rb),$$
with the action and the grading given by 
$$(xf)(p)=f(p)\omega(x),\quad(P^\sharp)_d=\Hom_\Rb(P[-d],\Rb).$$
We'll say that $P$ is $\sharp$-selfdual if $P^\sharp=P$.
The duality on $\Rb$-$\projb$ yields an $\Ac$-antilinear map
$$\Kb_I\to\Kb_I,\quad P\mapsto P^\sharp.$$
Set $\Bc=\ZZ((v))$. The $\Ac$-module $\Kb_I$ is equipped with a 
symmetric $\Ac$-bilinear form 
$$\gathered
\Kb_I\times \Kb_I\to \Bc,\quad
(P:Q)
=\gdim(P^\omega\otimes_{\Rb}Q).
\endgathered$$
Here $P^\omega$ is the right graded $\Rb$-module associated with $P$ and the
anti-automorphism $\omega$.
Finally, we equip $\Kb_I\otimes_\Ac\Kb_I$ with the algebra structure such that
$$(P\otimes Q,P'\otimes Q')\mapsto 
v^{-\mu\cdot\nu'}\phi_!(P,P')\otimes\phi_!(Q,Q'),$$
and with the $\Ac$-antilinear map such that $$P\otimes Q\mapsto 
(P\otimes Q)^\sharp=P^\sharp\otimes Q^\sharp.$$
The following is proved in \cite{KL}.

\proclaim{7.3.~Proposition}
The map $\phi_!$ turns $\Kb_I$ into an associative $\Ac$-algebra with 1, 
and it commutes with the duality $\sharp$.
The map $\phi^*$ is an algebra homomorphism
which turns $\Kb_I$ into a coassociative $\Ac$-coalgebra.
\endproclaim

The Cartan pairing is the perfect $\Ac$-bilinear form
$$\Kb_I\times \Gb_I\to \Ac,\quad
\la P:M\ra=\gdim\,\gHom_{\Rb}(P,M).$$
Consider the duality 
$$\Rb\text{-}\fmodb\to\Rb\text{-}\fmodb,\quad
M\mapsto M^\flat=\gHom(M,\kb),$$
where $\kb$ is considered as a graded $\kb$-space homogeneous of degree 0.
The action and the grading are given by 
$$(xf)(m)=f(\omega(x)m),\quad(M^\flat)_d=\Hom(M_{-d},\kb).$$
We'll say that $M$ is $\flat$-selfdual if $M^\flat=M$.

Finally, let $\Bc I^m$ be the free $\Bc$-module with basis $I^m$.
The {\it character} of a finitely generated graded
$\Rb_m$-module $M$ 
is  given by
$$\ch(M)=\sum_\ib\gdim(1_\ib M)\,\ib\in\Bc I^m.$$

\vskip3mm

\subhead 7.4.~The projective graded $\Rb_{m}$-module $\Rb_\yb$\endsubhead 
Fix  $\nu\in\NN I$ with $|\nu|=m$.
For $\yb=(\ib,\ab)$ in $Y^\nu$ we define
an object $\Rb_\yb$ in $\Rb_{m}$-$\proj$ as follows.

\vskip2mm

\itemitem{$\bullet$}
If $\ib=i^m$, $i\in I$, and $\ab=m$ then 
we set
$\Rb_{\yb}=\Fb_\nu[\ell_m]$. As a left graded $\Rb_{\nu}$-modules we
have a canonical isomorphism
$\Rb_{\nu}=\bigoplus_{w\in\Sen_m}\Rb_\yb[2\ell(w)-\ell_m]$.
We choose once for all an idempotent
$1_\yb$ in $\Rb_{m}$ such that
$\Rb_\yb=(\Rb_{m} 1_\yb)[\ell_m].$

\vskip1mm

\itemitem{$\bullet$}
If $\ib=(i_1,\dots, i_k)$ and $\ab=(a_1,\dots a_k)$ we define the
idempotent $1_{\yb}$ as the image of the element
$\bigotimes_{l=1}^k1_{(i_l)^{a_l},a_l}$ by the inclusion of graded
$\k$-algebras
$\bigotimes_{l=1}^k\Rb_{a_li_l}\!\!\subset\Rb_{\nu}$ in (7.2).
Then we set $\Rb_{\yb}=(\Rb_{\nu}1_{\yb})[\ell_\ab].$
\vskip2mm

\noindent The graded module $\Rb_{\yb}$ satisfies the following
properties.

\vskip2mm

\itemitem{$\bullet$}
Let $\ib'\in I^\nu$ be the sequence obtained by expanding the pair
$\ibt=(\ib,\ab)$. We have the following formula in
$\Rb_m$-$\proj$
$$\Rb_{\ib'}=\Rb 1_{\ib'}=
\bigoplus_{w\in\Sen_\ab}\Rb_{\yb}[2\ell(w)-\ell_\ab]
=:\la\ab\ra!\,\Rb_{\yb}.$$
As a consequence, since the graded module
$\Rb_{\ib'}$ is $\sharp$-selfdual we get 
$$\Rb_\yb[d]^\sharp=\Rb_\yb[-d],\quad\forall d\in\ZZ.$$


\itemitem{$\bullet$}
Given $\yb=(\ib,\ab)\in  Y_{\nu}$ and $\yb'=(\ib',\ab')\in Y_{\nu'}$ we
set $\yb\yb'=(\ib\ib',\ab\ab')$. We have an isomorphism of 
graded $\Rb_{\nu''}$-modules
$\phi_!(\Rb_{\yb}\otimes \Rb_{\yb'})=\Rb_{\yb''}.
$

\vskip3mm

\subhead 7.5.~Examples\endsubhead
For $\ib\in I^\nu$, $|\nu|=m$ and  $\yb\in Y^\nu$,
we define $\Lb_{\ib}=\top(\Rb_{\ib})$ and $\Lb_{\yb}=\top(\Rb_{\yb})$.
Observe that $\Lb_\ib$ is not a simple graded $\Rb_m$-module in general.

\vskip1mm

$(a)$ The graded $\kb$-algebra $\Rb_1$
is generated by elements $1_i$, $\varkappa_{i}$, 
$i\in I$,
satisfying the defining relations 
$1_i\,1_{i'}=\delta_{i,i'}1_i$ and
$\varkappa_{i}=1_{i}\varkappa_{i}1_i$.
Note that $\Rb_i=1_i\Rb_1=\Rb_11_i$ is a graded subalgebra of $\Rb_1$,
and that 
$\Lb_i=\Rb_i/(\varkappa_i)=\kb.$
Further, we have
$\ch(\Rb_{i})\in (1-v^2)^{-1}i$
and $\ch(\Lb_{i})=i$, where
the symbol $(1-v^2)^{-1}$ denotes the infinite sum $\sum_{r\geqslant 0}v^{2r}.$
\vskip1mm
$(b)$ Set $\nu=mi$ and $\yb=(i,m)$.
We'll abbreviate $\Lb_{mi}=\Lb_{i,m}=\Lb_\yb$. 
It is a simple graded $\Rb_{m}$-module. We have
$$
\ch(\Rb_{mi})\in v^{-\ell_m}\ZZ[[v^2]]\,i^m,\quad
\ch(\Lb_{mi})=\la m\ra !\,i^m.$$
The graded $\Rb_{m-1}\otimes\Rb_1$-module $\Lb_{mi}$ 
has a filtration by graded submodules whose associated 
graded is isomorphic to $[m]\,\Lb_{(m-1)i}\otimes\Lb_i$. The socle of
the graded $\Rb_{m-n}\otimes\Rb_n$-module $\Lb_{mi}$ is equal to
$\Lb_{(m-n)i}\otimes\Lb_{ni}$ for each $m\geqslant n\geqslant 0$.
See \cite{KL, ex.~2.2, prop.~3.11} for details.

\subhead 7.6.~Categorification of the global bases of $\fb$
\endsubhead Set $\Kc=\QQ(v)$. Let $\fb$ be the  
$\Kc$-algebra generated by elements
$\theta_i$, $i\in I$, with the defining relations
$$\sum_{a+b=1-i\cdot j}(-1)^a\theta_i^{(a)}\theta_j\theta_i^{(b)}=0,
\quad i\neq j,\quad
\theta_i^{(a)}=\theta_i^a/\la a\ra!,\quad a\geqslant 0.\leqno(7.3)$$
We have a weight decomposition
$\fb=\bigoplus_{\nu\in\NN I}\fb_\nu.$
Let ${}_\Ac\fb$ be the $\Ac$-submodule of $\fb$ generated by
all products of the elements $\theta_i^{(a)}$ with 
$a\in\ZZ_{\geqslant 0}$ and $i\in I$. 
We set $${}_\Ac\fb_\nu={}_\Ac\fb\cap\fb_{\nu}.$$
The element $\theta_i$ lies in ${}_\Ac\fb_i$ for each $i\in I$.
For each pair
$\yb=(\ib,\ab)$ in $Y^\nu$ with $\ib=(i_1,\dots, i_k)$,
$\ab=(a_1,\dots a_k)$ we write
$$\theta_{\yb}=
\theta_{i_1}^{(a_1)}\theta_{i_2}^{(a_2)}\cdots\theta_{i_k}^{(a_k)}.$$
We equip ${}_\Ac\fb$ with the unique $\Ac$-antilinear involution such that
$\bar\theta_i=\theta_i$ for each $i\in I$. We equip the tensor square of
$\fb$ with the $\Kc$-algebra structure such that
$$(x\otimes y)(x'\otimes y')= v^{-\mu\cdot\nu'}xx'\otimes yy',\quad
x\in \fb_\nu, \  x'\in \fb_{\nu'},\ y\in \fb_\mu,\ y'\in \fb_{\mu'}.$$
Consider  the $\Kc$-algebra homomorphism such that
$$r:\fb\to\fb\otimes\fb,\quad \theta_i\mapsto\theta_i\otimes 1+1\otimes \theta_i.$$
The $\Kc$-algebra $\fb$ comes equipped with a 
bilinear form $(\bullet:\bullet)$ which is
uniquely determined by the following conditions
\vskip2mm
\itemitem{$\bullet$}
$(1:1)=1$,
\vskip2mm
\itemitem{$\bullet$}
$(\theta_i:\theta_j)=\delta_{i,j}(1-v^2)^{-1}$ for all $i,j\in I$,
\vskip2mm
\itemitem{$\bullet$}
$(x:yy')=(r(x):y\otimes y')$ for all $x,y,y'$,
\vskip2mm
\itemitem{$\bullet$}
$(xx':y)=(x\otimes x': r(y))$ for all $x,x',y$.
\vskip2mm
\noindent
This bilinear form is symmetric and non-degenerate.
Let $\theta^i\in{}_\Ac\fb^*_i$ be the element dual to
$\theta_i$.
We may regard 
$${}_\Ac\fb^*=\bigoplus_\nu{}_\Ac\fb^*_\nu,\quad
{}_\Ac\fb^*_\nu=\Hom_\Ac({}_\Ac\fb_\nu,\Ac),$$ 
as an $\Ac$-submodule of $\fb$
via the bilinear form $(\bullet:\bullet)$. 
Let $\Gb^\low$ be the {\it canonical basis (=the lower global basis)} of $\fb$. 
It is indeed a $\Ac$-basis of ${}_\Ac\fb$. 
The {\it upper global basis of $\fb$} is the $\Kc$-basis
$\Gb^\up$ which is dual to $\Gb^\low$ with respect to the inner product
$(\bullet:\bullet)$.
We may regard $\Gb^\up$ as a $\Kc$-basis of $\fb^*$.
Let $B(\infty)$ be the set of isomorphism classes
of irreducible (non graded) $\Rb$-modules
such that the elements $\varkappa_1,\varkappa_2,\dots,\varkappa_m$ 
act nilpotently.

\proclaim{7.7.~Theorem} 
(a) There is an unique $\Ac$-algebra isomorphism
$\g:{}_\Ac\fb\to \Kb_I$ which intertwines $r$ and $\phi^*$, and such that
$\g(\theta_{\ibt})=\Rb_{\ibt}$ for each $\ibt.$

\vskip1mm

(b) 
We have $\Gb^\low=\{G^\low(b);\,b\in B(\infty)\}$,
where $\g (G^\low(b))$ is the unique 
$\sharp$-selfdual indecomposable projective graded module
whose top is isomorphic to $b$.
The map $\g$ takes the bilinear form $(\bullet:\bullet)$
and the involution $\bar{\bullet}$
on ${}_\Ac\fb$ to the bilinear form $(\bullet:\bullet)$
and the involution ${\bullet}^\sharp$ on $\Kb_I$.
\vskip1mm

(c) The transpose map $\Gb_I\to {}_\Ac\fb^*$
takes the $\Ac$-basis of $\Gb_I$ of the $\flat$-selfdual simple objects
to $\Gb^\up$.
We have ${}^t\g(\Lb_i)=\theta^i$ for all $i\in I$, and
$\Gb^\up=\{G^\up(b);\,b\in B(\infty)\}$
with $G^\up(b)={}^t\g\,\top\,\g\,G^\low(b)$.
\endproclaim

\noindent{\sl Proof :}
Claim $(a)$, and the second part of $(b)$
are due to \cite{KL, prop.~3.4}.
The first part of $(b)$ is due to \cite{VV} 
(the same result has also been anounced by R. Rouquier).
Part $(c)$ follows form $(b)$. For instance, the last claim in $(c)$
is as proved as follows.
Let $\la\bullet:\bullet\ra$ denote both the Cartan pairing and the canonical
pairing $${}_\Ac\fb\times{}_\Ac\fb^*\to\Ac.$$
Then we have
$$\la\bb':{}^t\g\,\top\,\g(\bb)\ra=
\la\g(\bb'):\top\,\g(\bb)\ra=\delta_{\bb,\bb'}.$$

\qed

\vskip2cm

\head 8.~Global bases of $\th\Vb(\l)$ and projective graded $\th\Rb$-modules
\endhead

Given an integer $m\geqslant 0$ we consider the graded $\kb$-algebra
$\th\Rb_m$ introduced in Sections 5.1, 6.4.

\subhead 8.1.~The Grothendieck groups of $\th\Rb_{m}$\endsubhead 
The graded $\kb$-algebra
$\th\Rb_{m}$ is free of finite type over its center by Proposition 5.7$(b)$.
Therefore any
simple object of $\th\Rb_{m}$-$\modb$ is finite-dimensional and
there is a finite number of isomorphism classes of
simple modules in $\th\Rb_{m}$-$\modb$.
Further, the Abelian group $G(\th\Rb_{m})$ is free with a
basis formed by the classes of the simple objects of
$\th\Rb_{m}$-$\modb$. 
For each $\nu$ the graded $\kb$-algebra
$\th\Rb_{\nu}$ has a graded dimension which lies in
$v^{d}\NN[[v]]$ for some integer $d$.
Therefore the Abelian group $K(\th\Rb_{m})$ is free with a
basis formed by the classes of the indecomposable projective
objects. 
Both $G(\th\Rb)$ and $K(\th\Rb)$
are free $\Ac$-modules where $v$ shifts the grading by $1$.
We consider the following $\Ac$-modules
$$\gathered
\th\Kb_I=\bigoplus_{m\geqslant 0}\th\Kb_{I,m},
\quad\th\Kb_{I,m}=K(\th\Rb_{m}),\vspace{1mm}
\th\!\Gb_I=\bigoplus_{m\geqslant 0}\th\!\Gb_{I,m},\quad
\th\!\Gb_{I,m}=G(\th\Rb_{m}).
\endgathered$$ 
From now on, to unburden the notation we may abbreviate
$\th\Rb=\th\Rb_m$, hoping it will not create any confusion.
For any $M,N$ in $\th\Rb$-$\modb$ we set
$$(M:N)
=\gdim(M^\omega\otimes_{\th\Rb}N),\quad
\la M:N\ra= \gdim\,\gHom_{\th\Rb}(M,N).\leqno(8.1)$$
Here $M^\omega$ is the right graded $\th\Rb$-module associated with $M$ and the
anti-automorphism $\omega$ introduced in Section 5.1.
The Cartan pairing is the perfect $\Ac$-bilinear form 
$$
\th\Kb_I\times \th\!\Gb_I\to \Ac,\quad
(P,M)\mapsto\la P: M\ra,$$
see (8.1).
First, we concentrate on the $\Ac$-module $\th\!\Gb_I$.
Consider the duality 
$$\th\Rb\text{-}\fmodb\to\th\Rb\text{-}\fmodb,\quad
M\mapsto M^\flat=\gHom(M,\kb),$$
with the action and the grading given by 
$$(xf)(m)=f(\omega(x)m),\quad(M^\flat)_d=\Hom(M_{-d},\kb).$$
We'll say that $M$ is $\flat$-selfdual if $M^\flat=M$.
The functor $\flat$ yields an $\Ac$-antilinear map
$$\th\!\Gb_I\to\th\!\Gb_I,\quad M\mapsto M^\flat.$$
We can now define the upper global basis of $\th\!\Gb_I$ as follows.
The proof is given in Section 8.26.

\proclaim{8.2.~Proposition/Definition}
Let $\th\!B(\l)$ be the set of isomorphism classes of 
simple objects in $\th\Rb$-$\fModb_0$.
For each $b$ in $\th\!B(\l)$ there is a unique $\flat$-selfdual 
irreducible graded $\th\Rb$-module
$\th\!G^\up(b)$ 
which is isomorphic to $b$ as a (non graded) $\th\Rb$-module.
We define a $\Ac$-basis 
$\th\!\Gb^\up(\l)$ of $\th\!\Gb_I$ by setting
$$\th\!\Gb^\up(\l)=\{\th\!G^\up(b);\,b\in \th\!B(\l)\},\quad
\th\!G^\up(0)=0.$$
\endproclaim

Now, we concentrate on the $\Ac$-module $\th\Kb_I$.
We equip $\th\Kb_I$ with the symmetric $\Ac$-bilinear form  
$$\gathered
\th\Kb_I\times \th\Kb_I\to \Bc,\quad
(P,Q)\mapsto (P:Q),
\endgathered\leqno(8.2)$$
see (8.1).
Consider the duality 
$$\th\Rb\text{-}\proj\to\th\Rb\text{-}\proj,\quad
P\mapsto P^\sharp=\gHom_{\th\Rb}(P,\th\Rb),$$
with the action and the grading given by 
$$(xf)(p)=f(p)\omega(x),\quad(P^\sharp)_d=\Hom_{\th\Rb}(P[-d],\th\Rb).$$
This duality functor yields an $\Ac$-antilinear map
$$\th\Kb_I\to\th\Kb_I,\quad P\mapsto P^\sharp.$$
Let $\Kc\to\Kc$, $f\mapsto \bar f$
be the unique involution such that
$\bar v=v^{-1}$.

\proclaim{8.3.~Definition}
For each $b$ in $\th\!B(\l)$ let $\th\!G^\low(b)$ 
be the unique indecomposable graded module
in $\th\Rb$-$\proj$ whose top is isomorphic to $\th\!G^\up(b)$.
We set $\th\!G^\low(0)=0$ and
$\th\!\Gb^\low(\l)=\{\th\!G^\low(b);\,b\in \th\!B(\l)\}$,
a $\Ac$-basis of $\th\!\Kb_I$.
\endproclaim

\proclaim{8.4.~Proposition} (a) We have
$\la \th\!G^\low(b):\th\!G^\up(b')\ra=\delta_{b,b'}$ 
for each $b,b'$ in $\th\!B(\l)$.

\vskip1mm

(b) We have
$\la P^\sharp:M\ra=\overline{\la P:M^\flat\ra}$ for each $P$, $M$.

\vskip1mm

(c) The graded $\th\Rb$-module $\th\!G^\low(b)$ is $\sharp$-selfdual
for each $b$ in $\th\!B(\l)$.
\endproclaim

\noindent{\sl Proof :}
Part $(a)$ is obvious because we have
$$\la \th\!G^\low(b):\th\!G^\up(b')\ra=
\gdim\,\gHom_{\th\Rb}(\th\!G^\low(b),\top\, \th\!G^\low(b'))
=\delta_{b,b'}.$$
Part $(c)$ is a consequence of $(b)$. Finally $(b)$ is proved as follows
$$\aligned
\la P^\sharp:M\ra
&=\gdim\,\gHom_{\th\Rb}(\gHom_{\th\Rb}(P,\th\Rb),M),
\vspace{2mm}
&=\gdim\,(P^\omega\otimes_{\th\Rb}M),
\vspace{2mm}
&=\overline{\gdim\,\gHom(P^\omega\otimes_{\th\Rb}M,\kb)},
\vspace{2mm}
&=\overline{\gdim\,\gHom_{\th\Rb}(P,M^\flat)},
\vspace{2mm}
&=\overline{\la P:M^\flat\ra}.
\endaligned$$
Here, the second equality holds because $P$ is a projective graded module
and the fourth one is adjointness of $\otimes$ and $\Hom$, see e.g., 
\cite{CuR,(2.19)}.

\qed

\vskip3mm

\subhead 8.5.~Example\endsubhead Set $\nu=i+\theta(i)$ and $\ib=i\theta(i)$. 
Set $\th\Rb_\ib=\th\Rb 1_\ib$ and $\th\Lb_\ib=\top(\th\Rb_\ib)$. 
Recall the description of $\th\Rb_1$ given in Example 5.9. 
Recall also that $\th\Rb_0=\kb$.
The global bases are given by the following.
First, the weight zero parts are given by
$\th\!\Gb^\low_0(\l)=\th\!\Gb^\up_0(\l)=\{\kb\}$.
Next, let us consider the weight $\nu$ parts.
\vskip2mm
\itemitem{$\bullet$}
If $\l_i+\l_{\theta(i)}\neq 0$ then 
$\th\!\Gb^\low_\nu(\l)=\{\th\Rb_\ib,\, \th\Rb_{\theta(\ib)}\}$ and 
$\th\!\Gb^\up_\nu(\l)=\{\th\Lb_\ib,\, \th\Lb_{\theta(\ib)}\}$.

\vskip2mm
\itemitem{$\bullet$} If $\l_i+\l_{\theta(i)}=0$ then 
$\th\!\Gb^\low_\nu(\l)=\{\th\Rb_\ib\}$,
$\th\!\Gb^\up_\nu(\l)=\{\th\Lb_\ib\}$, $\th\Rb_\ib=\th\Rb_{\theta(\ib)},$
and $\th\Lb_\ib=\th\Lb_{\theta(\ib)}$.

\vskip3mm

\subhead 8.6.~Definition of the operators $e_i$ and $f_i$\endsubhead 
First, let us introduce the following notation for a future use.
Given integers $m,m',n,n'\geqslant 0$ such that 
$$m+m'=n+n'=m'',$$
let $D_{m,m'}$ be the set of minimal representative in $W_{m''}$
of the cosets in 
$$W_{m,m'}\setminus W_{m''},\quad
W_{m,m'}=W_{m}\times\Sen_{m'}.$$
Recall that $W_{m''}$ is regarded as a Weyl group of type $B_{m''}$, 
see Section 5.4.
Write 
$$D_{m,m';n,n'}=D_{m,m'}\cap( D_{n,n'})^{-1}.$$
For each element $w$ of $D_{m,m';n,n'}$ we set 
$$W(w)=W_{m,m'}\cap w(W_{n,n'})w^{-1}.$$
We abbreviate
$$\th\Rb_{m,m'}=\th\Rb_{m}\otimes\Rb_{m'}.$$
For any integers $m'_1,m'_2,\dots,m'_r\geqslant 0$ 
we define the graded $\kb$-algebra
$$\th\Rb_{m,m'_1,m'_2,\dots,m'_r}$$ in the same way.
There is an unique
inclusion of graded $\k$-algebras
$$\gathered
\psi:\th\Rb_{m,m'}\to\th\Rb_{m''},\vspace{2mm}
1_\ib\otimes 1_{\ib'}\mapsto1_{\ib''}, \vspace{2mm}
1_{\ib}\otimes\varkappa_{\ib',l}\mapsto \varkappa_{\ib'',m+l},
\vspace{2mm}
1_\ib\otimes\sigma_{\ib',k}\mapsto \sigma_{\ib'',m+k},
\vspace{2mm}
\varkappa_{\ib,l}\otimes 1_{\ib'}\mapsto \varkappa_{\ib'',l},
\vspace{2mm}
\pi_{\ib,1}\otimes 1_{\ib'}\mapsto \pi_{\ib'',1},
\vspace{2mm}
\sigma_{\ib,k}\otimes 1_{\ib'}\mapsto \sigma_{\ib'',k},
\endgathered\leqno(8.3)$$
where $\ib\in\th\! I^m$, 
$\ib'\in I^{m'}$, and $\ib''=\theta(\ib')\ib\ib'$ is a sequence in
$\th\!I^{m''}$.

\proclaim{8.7.~Lemma}
The  
graded $\th\Rb_{m,m'}$-module $\th\Rb_{m''}$
is free of rank $2^{m'}\bigl(\smallmatrix m''\cr m\endsmallmatrix\bigr)$. 
\endproclaim

\noindent{\sl Proof :}
Set $\nu''=\theta(\nu')+\nu+\nu'$,
where $\nu$, $\nu'$ are vector dimensions in $\th\NN I$, $\NN I$ respectively,
such that $|\nu|=2m$ and $|\nu'|=m'$.  For each $w$ in $D_{m,m'}$
we have the element $\sigma_{\dot w}$ in 
$\th\Rb_{m''}$ defined in (5.2). 
Using filtered/graded arguments it is easy to see that
$$\th\Rb_{m''}=\bigoplus_{w\in D_{m,m'}}\th\Rb_{m,m'}\sigma_{\dot w}.$$

\qed

\vskip3mm

Now, we consider the triple of adjoint functors $(\psi_!,\psi^*,\psi_*)$ where
$$\psi^*\,:\,
\th\Rb_{m''}\text{-}\modb\to
\th\Rb_{m}\text{-}\modb\times\Rb_{m'}\text{-}\modb$$
is the restriction and $\psi_!$, $\psi_*$ are given by 
$$\aligned
&\psi_!\,:\,\cases
\th\Rb_{m}\text{-}\modb\times\Rb_{m'}\text{-}\modb\to
\th\Rb_{m''}\text{-}\modb,\vspace{2mm}
(M,M')\mapsto
\th\Rb_{m''}\otimes_{\th\Rb_{m,m'}}
(M\otimes M'),\endcases
\vspace{2mm}
&\psi_*\,:\,\cases
\th\Rb_{m}\text{-}\modb\times\Rb_{m'}\text{-}\modb\to
\th\Rb_{m''}\text{-}\modb,\vspace{2mm}
(M,M')\mapsto
\gHom_{\th\Rb_{m,m'}}(
\th\Rb_{m''},M\otimes M').\endcases
\endgathered
$$
The same discussion as for the triple
$(\phi_!,\phi^*,\phi_*)$ implies that 
$\psi_!$, $\psi^*$, $\psi_*$ are exact,
they commute with the shift of the grading, 
and they take finite dimensional modules to finite dimensional ones,
while the functors $\psi_!$, $\psi^*$
take projective modules to projective ones.
Thus the functor $\psi_!$ yields $\Ac$-bilinear maps 
$$\gathered
\Kb_I^\theta\times\Kb_I\to \th\Kb_I,
\quad
\th\!\Gb_I\times\Gb_I\to \th\!\Gb_I,
\endgathered$$ 
while $\psi^*$ yields maps in the inverse way.
For a graded $\th\Rb_m$-module $M$  we write 
$$\gathered
f_i(M)
=\th\Rb_{m+1}1_{m,i}
\otimes_{\th\Rb_{m}}M,\vspace{2mm}
e_i(M)
=\th\Rb_{m-1}
\otimes_{\th\Rb_{m-1,1}}1_{m-1,i}M.\endgathered\leqno(8.4)$$
Let us explain these formulas.
The symbols $1_{m,i}$ and $1_{m-1,i}$ are as in (6.3).
Note that $f_i(M)$ is a graded $\th\Rb_{m+1}$-module, while
$e_i(M)$ is a graded $\th\Rb_{m-1}$-module.
The tensor product in the definition of $e_i(M)$ is relative to the graded $\kb$-algebra homomorphism
$$\th\Rb_{m-1,1}=\th\Rb_{m-1}\otimes\Rb_1\to\th\Rb_{m-1}\otimes\Rb_i\to
\th\Rb_{m-1}\otimes(\Rb_i/(\varkappa_i))=\th\Rb_{m-1}.$$
In other words, let $e'_i(M)$ is the graded $\th\Rb_{m-1}$-module 
obtained by taking the direct summand $1_{m-1,i}M$
and restricting it to $\th\Rb_{m-1}$.
Observe that if $M$ is finitely generated then $e'_i(M)$ may not lie in $\th\Rb_{m-1}$-$\modb$.
To remedy this, since $e'_i(M)$ affords a $\th\Rb_{m-1}\otimes\Rb_i$-action
 we consider the graded $\th\Rb_{m-1}$-module 
$$e_i(M)=e'_i(M)/\varkappa_i e'_i(M).$$

\proclaim{8.8.~Definition}
The functors $e_i$, $f_i$ preserve the category $\th\Rb$-$\projb$, yielding
$\Ac$-linear operators on $\th\Kb_I$.
Let $e_i$, $f_i$ denote also the $\Ac$-linear operators on $\th\Gb_I$
which are the transpose of $f_i$, $e_i$ with respect to the Cartan pairing.
\endproclaim

Note that the symbols $e_i(M)$, $f_i(M)$ have different meaning if $M$ is viewed as an element of
$\th\Kb_I$ or if $M$ is viewed as an element of $\th\Gb_I$. In the first case they are given by (8.4), in the second one by the formulas in Lemma 8.9$(a)$ below. We hope this will not create any confusion.

\vskip3mm

\proclaim{8.9.~Lemma} 
(a) 
The operators $e_i$, $f_i$ on $\th\Gb_I$ are given by
$$\gathered
e_i(M)=1_{m-1,i}M,\quad
f_i(M)=
\hom_{\th\Rb_{m,1}}
(\th\Rb_{m+1},M\otimes\Lb_i),\quad
M\in\th\Rb_{m}\text{-}\fmodb.
\endgathered$$

(b) For $M,M''\in\th\Rb\text{-}\modb$ and $M'\in\Rb$-$\modb$ we have
$$(\psi_!(M,M'):M'')=(M\otimes M':\psi^*(M'')).$$
The bilinear form $(\bullet:\bullet)$ on $\th\Kb_I$ is such that
$$(e_i(P):P')=(1-v^2)(P:f_i(P')),\quad P,P'\in\th\Rb\text{-}\proj.$$

\vskip1mm

(c) We have $f_i(P)^\sharp=f_i(P^\sharp)$ for each
$P\in\th\Rb\text{-}\proj$.

\vskip2mm

(d) We have $e_i(M)^\flat=e_i(M^\flat)$ for each
 $M\in\th\Rb\text{-}\fmodb$.

\vskip2mm

(e) The operators $e_i$, $f_i$ on $\th\Kb_I$, $\th\Gb_I$ 
satisfy the relation $(7.3)$.
\endproclaim

\noindent{\sl Proof :} 
Let $M\in\th\Rb\text{-}\fmodb$ and 
$P\in\th\Rb\text{-}\proj$. We have
$f_i(P)=\psi_!(P,\Rb_i)$. Thus
$$
\aligned
\hom_{\th\Rb}(f_i(P),M)&=
\hom_{\th\Rb}(\psi_!(P,\Rb_i),M)
\vspace{2mm}
&=
\hom_{\th\Rb\otimes\Rb_i}(P\otimes\Rb_i,\psi^*(M))
\vspace{2mm}
&=\hom_{\th\Rb}(P,1_{m-1,i}M).
\endaligned$$
Next, we must prove that $f_i(M)=\psi_*(M,\Lb_i)$. We have
$$\aligned
\hom_{\th\Rb}(e_i(P),M)&=\hom_{\th\Rb\otimes\Rb_i}(1_{m-1,i}P,M\otimes(\Rb_i/\varkappa_i\Rb_i))
\vspace{2mm}
&=\hom_{\th\Rb\otimes\Rb_i}(\psi^*(P),M\otimes(\Rb_i/\varkappa_i\Rb_i))
\vspace{2mm}
&=\hom_{\th\Rb}(P,\psi_*(M,\Lb_i)).
\endaligned
$$
This proves part $(a)$.
The first claim of $(b)$ follows from the following identity
$$\aligned
(\psi_!(M,M'):M'')
&=\gdim\Bigl((M^\omega\otimes{M'}^\omega)\otimes_{\th\Rb_{m,m'}}\psi^*(M'')\Bigr),\vspace{2mm}
&=(M\otimes M':\psi^*(M'')).
\endaligned$$
The second one is proved as follows
$$\aligned
(1-v^2)(f_i(P):P')
&=(1-v^2)(\psi_!(P,\Rb_i):P')\vspace{2mm}
&=(P\otimes\Lb_i:1_{m-1,i}P')\vspace{2mm}
&=(P:e_i(P')).
\endaligned$$
Part $(c)$ follows from the following identities
$$\aligned
f_i(P)^\sharp
&=\hom_{\th\Rb_m}\bigl(
\th\Rb_m\otimes_{\th\Rb_{m-1,1}}(P\otimes\Rb_i),
\,\th\Rb_m\bigr),
\vspace{2mm}
&=\hom_{\th\Rb_{m-1,1}}\bigl(P\otimes\Rb_i,
\,\th\Rb_m\bigr),
\vspace{2mm}
&=\th\Rb_m\otimes_{\th\Rb_{m-1,1}}
\hom_{\th\Rb_{m-1,1}}\bigl(P\otimes\Rb_i,
\,\th\Rb_{m-1,1}\bigr),
\vspace{2mm}
&=\th\Rb_m\otimes_{\th\Rb_{m-1,1}}
\bigl(\hom_{\th\Rb_{m-1}}(P,
\,\th\Rb_{m-1})\otimes\Rb_i\bigr),
\vspace{2mm}
&=f_i(P^\sharp).
\endaligned$$
Here the second equality is Frobenius reciprocity and the third one follows from
Lemma 8.7, see e.g., \cite{CuR, (2.29)}.
Part $(d)$ follows from $(c)$ and Proposition 8.4$(b)$.
To prove $(e)$ it is enough to check that the operators $e_i$, $f_i$ 
on $\th\Kb_I$ satisfy the relation $(7.3)$.
For $f_i$ it is enough to observe that
$$f_i(P)=\psi_!(P,\Rb_i),\quad\forall P\in\th\Rb\text{-}\proj.$$ 
Then the claim follows from Theorem 7.7 and the associativity of induction.
For $e_i$ it is enough to observe that the transposed operator is given by
$$f_i(M)=\psi_*(M,\Lb_i),\quad M\in\th\Rb\text{-}\fmodb$$ and to use the associativity of coinduction.

\qed

\vskip3mm

\subhead 8.10.~Shuffles, projectives, and characters\endsubhead
For each sequence $\ib$ in $\th\! I^m$ we define a projective graded module 
in $\th\Rb_{m}$-$\proj$ by setting $\th\Rb_\ib=\th\Rb_m1_\ib$. 
More generally, for $\yb\in\th Y^m$ we define an object 
$\th\Rb_\yb$ of $\th\Rb_{m}$-$\proj$  as follows.
Write $$\yb=(\theta(\jb)\jb,\theta(\bb)\bb),\quad\jb\in I^m,\quad\bb\in\ZZ^m.$$ 
We may abbreviate $\yb=\theta(\zb)\zb$ where $\zb=(\jb,\bb)$. 
Define the idempotent 
$1_\yb$ as the image of the idempotent $1_{\zb}$ by the
inclusion $\psi:\Rb_m\to\th\Rb_m$ given by setting $m,m',m''$ equal to
$0,m,m$ in (8.3).  Then set
$$\th\Rb_\yb=(\th\Rb_m1_\yb)[\ell_{\bb}].$$
The graded module $\th\Rb_\yb$ 
satisfies the same properties as the 
projective graded $\Rb_m$-module $\Rb_\zb$ introduced 
in Section 7.4. In particular $\th\Rb_\yb$ is $\sharp$-selfdual, 
and if the sequence
$\ib$ in $\th\!I^m$ is the expansion of the pair $\yb$ then we have
$$\th\Rb_{\ib}=\la\bb\ra !\,\th\Rb_\yb.\leqno(8.5)$$ 
For $\yb=(\ib,\ab)$ in $\th Y^m$ and $\yb'=(\ib',\ab')$ in $Y^{m'}$ 
we have $$\psi_!(\th\Rb_\yb,\Rb_{\yb'})=\th\Rb_{\yb''},\quad
\yb''=\theta(\yb')\yb\yb'=(\theta(\ib')\ib\ib',\theta(\ab')\ab\ab').$$
Write $\ib'=(i'_1,\dots, i'_k)$,
$\ab'=(a'_1,\dots a'_k)$,  and
$$f_{\yb'}=
f_{i'_1}^{(a'_1)}f_{i'_2}^{(a'_2)}\cdots f_{i'_k}^{(a'_k)}.$$
Lemma 8.9$(a)$ yields
$$f_{\yb'}(P)=\psi_!(P,\Rb_{\yb'}),\quad P\in\th\Rb_m\text{-}\proj.$$
In particular,  we have
$$f_{\yb'}(\kb)=\th\Rb_{\theta(\yb')\yb'}.$$

\proclaim{8.11.~Definition}
A {\it shuffle} of a pair of sequences $(\ib,\ib')$ in $\th\!I^m\times I^{m'}$
is a sequence $\ib''$ in $\th\!I^{m''}$ 
together with a subsequence of $\ib''$ isomorphic to $\ib$
and such that the complementary subsequence 
is equal to $\theta(\ib')\ib'$ modulo $\theta$. 
\endproclaim

\noindent
Let $Sh(\ib,\ib')$ be the set of shuffles of $\ib$, $\ib'$.
The assignment $w\mapsto w^{-1}(\theta(\ib')\ib\ib')$ 
gives a bijection from $D_{m,m'}$ to $Sh(\ib,\ib')$.
To a shuffle $\ib''$ in $Sh(\ib,\ib')$ associated with an element
$w$ of $D_{m,m'}$ we assign the following degree 
$$\deg(\ib,\ib';\ib'')=\deg(\sigma_{\dot w}1_{\ib''}).$$
This degree does not depend of the choice of the reduced decomposition 
$\dot w$ of $w$. 
Let $\th\!\Bc I^m$ be the free $\Bc$-module with basis $\th\!I^m$.
For any $f$ in $\th\!\Bc I^m$ we write
$$f=\sum_\ib f(\ib)\,\ib.$$

\proclaim{8.12.~Definitions}
(a) For any finitely generated graded
$\th\Rb_m$-module $M$ we define the {\it character of $M$} 
as the element of $\th\!\Bc I^m$ given by
$$\ch(M)=\sum_\ib\gdim(1_\ib M)\,\ib.$$
(b) For any elements
$f\in\th\!\Bc I^m$, $g\in\Bc I^{m'}$ we define their product 
$f\circledast g\in\th\!\Bc I^{m''}$ by
$$(f\circledast g)(\ib'')=
\sum_{\ib,\ib'}v^{\deg(\ib,\ib';\ib'')}f(\ib)g(\ib').$$
Here the sum is over all ways to represent $\ib''$ as a shuffle of $\ib$ and $\ib'$. 
\endproclaim

\noindent

\proclaim{8.13.~Proposition} For any $M\in\th\Rb_m$-$\modb$
and any $M'\in\Rb_{m'}$-$\modb$ we have
$$\ch(\psi_!(M,M'))=\ch(M)\circledast\ch(M').$$
\endproclaim

\noindent{\sl Proof :}
We have $\ch(M)=\sum_\ib(\th\Rb_\ib: M)\,\ib$. 
Thus Lemma 8.9$(b)$  yields
$$\ch(\psi_!(M,M'))=\sum_{\ib''}(\psi^*(\th\Rb_{\ib''}):M\otimes M')\,\ib''.$$
Next, we have the following formula 
$$\psi^*(\th\Rb_{\ib''})=
\bigoplus_{\ib,\ib'}\th\Rb_{\ib}\otimes\Rb_{\ib'}[\deg(\ib,\ib';\ib'')],$$
where the sum runs over all sequences $\ib$, $\ib'$ such that
$\ib''$ lies in $Sh(\ib,\ib')$. This formula is a consequence
of the Mackey's induction-restriction theorem. 
The details are left to the reader. See e.g., the proof of Theorem 8.31 below.
Therefore we get
$$\aligned
\ch(\psi_!(M,M'))&=
\sum_{\ib,\ib'}\sum_{\ib''\in Sh(\ib,\ib')}(\th\Rb_{\ib}:M)\,(\Rb_{\ib'}:M')\,v^{\deg(\ib,\ib';\ib'')}\,\ib''\vspace{2mm}
&=\ch(M)\circledast\ch(M').\endaligned$$

\qed

\vskip3mm

\proclaim{8.14.~Proposition}
We have
$$f_i(\th\Rb_\ib)=\th\Rb_{\theta(i)\ib i},\quad
e_i(\th\Rb_\ib)=
\bigoplus_{\ib'}\th\Rb_{\ib'}[\deg(\ib',i;\ib)].$$
Here the sum runs over all sequences in $\th\!I^{m-1}$ such that
$\ib$ lies in $Sh(\ib',i)$.
\endproclaim

\noindent{\sl Proof :} Left to the reader.

\qed

\vskip3mm

\subhead 8.15.~Example\endsubhead Set $\nu=i+\theta(i)$ and
$\ib=i\theta(i)$. Using the description of
$\th\Rb_0$, $\th\Rb_1$ given in Example 5.9
we can compute $e_j$, $f_j$ and ch.
Let us regard $\kb$ as an object of
$\th\Rb_0$-$\projb$. We have
$$
f_i(\kb)=\th\Rb_{\theta(\ib)},\quad
e_j(\th\Rb_\ib)=\cases
v^{\lambda_i+\lambda_{\theta(i)}}\,\kb&\roman{if}\ j=i,
\vspace{2mm}
\kb&\roman{if}\ j=\theta(i),
\vspace{2mm}
0&\roman{else}.\endcases
$$
Next, observe that $\th\Rb_{0,1}=\Rb_1$ and that the inclusion 
in (8.3) yields following formula
$\th\Rb_1=\Rb_1\oplus\Rb_1\pi_1$. Thus we get
$$\gathered
\psi_*(\kb,\Rb_i)=\hom_{\Rb_1}(\th\Rb_1,\Rb_i)
=\th\Rb_\ib[\l_i+\l_{\theta(i)}],
\vspace{2mm}
\psi_!(\kb,\Rb_i)=\th\Rb_1\otimes_{\Rb_1}\Rb_i
=\th\Rb_{\theta(\ib)}.
\endgathered
$$
In particular, we have the following.
\vskip2mm
\itemitem{$\bullet$}
If $\l_i+\l_{\theta(i)}\neq 0$ then $e_j(\th\Lb_{\theta(\ib)})=\kb$ if $j=i$ and 0 else, and
$f_i(\kb)=v^{\lambda_i+\lambda_{\theta(i)}}\th\Lb_\ib+\th\Lb_{\theta(\ib)}$.
Further $\ch(\th\Lb_{\theta(\ib)})=\theta(\ib)$ and $\ch(\th\Lb_\ib)=\ib.$
\vskip2mm
\itemitem{$\bullet$} If $\l_i+\l_{\theta(i)}=0$ then 
$e_j(\th\Lb_{\theta(\ib)})=\kb$ if $j=i,\theta(i)$ and 0 else, and
$f_i(\kb)=\th\Lb_\ib.$
Further
$\ch(\th\Lb_\ib)=\ib+\theta(\ib).$

\vskip3mm

\subhead 8.16.~Induction of $\Hb_m$-modules versus induction of 
$\th\Rb_{m}$-modules\endsubhead
Recall the functors $E_i$, $F_i$ on
$\Hb\text{-}\fModb_I$ defined in (6.4).
We have also the functors
$$\Psi:\th\Rb_{m}\text{-}\fModb_0\to\Hb_m\text{-}\fModb_I,\quad
\forb: 
\th\Rb_m\text{-}\fmodb\to
\th\Rb_m\text{-}\fModb_0,$$
where $\forb$ is the forgetting of the grading.
Finally we define functors
$$\gathered
E_i:\th\Rb_{m}\text{-}\fModb_0\to\th\Rb_{m-1}\text{-}\fModb_0,\quad
E_iM=1_{m-1,i}M,
\vspace{2mm}
F_i:\th\Rb_{m}\text{-}\fModb_0\to\th\Rb_{m+1}\text{-}\fModb_0,\quad
F_iM=\psi_!(M,\Lb_i).
\endgathered\leqno(8.6)$$

\proclaim{8.17.~Proposition} There are canonical isomorphisms of functors
$$E_i\circ\Psi=\Psi\circ E_{i},\quad
F_i\circ\Psi=\Psi\circ F_{i},\quad
E_i\circ\forb=\forb\circ e_i,\quad
F_i\circ\forb=\forb\circ f_{\theta(i)}.$$
\endproclaim

\noindent{\sl Proof :}
Recall that $\kb_i$ is the 1-dimensional $\kb[X_{m+1}^{\pm 1}]$-module
such that $X_{m+1}\mapsto i^{-1}$, that
$\Lb_i$ is the 1-dimensional $\Rb_1$-module such that
$1_i\mapsto 1$ and $\varkappa_1\mapsto 0$, and that
$\Psi$ identifies $X_{m+1}$ and the element
$1\otimes i^{-1}f(\varkappa_{m+1})$ in $\th\Rb_{m,1}$, where the
function $f$ is as in Theorem 6.5.
The first two isomorphisms are obvious consequences of (6.4), (8.6), 
because $E_iM$ is the generalized $i^{-1}$-eigenspace of $M$ 
with respect to the action of $X_{m+1}$, and 
$F_iM$ is induced from the $\Hb_m\otimes\kb[X_{m+1}^{\pm 1}]$-module
$M\otimes\kb_{i}$.
The third isomorphism follows from (8.6) and  Lemma 8.9$(a)$.
Now, we concentrate on the last isomorphism.
Lemma 8.9$(a)$ yields
$$f_{\theta(i)}(M)=\psi_*(M,\Lb_{\theta(i)}),\quad
M\in\th\Rb\text{-}\fmodb.$$
For any graded
$\Rb$-module $N$ let $N^\kappa$ be equal to $N$, 
with the $\Rb$-action twisted by the
involution $\kappa$ in (7.1). Note that $\Lb_i^\kappa=\Lb_{\theta(i)}$.
Therefore the proposition follows from the following lemma. 

\proclaim{8.18.~Lemma} For each $M\in\th\Rb_m\text{-}\fmodb$ and 
$N\in\Rb_{m'}\text{-}\fmodb$
there is an isomorphism of (non-graded)
$\th\Rb_{m''}$-modules $\psi_!(M,N)=\psi_*(M,N^\kappa).$
\endproclaim

\noindent{\sl Proof :} 
Recall that $\th\Rb_{m,m'}=\th\Rb_m\otimes\Rb_{m'}$.
The involution $\kappa:\Rb_{m'}\to\Rb_{m'}$ in (7.1) yields an involution
of $\th\Rb_{m,m'}$. Let us denote it by $\kappa$ again.
Let $\th\Rb_{m,m'}^\kappa$ be the
$(\th\Rb_{m,m'},\th\Rb_{m,m'})$-bimodule which is equal to
$\th\Rb_{m,m'}$ as a  right $\th\Rb_{m,m'}$-module, and such that the
left $\th\Rb_{m,m'}$-action is twisted by $\kappa$.
It is enough to prove that there is an isomorphism of (non-graded)
$(\th\Rb_{m''},\th\Rb_{m,m'})$-bimodules
$$\th\Rb_{m''}\to\hom_{\th\Rb_{m,m'}}(\th\Rb_{m''},\th\Rb_{m,m'}^\kappa).$$
The bimodule structure on the right hand side is given by
$$(xfy)(z)=f(zx)y,\quad
x,z\in\th\Rb_{m''},\ y\in\th\Rb_{m,m'}.$$
Lemma 8.7 yields an isomorphism
$$\th\Rb_{m''}=\bigoplus_{w\in D_{m,m'}}\th\Rb_{m,m'}\sigma_{\dot w}$$ 
of graded $\th\Rb_{m,m'}$-modules.
The longest double coset representative in $D_{m,m';m,m'}$ 
is the coset of the involution $u\in W_{m''}$ given by
$$u=w_{m'}\,\eps_{m+1}\dots\eps_{m''},$$
with $w_{m'}$ the longest element of $\Sen_{m'}$.
There is an unique morphism of $(\th\Rb_{m,m'},\th\Rb_{m,m'})$-bimodules
$$h:\th\Rb_{m,m'}\to
\hom_{\th\Rb_{m,m'}}(\th\Rb_{m''},\th\Rb_{m,m'}^\kappa)
,$$
taking 1 to the map
$$y\,\sigma_{\dot w}\mapsto
\kappa(y)\,\delta_{w,u},\quad y\in\th\Rb_{m,m'},\quad
w\in D_{m,m';m,m'}.$$
Since the right hand side is a left $\th\Rb_{m''}$-module,
by Frobenius reciprocity
$h$ yields a morphism of  $(\th\Rb_{m''},\th\Rb_{m,m'})$-bimodules
$$\th\Rb_{m''}\to\hom_{\th\Rb_{m,m'}}(\th\Rb_{m''},\th\Rb_{m,m'}^\kappa).$$
This map is invertible.
The proof is the same as in
\cite{M, sec.~3}, \cite{LV, thm.~2.2}.

\qed

\vskip3mm

\proclaim{8.19.~Proposition} 
(a) The functor $\Psi$ yields an isomorphism of Abelian groups
$$\bigoplus_{m\geqslant 0}[\th\Rb_{m}\text{-}\fModb_0]=
\bigoplus_{m\geqslant 0}[\Hb_m\text{-}\fModb_I].$$
The functors $E_i$, $F_i$ yield endomorphisms of 
both sides which are intertwined by $\Psi$. 

\vskip1mm

(b) The forgetful functor $\forb$ factors to a group isomorphism 
$$\th\!\Gb_I/(v-1)=
\bigoplus_{m\geqslant 0}[\th\Rb_{m}\text{-}\fModb_0].$$
\endproclaim

\noindent{\sl Proof :}
Claim $(a)$ follows from Corollary 6.7 and Proposition 8.17.
Claim $(b)$ follows from Proposition 8.2.

\qed

\vskip3mm

\subhead 8.20.~The crystal operators on $\th\Gb_I$ and $\th\!B(\l)$\endsubhead 
Fix a vertex $i$ in $I$. 
For each irreducible graded module $M\in\th\Rb$-$\fmodb$ we define
$$\gathered
\tilde e_i(M)=\soc\, (e_i(M)),\quad\tilde f_i(M)=\top\, \psi_!(M,\Lb_i),\vspace{2mm}
\eps_i(M)=\max\{n\geqslant 0;e_i^n(M)\neq 0\},\quad
\eps_i(M)\in\NN\cup\{\infty\}.
\endgathered$$
For each positive integers $m\geqslant n$  we consider the functor
$$\Delta_{ni}:\th\Rb_m\text{-}\fmodb\to\th\Rb_{m-n}\text{-}\fmodb
\times\Rb_{ni}\text{-}\fmodb,\quad M\mapsto 1_{m-n,ni}M.$$
Given an irreducible graded module $M\in\th\Rb_m$-$\fmodb$ we have, see Lemma 8.9$(a)$,
$$\eps_i(M)=\max\{n\geqslant 0;\Delta_{ni}(M)\neq 0\},\quad e_i(M)=\Delta_i(M).$$

\proclaim{8.21.~Proposition}
Let $M$ be an irreducible graded $\th\Rb_m$-module and $n$ be an integer $\geqslant 0$.
Set $\eps=\eps_i(M)$,
$M^+=\psi_!(M,\Lb_{ni})$ and $M^-=\Delta_{ni}(M)$.

\vskip1mm

(a)  If $\eps=0$ then $\Delta_{ni}(M^+)=M\otimes\Lb_{ni}$, $\top(M^+)$ is irreducible,
 $\eps_i(\top(M^+))=n$, 
all other composition factors $L$ of $M^+$ have
$\eps_i(L)<n$. 
\vskip1mm

(b)  If $\eps\geqslant n$ then any irreducible submodule  of $M^-$ is of the form
$N\otimes\Lb_{ni}$ with $\eps_i(N)=\eps-n$.
 If $\eps=n$ then $M^-$ is irreducible. 
If $\eps\geqslant n$ then
$\soc(M^-)$ is  irreducible.
In particular $\tilde e_i(M)$ is irreducible if $\eps\neq 0$ and 0 else. Finally  we have 
$\soc(M^-)=\tilde e_i^n(M)\otimes\Lb_{ni}$. 
\vskip1mm

(c) $\top(M^+)$ is irreducible,
 $\eps_i(\top(M^+))=\eps+n$, 
and all other composition factors $L$ of $M^+$ have
$\eps_i(L)<\eps+n$. In particular $\tilde f_i(M)$ is irreducible. 
Finally we have $\top(M^+)=\tilde f_i^n(M)$. 
\endproclaim

\noindent

\noindent{\sl Proof :} 
Part $(a)$ is the analogue of \cite{K, lem.~5.1.3},  \cite{KL, lem.~3.7}.
More precisely, note first that we have
$$\ch(\Delta_{ni}(M^+))=\sum_\ib\gdim(1_{\theta(i^n)\ib i^n}M^+)\,\theta(i^n)\ib i^n.
\leqno(8.7)$$
Hence, since $\eps=0$ Proposition 8.13 implies that $$\dim(\Delta_{ni}(M^+))=\dim(M\otimes\Lb_{ni}).$$
Since $\Delta_{ni}(M^+)$ contains a copy of $M\otimes\Lb_{ni},$ we get the first claim of $(a)$.
By Frobenius reciprocity, a copy of
$M\otimes\Lb_{ni}$, possibly with a grading shift, appears as a submodule of $\Delta_{ni}(M')$ 
for any nonzero quotient
$M^+\to M'$. Since $$\Delta_{ni}(M^+)=M\otimes\Lb_{ni},$$ this implies that 
$\top(M^+)$ is irreducible with $\eps_i(\top(M^+))\geqslant n$, that
$$\Delta_{ni}(M^+)=\Delta_{ni}(\top(M^+)),$$ and that $\Delta_{ni}(L)=0$ for all other composition factors
$L$ of $M^+$.
Finally we have
$\eps_i(\top(M^+))=n$, because $\eps=0$. 

Now we prove  $(b)$. The first claim is the analogue of \cite{K, lem.~5.1.2}. Indeed, any irreducible submodule
of $M^-$ is of the form $N\otimes \Lb_{ni}$ with $N$ irreducible. We have $\eps_i(N)\leqslant \eps-n$
by definition of $\eps_i$. For the reverse inequality, Frobenius reciprocity and the irreducibility of $M$ imply that $M$ is a quotient of $\psi_!(N,\Lb_{ni})$. So applying the exact functor $\Delta_{\eps i}$  we see that
$\Delta_{\eps i}(M)$ is a quotient of $\Delta_{\eps i}\psi_!(N,\Lb_{ni})$. In particular
$$\Delta_{\eps i}\psi_!(N,\Lb_{ni})\neq 0.$$
By Proposition 8.13 and (8.7) we have also
$\Delta_{(\eps-n) i}(N)\neq 0$. Thus $\eps_i(N)=\eps-n$.
The second claim of $(b)$ is the analogue of \cite{K, lem.~5.1.4}. Indeed, if $\eps=n$ then
any irreducible submodule of $M^-$ is
of the form $N\otimes\Lb_{ni}$ with $\eps_i(N)=0$. Once again Frobenius reciprocity and the irreducibility of $M$ imply that $M$ is a quotient of $\psi_!(N,\Lb_{ni})$. Hence $M^-$ is a quotient of 
$\Delta_{n i}\psi_!(N,\Lb_{ni})$. But the later is isomorphic to $N\otimes\Lb_{ni}$ by $(a)$. 
Next, the third claim of $(b)$
is the analogue of \cite{K, lem.~5.1.6}, \cite{KL, prop.~3.10}. 
Indeed, suppose that $N\otimes\Lb_{ni}\subset\soc(M^-)$.
Then $\eps_i(N)=\eps-n$ by the first part of $(b)$. Thus $N$ contributes a non-trivial submodule to
$\Delta_{\eps i}(M)$. But $\Delta_{\eps i}(M)$ is an irreducible graded
$\th\Rb_{m-\eps,\eps}$-module
by the second part of $(b)$.
Thus the socle of $\Delta_{\eps i}(M)$ as a graded
$\th\Rb_{m-\eps,\eps-n,n}$-module is
$N\otimes\Lb_{(\eps-n)i}\otimes\Lb_{ni}$ by Example 7.5.
Hence $\soc(M^-)$ must equal $N\otimes\Lb_{ni}$.
Finally, the last claim of $(b)$
is the analogue of \cite{K, lem.~5.2.1$(i)$},  \cite{KL, lem.~3.13}.
Indeed, note first that if $n>\eps$ then $$\soc(M^-)=\tilde e_i^n(M)=0.$$
Assume now that $\eps\geqslant n$.
Observe that $$\tilde e_i(M)=\soc(\Delta_i(M))$$ is irreducible or zero by the third part of $(b)$. 
Hence $\tilde e_i(M)\otimes\Lb_i$ is a submodule of $\Delta_i(M)$.
Applying this $n$ times we deduce that
$\tilde e_i^n(M)\otimes(\Lb_i)^{\otimes n}$ is a submodule of $\Delta_{ni}(M)$
as a graded $\th\Rb_{m-n,1^n}$-module.
Hence $\tilde e_i^n(M)\otimes\Lb_{ni}$ is a submodule of $\Delta_{ni}(M)$ by Frobenius reciprocity.

Finally,  we prove $(c)$. It
is the analogue of \cite{K, lem.~5.2.1$(ii)$},  \cite{KL, lem.~3.13}.
Indeed, by $(b)$ the graded module $\Delta_{\eps i}(M)$ is
of the form $N\otimes\Lb_{ni}$, with $N$ irreducible such that $\eps_i(N)=0$.
Thus, by Frobenius reciprocity $M$ is a quotient of $\psi_!(N,\Lb_{\eps i})$.
So the transitivity of induction implies that $M^+$ is a quotient of $\psi_!(N,\Lb_{(\eps+n) i})$.
Hence all claims except the last one follow from $(a)$. Finally, by exactness of the induction
$\tilde f_i^n(M)$ is a quotient of $M^+$, hence they are equal by simplicity of the top.

\qed

\vskip3mm

For each irreducible module $b$ in $\th\!B(\l)$
we define
$$\tilde E_i(b)=\soc(E_ib),\quad\tilde F_i(b)=\top(F_ib),\quad
\eps_i(b)=\max\{n\geqslant 0;E_i^nb\neq 0\}.\leqno(8.8)$$
Hence,  we have
$$\forb\circ\tilde e_i=\tilde E_i\circ\forb,\quad
\forb\circ\tilde f_i=\tilde F_i\circ\forb,\quad
\eps_i=\eps_i\circ\forb.$$

\proclaim{8.22.~Proposition} For each $b,b'$ in $\th\!B(\l)$ 
we have

\vskip1mm

(a) $\tilde F_i(b)\in \th\!B(\l)$,

\vskip1mm

(b) $\tilde E_i (b)\in \th\!B(\l)\cup\{0\}$,

\vskip1mm

(c) $\tilde F_i(b)=b'\iff\tilde E_i(b')=b$,

\vskip1mm

(d) $\eps_i(b)=\max\{n\geqslant 0;\tilde E_i^n(b)\neq 0\}$,

\vskip1mm

(e) $\eps_i(\tilde F_i(b))=\eps_i(b)+1$,

\vskip1mm

(f) if $\tilde E_i(b)=0$ for all $i$ then $b=\kb$.
\endproclaim

\noindent{\sl Proof :}
Parts $(a)$, $(b)$, $(d)$, $(e)$ and $(f)$ are immediate consequences of Proposition 8.21.
Part $(c)$ is proved as in \cite{K, lem.~5.2.3}. 
More precisely, let $M,N$ be irreducible graded modules.
By Proposition 8.21$(c)$ we have
$$\aligned
\tilde f_i(M)=N&\iff\Hom_{\th\Rb}(\psi_!(M,\Lb_i),N)\neq 0
\vspace{2mm}
&\iff\Hom_{\th\Rb}(M\otimes\Lb_i,e_i(N))\neq 0
\vspace{2mm}
&\iff\Hom_{\th\Rb}(M\otimes\Lb_i,\soc (e_iN))\neq 0
\vspace{2mm}
&\iff M=\tilde e_i(N).
\endaligned$$
Note that the proposition can also be deduced from \cite{M, Section 4} 
and Proposition 8.17. 
\qed

\vskip3mm

\proclaim{8.23.~Proposition}
The following identity holds in $\th\Kb_I$
$$f_i \th\!G^\low(b)=\la\eps_i(b)+1\ra\th\!G^\low(\tilde F_ib)+
\sum_{b'}f_{b,b'}\th\!G^\low(b'),\quad\forall b\in \th\!B(\l),$$ 
where $b'$ runs over the elements
of $\th\!B(\l)$ such that $\eps_i(b')>\eps_i(b)+1$, and
$f_{b,b'}\in\Ac$.
\endproclaim

\noindent{\sl Proof :}
We claim that there are elements 
$f_{b',b}$ in $\Ac$ such that
$$e_i \th\!G^\up(b)=\la\eps_i(b)\ra\th\!G^\up(\tilde E_ib)+
\sum_{b'}f_{b',b}\th\!G^\up(b'),\leqno(8.9)$$ 
where $b'$ runs over the elements of $\th\!B(\l)$ with $\eps_i(b')<\eps_i(b)-1$.
Taking the transpose with respect to the Cartan pairing, Proposition 8.4$(a)$,
Definition 8.8, and Proposition 8.22 yield 
$$\aligned
f_i \th\!G^\low(b)
&=\la\eps_i(\tilde F_ib)\ra\th\!G^\low(\tilde F_ib)+
\sum_{b'}f_{b,b'}\th\!G^\low(b'),
\vspace{1mm}
&=\la\eps_i(b)+1\ra\th\!G^\low(\tilde F_ib)+
\sum_{b'}f_{b,b'}\th\!G^\low(b'),
\endaligned$$ 
where $b'\in \th\!B(\l)$ with $\eps_i(b)+1<\eps_i(b')$.
Now, let us prove (8.9).
This is the analogue of \cite{K, lem.~5.5.1$(i)$}.
Fix an irreducible $\th\Rb_m$-module $b$. Set 
$$\eps=\eps_i(b),\quad
M=\th\!G^\up(b),\quad
N=\th\!G^\up(\tilde E_i^\eps b).$$
We can assume that $\eps>0$, because else (8.9) is obvious.
Note that $\eps_i(N)=0$ by Proposition 8.22.
By Frobenius reciprocity and Proposition 8.21$(b)$
there is a short exact sequence of graded modules
$$0\to R\to \psi_!(N,\Lb_{\eps i})\to M\to 0.$$
Applying the functor $e_i$ we obtain the following exact sequence
of graded modules
$$0\to e_iR\to e_i\psi_!(N,\Lb_{\eps i})\to e_iM\to 0.$$
Note that 
$$e_i\psi_!(N,\Lb_{\eps i})=1_{m-1,i}\th\Rb_m1_{m-\eps,\eps i}\otimes_{
\th\Rb_{m-\eps,\eps}}(N\otimes\Lb_{\eps i}).$$
Note also that
$$\gathered
D_{m-1,1;\,m-\eps,\eps}=\{e,x,y\},
\vspace{2mm}
x=s_{m-1}\dots s_{m-\eps+1}s_{m-\eps},\quad
y=s_{m-1}\dots s_{m-\eps+1}\eps_{m-\eps+1} ,\vspace{2mm}
W(e)=W_{m-\eps,\eps-1,1},\quad W(x)=W_{m-\eps-1,1,\eps},\quad
W(y)=W_{m-\eps,\eps-1,1}.
\endgathered$$
By Proposition 5.5 we can filter the graded 
$(\th\Rb_{m-1,1},\th\Rb_{m-\eps,\eps})$-bimodule
$$1_{m-1,i}\th\Rb_m1_{m-\eps,\eps i}.$$
This filtration is the same as in the Mackey induction-restriction theorem. 
Compare Lemma 8.32 below and the references there.
The associated graded is a direct sum of graded 
$(\th\Rb_{m-1,1},\th\Rb_{m-\eps,\eps})$-bimodules
labelled by elements of $\{e,x,y\}$
$$\gr(1_{m-1,i}\th\Rb_m1_{m-\eps,\eps i})=P_e\oplus P_x\oplus P_y.$$
We have
$$P_x1_{m-\eps-1,\eps i +i}=P_x.$$
Thus, since $e_i(N)=1_{m-\eps-1,i}N=0$, we have also
$$P_x\otimes_{
\th\Rb_{m-\eps,\eps}}(N\otimes\Lb_{\eps i})=0.$$
Next, we have
$$P_y1_{m-\eps,\theta(i),\eps i -i}=P_y.$$
Since $1_\nu\Lb_{\eps i}=0$ if $\nu\neq\eps i$ 
we have also
$$P_y\otimes_{
\th\Rb_{m-\eps,\eps}}(N\otimes\Lb_{\eps i})=0.$$
Finally, we have
$$P_e=1_{m-1,i}\th\Rb_{m-1,1}\otimes_{\Rb'}\th\Rb_{m-\eps,\eps}1_{m-\eps,\eps i},
\quad\Rb'=\th\Rb_{m-\eps,\eps-1,1}.$$
Therefore, we obtain
$$e_i\psi_!(N,\Lb_{\eps i})=
\th\Rb_{m-1,1}\otimes_{
\Rb'}(N\otimes\Lb_{\eps i}).$$
By Example 7.5 the graded $\Rb_{\eps-1,1}$-module
$\Lb_{\eps i}$ has a filtration by graded submodules whose associated
graded is isomorphic to $$\la\eps\ra\Lb_{(\eps-1)i}\otimes\Lb_i.$$
Therefore, up to some filtration, we have
$$e_i\psi_!(N,\Lb_{\eps i})=
\la\eps\ra\psi_!(N,\Lb_{(\eps-1)i})\otimes\Lb_i.$$
Now, by Proposition 8.21$(a),(c)$ we have 
$$\top\,\psi_!(N,\Lb_{(\eps-1)i})=\tilde f_i^{\eps-1}(N)=\tilde e_i(M)$$
and all other composition factors $L$ of $\psi_!(N,\Lb_{(\eps-1)i})$
have $\eps_i(L)<\eps-1$.
Moreover, by Proposition 8.21$(a)$ all composition factors $L$ of $R$ 
have $\eps_i(L)<\eps$. 
Thus, by Proposition 8.21$(b)$ all composition factors of $e_i(R)$ are
of the form $L\otimes\Lb_i$ with $\eps_i(L)<\eps-1$. 
Therefore, we obtain
$$e_i(M)=
\la\eps\ra\,\tilde e_i(M)+\sum_rf_rN_r,\quad f_r\in\Ac,$$
where $N_r$ is an irreducible graded module with $\eps_i(N_r)<\eps-1$.

\qed

\vskip3mm

\subhead 8.24.~Example\endsubhead
Set  $\nu=i+\theta(i)$ and $\ib=i\theta(i)$.
Let us compute $\tilde e_j$ and $\eps_j$.
\vskip2mm
\itemitem{$\bullet$}
If $\l_i+\l_{\theta(i)}\neq 0$ then 
$\tilde e_j(\th\Lb_\ib)=\kb$ 
if $j=\theta(i)$ and 0 else. 
We have $\eps_j(\th\Lb_\ib)=1$ if $j=\theta(i)$ and 0 else.
\vskip2mm
\itemitem{$\bullet$} If $\l_i+\l_{\theta(i)}=0$ then 
$\tilde e_j(\th\Lb_\ib)=\kb$ 
if $j=i,\theta(i)$ and 0 else. 
We have $\eps_j(\th\Lb_\ib)=1$ if $j=i,\theta(i)$ and 0 else.

\vskip3mm

\subhead 8.25.~Remark\endsubhead 
If $M$ is an irreducible graded $\th\Rb_m$-module such that
$\eps_i(M)=m$ then Proposition 8.21$(b)$ implies that
$\tilde e_i^m(M)=\kb$ and $\Delta_{mi}(M)=\kb\otimes\Lb_{mi}$.
Further, by Proposition 8.22$(c)$ there is a unique $M$ as above, up to isomorphism,
such that $\tilde e_i^m(M)=\kb$. We claim that $M\simeq\th\Lb_{mi}$, the top
of  $\th\Rb_{mi}$. 
By Proposition 8.21$(a)$ we have
$M=\top(\psi_!(\kb,\Lb_{mi}))$.
First, recall that $\th\Rb_{mi}=\psi_!(\kb,\Rb_{mi})$.
Thus, since $\psi_!$ is exact, there is a surjective map
$\th\Lb_{mi}\to M$. So it is enough to check that
$\th\Lb_{mi}$ is irreducible (left to the reader).
We'll not need this.

\vskip3mm

\subhead 8.26.~Proof of Proposition 8.2\endsubhead 
First, we prove the following.

\proclaim{8.27.~Proposition}
The character map $\ch:\th\!\Gb_{I,m}\to\th\!\Bc I^m$ is injective.
\endproclaim

\noindent{\sl Proof :} The proof is similar to that of \cite{K, thm.~5.3.1}.
We must prove that the characters of the irreducible graded modules 
in $\th\Rb_m$-$\fmodb$ are linearly independent.
We proceed by induction on $m$, the case $m=0$ being trivial.
Suppose $m>0$ and there is a non-trivial $\Ac$-linear dependence
$$\sum_Mc_M\,\ch(M)=0.\leqno(8.10)$$
We'll show by downward induction on $\eps_i(M)$ that $c_M=0$
for each graded $\th\Rb_m$-module $M$ which enter in (8.10).
Fix $M$ as above.  We have $\eps_i(M)\leqslant m$.
First, assume that $\eps_i(M)=m$. Note that $M$ is the unique irreducible graded 
$\th\Rb_m$-module such that
$\Delta_{mi}(M)=0$. Indeed we have $M=\th\Lb_{mi}$,
see Remark 8.25.
Applying $\Delta_{mi}$ to the $\th\Rb_m$-modules which enter in (8.10) 
and using the formula
$$\ch(\Delta_{mi}(M))=\sum_\ib\gdim(1_{\theta(i^m)\ib i^m}M)\,\theta(i^m)\ib i^m,
$$
we deduce that the
coefficient $c_M$ is zero.
Now, assume that $\eps_i(M)=k<m$ and that we have shown that $c_N=0$ for all
$N$ with $\eps_i(N)>k$. Applying $\Delta_{ki}$ to the $\th\Rb_m$-modules 
which enter in (8.10) we get
$$\sum_Nc_N\,\ch(\Delta_{ki}N)=0,$$
where $N$ runs over all irreducible graded 
$\th\Rb_m$-modules with $\eps_i(N)=k$. For such a graded module
$N$ we have
$\Delta_{ki}(N)=\tilde e_i^k(N)\otimes\Lb_{ki}$ by Proposition 8.21$(b)$.
Further $\tilde e_i^k(N)\neq\tilde e_i^k(N')$ if $N\neq N'$ by Proposition 8.22$(c)$.
So we conclude by the induction hypothesis.

\qed

\vskip3mm

Now, we can prove Proposition 8.2.
Forgetting the grading takes irreducible graded 
$\th\Rb$-modules to irreducible modules,
and any irreducible module in $\th\Rb$-$\fModb_0$ comes from an  
irreducible graded module in $\th\Rb$-$\fmodb$ which is unique up to isomorphism
and up to grading shift, see e.g., \cite{NV, thm.~4.4.4$(v)$, thm.~9.6.8}. Thus it is enough to prove
that for any irreducible graded module $M$ there is an integer 
$d$ such that $M[d]$ is $\flat$-selfdual.
This is proved as in \cite{KL, p.~342}.
More precisely, by definition of the duality functor $\flat$,
for any graded module $M$ in $\th\Rb$-$\fmodb$ we have
$$\ch(M^\flat)=\ch(M)\ \mod\ (v-1),\quad
\gdim(1_\ib M^\flat)=\overline{\gdim(1_\ib M)}.$$ 
Thus if $M$ is irreducible then we have $M^\flat=M[d]$ for some integer $d$.
We must prove that $d$ is even. 
It is enough to prove the following.

\proclaim{8.28.~Lemma} If $M\in\th\Rb$-$\fmodb$ 
is irreducible then for each $\ib$ we have
$$\gdim(1_\ib M)\in v\ZZ[v^2,v^{-2}]\cup\ZZ[v^2,v^{-2}].$$
\endproclaim

\noindent{\sl Proof :}
Indeed, we'll prove that this identity holds for 
the projective module $M=\th\Rb_{\jb}$ where $\jb$ is any sequence
in $\th\!I^m$. This implies our claim.
Set $\jb=(j_{1-m},\dots,j_{m-1},j_m)$. Proposition 8.13 yields
$$\ch(\th\Rb_\jb)=\ch(\th\Rb_{j_0j_1})\circledast\ch(\Rb_{j_2})\circledast
\cdots\circledast\ch(\Rb_{j_m}).$$
Examples 5.9$(b)$, 7.5$(a)$ yield
$$\gathered
\ch(\th\Rb_{j_0j_1})=(1-v^2)^{-1}
\bigl(j_0j_1+v^{\l_{j_0}+\l_{j_1}}j_1j_0\bigr),
\quad
\ch(\Rb_{j_k})=(1-v^2)^{-1}\,j_k.\endgathered$$
So, by Definition 8.12$(b)$, it is enough to check that for each reflection
$w$ of $W_m$ which fixes the sequence $\jb$ 
the degree of $\sigma_{\dot w}1_\jb$ is even.
This reduces to the following computation (left to the reader). 
Fix $k\neq l$ such that
$j_k=j_l$. 
If one of the following holds 
$$\gathered
1\leqslant k<l,\quad
\dot w=s_{k}\dots s_{l-2}s_{l-1}s_{l-2}\dots s_{k},
\vspace{2mm}
1\leqslant 1-k<l,\quad
\dot w=s_{-k}\cdots s_1\eps_1s_1\dots s_{l-1}\dots s_{1}
\eps_1s_{1}\dots s_{-k},\endgathered$$
then
$\deg(\sigma_{\dot w}1_\jb)$ is even.

\qed

\vskip3mm

\subhead 8.29.~The algebra $\th\Bcb$ and 
its representation in $\th\Vb(\l)$\endsubhead 
Following \cite{EK1,2,3} we define a $\Kc$-algebra $\th\Bcb$
as follows.

\proclaim{8.30.~Definition} Let $\th\Bcb$ be the $\Kc$-algebra
generated by $e_i$, $f_i$ and invertible elements $t_i$, $i\in I$,
satisfying the following defining relations

\vskip2mm
\itemitem{$(a)$} 
$t_it_j=t_jt_i$ and  $t_{\theta(i)}=t_i$ for all $i,j$,

\vskip2mm
\itemitem{$(b)$} 
$t_ie_jt_i^{-1}=v^{i\cdot j+\theta(i)\cdot j}e_j$ 
and $t_if_jt_i^{-1}=
v^{-i\cdot j-\theta(i)\cdot j}f_j$ for all $i,j$, 

\vskip2mm
\itemitem{$(c)$}
$e_if_j=v^{-i\cdot j}f_je_i+
\delta_{i,j}+\delta_{\theta(i),j}t_i$ for all $i,j$,

\vskip2mm
\itemitem{$(d)$}
formula $(7.3)$ holds with $\theta_i=e_i$, or with $\theta_i=f_i$.
\endproclaim

\noindent
We define a representation of
$\th\Bcb$ as follows.
The $\Kc$-vector space $$\th\Vb(\l)=\Kc\otimes_\Ac\th\Kb_I$$
is equipped with $\Kc$-linear operators $e_i$, $e'_i$, $f_i$
and with a $\Kc$-bilinear form in the obvious way. 
Let $\phi_\l$ be the class of $\kb$ in $\th\Kb_I$, where $\kb$ is regarded
as the trivial module over the $\kb$-algebra $\th\Rb_{0}$.
Let $\l$ be as in (6.2).
We can now prove the following theorem, which is the main result of Section 8.

\proclaim{8.31.~Theorem}
(a) The operators $e_i$, $f_i$ 
define a representation of $\th\Bcb$ on $\th\Vb(\l)$.
The $\th\Bcb$-module
$\th\Vb(\l)$ is irreducible and 
for $i\in I$ we have
$$e_i\phi_\l=0,\quad
t_i\phi_\l=v^{\l_i+\l_{\theta(i)}}\phi_\l,\quad
\{x\in \th\Vb(\l);\,e_jx=0,\,\forall j\}=\Kc\, \phi_\l.$$

(b) There is a unique symmetric $\Kc$-bilinear form $(\bullet :\bullet)$
on $\th\Vb(\l)$ such that $(\phi_\l:\phi_\l)=1$ and
$(e'_ix:y)=(x:f_iy)$ for all $i\in I$, $x,y\in\th\Vb(\l)$,
and it is non-degenerate.

\vskip2mm

(c) There is a unique $\Kc$-antilinear map
$\th\Vb(\l)\to\th\Vb(\l)$ such that $P\mapsto P^\sharp$
for all graded projective module $P$.
It is the unique $\Kc$-antilinear map such that 
$\phi_\l^\sharp=\phi_\l$ and
$(f_ix)^\sharp=f_i (x^\sharp)$ for all $x\in\th\Vb(\l)$.
\endproclaim

\noindent{\sl Proof :}
For each $i$ in $I$ we define the $\Ac$-linear operator $t_i$ 
on $\th\Kb_I$ by setting
$$t_i P=v^{\l_i+\l_{\theta(i)}-\nu\cdot(i+\theta(i))}P,\quad
\forall P\in\th\Rb_{\nu}\text{-}\proj.$$
We must prove that the operators $e_i$, $f_i$, and $t_i$ satisfy the
relations in Definition 8.30. Relation $(c)$ is the only non trivial one, see Lemma 8.9$(e)$. 
To check it we need a version of the Mackey's induction-restriction theorem.
Note that we have
$$\gathered D_{m,1;m,1}=\{e,s_{m},\eps_{m+1}\},\vspace{2mm}
W(e)=W_{m,1},\quad W(s_{m})=W_{m-1,1,1},\quad
W(\eps_{m+1})=W_{m,1}.\endgathered$$

\proclaim{8.32.~Lemma} Fix $i$, $j$ in $I$. Let $\mu$, $\nu$ 
in $\th\NN I$ be such that
$\nu+i+\theta(i)=\mu+j+\theta(j)$.
Put $m=|\nu|/2=|\mu|/2$.
The graded 
$(\th\Rb_{m,1},\th\Rb_{m,1})$-bimodule 
$1_{\nu,i}\th\Rb_{m+1}1_{\mu,j}$ has a filtration by graded bimodules 
whose associated graded is isomorphic to :

\vskip2mm

\itemitem{(a)} 
$\bigl(\th\Rb_\nu\otimes\Rb_i\bigr)\oplus
\Bigl((\th\Rb_m1_{\nu',i}\otimes \Rb_{i})\otimes_{\Rb'}
(1_{\nu',i}\th\Rb_m\otimes \Rb_i)\Bigr)[d]$ if $j=i$,

\vskip2mm

\itemitem{(b)} 
$\bigl(\th\Rb_\nu\otimes\Rb_{\theta(i)}\bigr)[d']\oplus
\Bigl((\th\Rb_m1_{\nu',\theta(i)}\otimes \Rb_{i})\otimes_{\Rb'}
(1_{\nu',i}\th\Rb_m\otimes \Rb_{\theta(i)})\Bigr)[d]$ if $j=\theta(i)$,

\vskip2mm

\itemitem{(c)} 
$\Bigl((\th\Rb_m1_{\nu',j}\otimes \Rb_{i})\otimes_{\Rb'}
(1_{\nu',i}\th\Rb_m\otimes \Rb_j)\Bigr)[d]$ if $j\neq i,\theta(i)$.

\vskip2mm

\noindent
Here we have set $\nu'=\nu-j-\theta(j)$,
$\Rb'=\th\Rb_{m-1,1,1}$,
$d=\deg(\sigma_{m}1_{\nu',i,j})$,
and $d'=\deg(\pi_{m+1}1_{\nu,\theta(i)})$

\endproclaim

\noindent
The proof is similar to the proof of \cite{M, thm.~1}, \cite{KL, prop.~2.18}. 
It is left to the reader.
Note that we have the following formulas, see Remark 5.2,
$$\deg(\pi_{m+1}1_{\nu,\theta(i)})=
\l_i+\l_{\theta(i)}-\nu\cdot(i+\theta(i))/2,
\quad
\deg(\sigma_{m}1_{\nu',i,j})=-i\cdot j.$$
Now, recall that $P$ lie in $\th\Rb_{\nu}\text{-}\proj$ and that
$$f_j(P)=\th\Rb_{m+1}1_{m,j}\otimes_{\th\Rb_{m,1}}(P\otimes\Rb_1),\quad
e'_i(P)=1_{m-1,i}P,$$
where $1_{m-1,i}P$ is regarded as a $\th\Rb_{m-1}$-module. 
Therefore we have
$$\gathered
e'_if_j(P)=1_{m,i}\th\Rb_{m+1}1_{m,j}\otimes_{\th\Rb_{m,1}}(P\otimes\Rb_1),
\vspace{2mm}
f_je'_i(P)=\th\Rb_{m}1_{m-1,j}\otimes_{\th\Rb_{m-1,1}
}(1_{m-1,i}P\otimes\Rb_1).
\endgathered$$
Therefore we have the following identities

\vskip2mm

\itemitem{$\bullet$}  $e'_if_i(P)=P\otimes\Rb_i+f_ie'_i(P)[-2]$,

\vskip2mm

\itemitem{$\bullet$}  
$e'_if_{\theta(i)}(P)=P\otimes\Rb_{\theta(i)}
[\l_i+\l_{\theta(i)}-\nu\cdot(i+\theta(i))/2]
+f_{\theta(i)}e'_i(P)[-i\cdot\theta(i)]$,

\vskip2mm

\itemitem{$\bullet$}  $e'_if_j(P)=f_je'_i(P)[-i\cdot j]$ if $i\neq j,\theta(j)$.
\vskip2mm

\noindent
Note that Lemma 8.32 implies these relations up to some filtration. Hence,
since the associated graded is projective, they hold in full generality.
This proves the first claim of part $(a)$ of the theorem.
Next, recall the following fact, see
\cite{EK1, prop.~2.5},
\cite{EK3, prop.~2.11}.

\proclaim{8.33.~Claim}
There is a $\th\Bcb$-module
generated by a non-zero vector $\phi_\l$ such that
$$e_i\phi_\l=0,\quad
t_i\phi_\l=v^{\l_i+\l_{\theta(i)}}\phi_\l, \quad
\{x;\,e_jx=0,\,\forall j\}=\Kc\, \phi_\l,\quad
i\in I.$$
This $\th\Bcb$-module is irreducible and it is unique up to an isomorphism.
\endproclaim

\noindent
So we must check that 
the $\th\Bcb$-module $\th\Vb(\l)$ 
satisfies the axioms above.
It is generated by $\phi_\l$ 
by Lemma 8.34 below. The other axioms are obvious.

Part $(b)$ of the theorem 
follows from \cite{EK2, prop.~4.2(ii)} and Lemma 8.9$(b)$.
The bilinear form $(\bullet:\bullet)$ is the same as the bilinear form
obtained from (8.2) by base change from $\Bc$ to $\Kc$.

Finally, for part $(c)$ of the theorem it is enough to check that
$(f_iP)^\sharp=f_i(P^\sharp)$ 
for any graded module $P$ in $\th\Rb$-$\proj$.
By Lemma 8.34 below we may assume that
$P=\th\Rb_\yb$ for some $\yb$.
By (8.5) we can also assume that
$\yb=\ib$ with $\ib\in\th\! I^\nu$.
Then the claim follows from the formulas in Proposition 8.14,
because $\th\Rb_\ib$ is $\sharp$-selfdual for any $\ib$, see Section 8.10.

\qed

\vskip3mm

\proclaim{8.34.~Lemma} 
Any object of $\th\Rb_m$-$\proj$ 
is of the form
$\th\Rb_m\otimes_{\Rb_m}P$ 
for some $P$ in $\Rb_m$-$\proj$.
The $\Ac$-module $\th\Kb_{I,m}$ is spanned by the 
$\th\Rb_\yb$'s with
$\yb$ in $\th Y^m$.
\endproclaim

\noindent{\sl Proof :}
Let $b$ be a simple object  in $\th\Rb_\nu$-$\fModb_0$ with $|\nu|=m$.
We'll view it as an element of $\th\!B(\l)$.
An easy induction using Proposition 8.23 implies that for each integer $a\geqslant 1$ we have
$$f_i^{(a)} \th\!G^\low(b)=
\left\langle\matrix\eps_i(b)+a\cr a\endmatrix\right\rangle\,
\th\!G^\low(\tilde F^a_ib)+
\sum_{b'}f_{b,b'}\th\!G^\low(b'),$$ where $b'$ runs over the set of elements
of $\th\!B(\l)$ such that $\eps_i(b')>\eps_i(b)+a$ and 
$f_{b,b'}$ lies in $\Ac$. Therefore, for any
$i$ such that $\eps_i(b)\geqslant 1$, we have
$$f_i^{(\eps_i(b))} \th\!G^\low(\tilde E_i^{(\eps_i(b))}b)=
\th\!G^\low(b)+
\sum_{b'}f_{b,b'}\th\!G^\low(b')$$
by Proposition 8.22$(c),(e)$.
Here $b'$ runs over the set of elements
of $\th\!B(\l)$ such that $\eps_i(b')>\eps_i(b)$.
Thus a simultaneous induction on $\nu$ and descending induction on $\eps_i(b)$ implies that
$\th\!G^\low(b)$ lies in the $\Ac$-span of the elements $f_\yb(\kb)$ with $\yb\in Y^m$.
We are done, because $f_\yb(\kb)=\th\Rb_{\theta(\yb)\yb}$, see Section 8.10.

\qed

\vskip3mm

\subhead 8.35.~Remark\endsubhead
The $\th\Bcb$-module $\th\Vb(\l)$ is the same as the $\th\Bcb$-module
$V_\theta(\l+\theta(\l))$ in \cite{EK1, prop.~2.5}.
Let $(\bullet:\bullet)_{_{KE}}$ be the bilinear form on $\th\Vb(\l)$ 
considered in loc.~cit. We have
$$(P:Q)=(1-v^2)^{-m}(P:Q)_{_{KE}},\quad\forall P,Q\in\th\Rb_m\text{-}\proj.$$ 
Note that 
$(\bullet:\bullet)_{_{KE}}$ 
is a symmetric $\Ac$-bilinear form  
$\th\Kb_I\times \th\Kb_I\to \Ac$, and that Theorem 8.31$(b)$ yields
$$(e_ix:y)_{_{KE}}=(x:f_iy)_{_{KE}},\quad i\in I,\,x,y\in \th\Vb(\l).$$

\vskip3mm

\subhead 8.36.~Results over an arbitrary field $\kb$\endsubhead
Recall that $p,q\in\kb^\times$ and that $I$ is a $\ZZ\rtimes\ZZ_2$-invariant subset of $\kb^\times$.
We associate to $I$ a quiver with an involution $(\Gamma,\theta)$ as in Section 6.2.
Fix $\l\in\NN I$ as in (6.2), i.e., we set
$$\l=\sum_{i}i,\quad i\in I\cap\{q,-q\}.$$
The graded $\kb$-algebra $\th\Rb_m$, defined in Sections 5.1, 6.4, 
and the operators $e_i$, $f_i$ on
$\th\Kb_I$, defined in (8.4), make sense over any field $\kb$ (not necessarily algebraically closed
nor of characteristic zero). For any $\kb$ there is again a
$\th\Bcb$-module isomorphism $$\th\Vb(\l)=\Kc\otimes_\Ac\th\Kb_I,$$ where $\th\Vb(\l)$ is 
the Enomoto-Kashiwara's $\th\Bcb$-module. To prove this it is enough to check the axioms in 
Claim 8.33. The proof is the same as in characteristic zero.
Note that the $\Kc$-algebra $\th\Bb$ and the $\th\Bb$-module $\th\Vb(\l)$ depend only on $(\Gamma,\theta)$
(i.e., on $I$ and $p$), and on $\l$ (i.e., on $q$).
Therefore, for each $m$, the number of simple graded
$\th\Rb_m$-modules is the same for any field $\kb$ as long as $\Gamma,\theta,\l$ remain unchanged. 
In particular all simple graded
$\th\Rb_m$-modules are absolutely irreducible. Recall that the simple graded
$\th\Rb_m$-modules are finite dimensional, because $\th\Rb_m$ is finitely generated
over its center. Therefore all simple graded
$\th\Rb_m$-modules are {\it split simple}, i.e., with a one dimensional endomorphism
$\kb$-algebra, see e.g, \cite{L,thm.~7.5}. Note that, for $\kb$ is algebraically closed,
we already use this when claiming that the Cartan pairing is perfect.

The discussion above, Theorem 6.5, and Remark 6.10 imply that the number of simple graded
$\Hb_m$-modules in $\Modb_I$ is the same for any field $\kb$ of characteristic $\neq 2$ as long as 
$\Gamma,\theta,\l$ remain unchanged. 
 
\vskip2cm

\head 9.~Presentation of the graded $\kb$-algebra $\th\Zb_{\Lbb,\Vb}^\delta$
\endhead

Fix a quiver $\Gamma$ with set of vertices $I$ and set of arrows $H$.
Fix an involution $\theta$ on $\Gamma$. 
Assume that $\Gamma$ has no 1-loops and that
$\theta$ has no fixed points. Fix a
dimension vector $\nu$ in $\th\NN I$ and a
dimension vector $\l$ in $\NN I$. Set $|\nu|=2m$.
Fix an object $(\Vb,\varpi)$ in $\th\Vcb_\nu$ 
and an object $\Lbb$ in $\Vcb_\l$.
In this section we give a proof of Theorem 5.8.
By Theorem 4.17 and Corollary 5.6 there is a unique injective graded
$\k$-algebra homomorphism
$$\gathered
\Phi:\th\Rb(\Gamma)_{\l,\nu}\to\th\Zb^\delta_{\Lbb,\Vb},
\vspace{2mm}
1_\ib\mapsto1_{\Lbb,\Vb,\ib}, \quad
\varkappa_{\ib,l}\mapsto\varkappa_{\Lbb,\Vb,\ib}(l), \quad
\sigma_{\ib,k}\mapsto\sigma_{\Lbb,\Vb,\ib}(k), \quad
\pi_{\ib,1}\mapsto\pi_{\Lbb,\Vb,\ib}(1),
\vspace{2mm}\ib\in\th\! I^\nu,\quad k=1,\dots,m-1,\quad l=1,2,\dots,m.
\endgathered$$ We must prove
that $\Phi$ is a surjective map. 
Note that both algebras have 1, because the set $\th\! I^\nu$ is finite.
Since the grading does not matter, we can replace
$\th\Zb_{\Lbb,\Vb}^\delta$ by
$\th\Zb_{\Lbb,\Vb}$.
To unburden the notation we abbreviate
$$\gathered
\th\Rb_\nu=\th\Rb(\Gamma)_{\l,\nu},\quad
\th\Fb_\nu=\th\Fb_{\Lbb,\Vb},
\quad\th\Zb_\nu=\th\Zb_{\Lbb,\Vb}, 
\vspace{2mm}
\varkappa_{\nu,\ib}(l)=\varkappa_{\Lbb,\Vb,\ib}(l),\quad
\sigma_{\nu,\ib}(k)=\sigma_{\Lbb,\Vb,\ib}(k), \quad
\pi_{\nu,\ib}(1)=\pi_{\Lbb,\Vb,\ib}(1),
\vspace{2mm}
\th\!\Zc_\nu=\th\!\Zc_{\Lbb,\Vb},
\quad\th\! Z_\nu=\th\! Z_{\Lbb,\Vb},\quad etc.
\endgathered$$ 

\vskip2mm

\subhead 9.1.~The filtration of $\th\Zb_{\nu}$\endsubhead
Recall that $W_m$ is regarded as a Coxeter group of type $\B_m$ 
with the set of simple reflections $\{s_1,s_2,\dots,s_m\}$.
From now on let $\leqslant$ and $\ell$ be the 
corresponding Bruhat order and length function.
For a future use we set also
$$\gathered
\Delta^+=\Delta^+_s\sqcup\Delta^+_l,\quad
\Delta^+_s=\{\chi_k\pm \chi_l;\,1\leqslant l<k\leqslant m\},\quad\Delta_l^+=\{\chi_1,\chi_2,\dots,\chi_m\}.
\endgathered$$
Note that $\leqslant$, $\ell$, and $\Delta^+$ differ from the order, the length function,
and the set of positive roots introduced in Section 4.2. We hope this will not create any confusion.
We can now introduce a filtration of
$\th\! Z_{\nu}$ by closed subsets. 
We define 
$$
\th\! O_\nu^{\leqslant x}=\bigcup_{w\leqslant x}\th\! O_\nu^w,\quad
\th\! Z_\nu^{\leqslant x}=q^{-1}(\th\!O^{\leqslant x}_\nu),\quad
\th\!\Zc_\nu^{\leqslant x}=
H_*^{\th\! G_\nu}(\th\! Z_\nu^{\leqslant x},\k).
$$

\proclaim{9.2.~Lemma}
(a) 
The set $\th\! Z_{\nu}^{\leqslant x}$ is closed.
The variety $\th\! Z_{\nu}^{x}$ is smooth if $\ell(x)=1$.

\vskip1mm
(b) The direct image by the inclusion 
$\th\! Z_{\nu}^{\leqslant x}\subset\th\! Z_{\nu}$ is an injection
$\th\!\Zc_{\nu}^{\leqslant x}\subset\th\!\Zc_{\nu}$.
\vskip1mm
(c) The
convolution product maps
$\th\!\Zc_{\nu}^{\leqslant x}\times\th\!\Zc_{\nu}^{\leqslant y}$ into
$\th\!\Zc_{\nu}^{\leqslant xy}$ for each $x$, $y$ such that 
$\ell(xy)=\ell(x)+\ell(y)$. 
\vskip1mm
(d) The unit of $\th\!\Zc_\nu$ lies in 
$\th\!\Zc_{\nu}^{e}$.
\endproclaim

\noindent{\sl Proof :}
To avoid confusions, let $\leqslant_D$ and $\ell_D$ be the 
Bruhat order and length function 
introduced in Section 4.2. The claims in the lemma are well-known if
we replace $\leqslant$, $\ell$ by $\leqslant_D$, $\ell_D$ respectively.
Therefore, it is enough to prove the following : assume that
$\th\! O^{v}_{\nu,x,y}$ and $\th\!O^{w}_{\nu,x,y}$ are non empty.
Then $\th\!O^{v}_{\nu,x,y}\subset\th\!\bar O^{w}_{\nu,x,y}$ iff
$v\leqslant w$. Up to the action of a well-chosen diagonal element we may assume that $x=e$.
We can also assume that $y$ is minimal in the coset $W_\nu y$.
Since $\th\! O^{v}_{\nu,e,y}$ and $\th\!O^{w}_{\nu,e,y}$ are non empty, we have
$v,w\in W_\nu y$.
Finally, on the coset $W_\nu y$
the orders $\leqslant$, $\leqslant_D$ are the same because
$W_\nu\subset\Sen_m$, see (4.2) and the last remark in Section 4.2.

\qed

\vskip3mm

Let $\th\Zb^{\leqslant x}_\nu$ be the image of
$\th\!\Zc^{\leqslant x}_\nu$ by the isomorphism 
$\th\Zb_\nu=\th\!\Zc_\nu$ in Proposition 3.1$(b)$. 

\vskip3mm
\subhead 9.3.~PBW theorem for $\th\Zb_{\nu}$\endsubhead 
Recall that 
$$\th\Fb_{\nu}=\bigoplus_{\ib\in\th\!I^\nu}
\kb[x_\ib(1),x_\ib(2),\dots,x_\ib(m)],$$ 
see Section 4.11. The graded $\kb$-algebra
$\th\Zb_{\nu}$ has a natural structure of graded
$\th\Fb_{\nu}$-module such that 
$x_\ib(l)$  acts  as  $\varkappa_{\nu,\ib}(l)$ 
under the inclusion $\th\Zb_\nu\subset\End(\th\Fb_\nu)$ in Theorem 4.17.  
Recall that $\th\Zb_{\nu}=\th\!\Zc_{\nu}$. 
The following is immediate, see e.g., \cite{CG, sec.~5.5}.

\proclaim{9.4.~Lemma} We have $\th\Zb_{\nu}^{\leqslant x}=
\bigoplus_{w\leqslant x}\th\Fb_{\nu}\,[\th\! Z_{\nu}^{w}]$ for each $x$.
In particular $\th\Zb_{\nu}$ is a free graded
$\th\Fb_{\nu}$-module of rank $2^mm!$.
\endproclaim

\vskip3mm

The map $\Phi$ is a 
graded $\th\Fb_\nu$-module homomorphism. For each $x$ we set
$$\th\Rb^{\leqslant x}_\nu=\sum_{w\leqslant x}
\th\Fb_\nu\,1_\nu\sigma_{\dot w},$$ 
where $\sigma_{\dot w}$ is defined as in (5.2).
It is a graded $\th\Fb_\nu$-submodule
of $\th\Rb_{\nu}$. 
We abbreviate $\th\Rb^{e}_\nu=\th\Rb^{\leqslant e}_\nu$. The
proof of the surjectivity of $\Phi$ consists of two steps. First we prove that
$\Phi(\th\Rb^{\leqslant x}_\nu)\subset\th\Zb_{\nu}^{\leqslant x}$. Then
we prove that this inclusion is an equality.

\vskip3mm

\subhead 9.5.~Step 1\endsubhead
Since $\Phi$ is a $\th\Fb_\nu$-module
homomorphism it is enough to prove that the element
$\Phi(\sigma_{\dot x})$ lies in $\th\Zb_{\nu}^{\leqslant x}$. By
an easy induction on the length of $x$ it is
enough to observe that we have
$$\gathered
\Phi(1)\in\th\Zb_{\nu}^{e},\quad
\Phi(\sigma_{k})\subset\th\Zb_{\nu}^{\leqslant
s_k},\quad  k=1,2,\dots,m.\endgathered$$ 
This follows from the definition of  the elements
$\sigma_\nu(k)$, $\pi_\nu(1)$ of $\th\Zb_{\nu}$
in Section 4.14. Recall that $s_m=\eps_1$ and $\sigma_m=\pi_1$, see (5.2) for details.

\vskip3mm

\subhead 9.6.~Step 2\endsubhead
Note that $\th\Zb^e_{\nu}$ is the free
$\th\Fb_\nu$-module of rank one generated by $[\th\! Z^e_{\nu}]$. 
Therefore we have
$$\Phi(\th\Rb^{e}_\nu)=
\th\Fb_\nu\,[\th\!Z_\nu^e]= \th\Zb_\nu^{e}.$$ 
To complete the proof of Step 2 we are reduced to prove
the following.

\proclaim{9.7.~Lemma} If $\ell(s_kw)=\ell(w)+1$ and $k=1,2,\dots m$,
then we have the following formula in 
$\th\Zb_\nu^{\leqslant s_kw}/\th\Zb_\nu^{<s_kw}$ :
$$[\th\!Z^{s_k}_\nu]\star[\th\!Z^w_\nu]=[\th\!Z^{s_kw}_\nu].$$
\endproclaim

\noindent{\sl Proof :} 
By Lemmas 9.2$(c)$ and 9.4 there
is an unique element $c$ in $\th\Fb_\nu$ such that
$$[\th\! Z^{s_k}_\nu]\star[\th\! Z^w_\nu]=
c\star[\th\!Z^{s_kw}_\nu]\ \roman{in}\ \th\Zb_\nu^{\leqslant
s_kw}/\th\Zb_\nu^{<s_kw}.$$ Let us prove that $c=1$.
For each $x,y,z$ there is a
unique element $\Lambda^x_{y,z}$ in $\Qb$ such that
$$[\th\!Z^x_\nu]=\sum_{y,z}\Lambda^x_{y,z}\psi_{y,z},$$
see Section 4.12.
Since $\phi_{\nu,y,yx}$ is a smooth point of
$\th\!Z^x_\nu$ we have also
$$\Lambda^x_{y,yx}=\Eu(\th\!Z^x_\nu,\phi_{\nu,y,yx})^{-1}.$$
Hence, in the expansion of the element $[\th\!Z^{s_kw}_\nu]$ in the
$\Qb$-basis $(\psi_{y,z})$  the coordinate along the vector
$\psi_{x,xs_kw}$ is equal to
$$\Lambda_{x,xs_kw}^{s_kw}=\Eu(\th\!Z^{s_kw}_\nu,\phi_{\nu,x,xs_kw})^{-1}.$$
On the other hand, since $\Lambda^w_{x,xs_kw}=0$ and
$$[\th\!Z^{s_k}_\nu]=\sum_{x}(\Lambda^{s_k}_{x,x}\psi_{x,x}+
\Lambda^{s_k}_{x,xs_k}\psi_{x,xs_k}),$$
the coordinate of $[\th\!Z^{s_k}_\nu]\star[\th\!Z^w_\nu]$ along 
$\psi_{x,xs_kw}$ is equal to
$$\Lambda_{x,xs_k}^{s_k}\Lambda^w_{xs_k,xs_kw}\Lambda_{xs_k}=
\Eu(\th\!Z^{s_k}_\nu,\phi_{\nu,x,xs_k})^{-1}
\Eu(\th\!Z^{w}_\nu,\phi_{\nu,xs_k,xs_kw})^{-1}\Lambda_{xs_k}.$$ Thus we must
check that
$$\Eu(\th\!Z^{s_k}_\nu,\phi_{\nu,x,xs_k})\,
\Eu(\th\!Z^{w}_\nu,\phi_{\nu,xs_k,xs_kw})=
\Eu(\th\!Z^{s_kw}_\nu,\phi_{\nu,x,xs_kw})\, \Lambda_{xs_k}.$$ 
This follows from the lemma below.

\proclaim{9.8.~Lemma} (a) For $x,y\in W$ we have
$$\aligned
\eu(\th\!O^y_\nu,\phi_{\nu,x,xy})&=
\eu(\th\!\nen_{\nu,x}\oplus\th\!\men_{\nu,xy,x}),
\hfill\vspace{2mm}
\Eu(\th\!Z^{y}_\nu,\phi_{\nu,x,xy})&=
\Eu(\th\!O^{y}_\nu,\phi_{\nu,x,xy})\,\Eu(\th\!\een_{\nu,x,xy}^*), 
\hfill\vspace{2mm}
\Lambda_x=\Eu(\th\!Z^{e}_\nu,\phi_{\nu,x,x})
&=\Eu(\th\! F_\nu,\phi_{\nu,x})\,\Eu(\th\!\een_{\nu,x}^*).
\endaligned$$

(b) For $w,x,y\in W$ such that $\ell(xy)=\ell(x)+\ell(y)$ we
have
$$\aligned
\eu(\th\!O_\nu^{xy},\phi_{\nu,w,wxy})\,\eu(\th\!F_\nu,\phi_{\nu,wx})&=
\eu(\th\!O^x_\nu,\phi_{\nu,w,wx})\,\eu(\th\!O_\nu^{y},\phi_{\nu,wx,wxy}),
\hfill\vspace{2mm}
\Eu(\th\!\een^*_{\nu,w,wxy}\oplus\th\! \een^*_{\nu,wx})&=
\Eu(\th\!\een^*_{\nu,w,wx}\oplus\th\! \een^*_{\nu,wx,wxy}).
\endaligned$$
\endproclaim

\noindent{\sl Proof :} Part $(a)$ is left to the reader. 
Compare Proposition 4.13 where similar formulas are given. Let us
prove $(b)$. Clearly we can assume $w=e$.
Set $$\Delta(y)^-=y(\Delta^+)\cap\Delta^-,
\quad\Delta(y)^+=y(\Delta^-)\cap\Delta^+.$$ 
Let the symbol $\sqcup$ denote a disjoint union. Recall that
$$\ell(xy)=\ell(x)+\ell(y)\Rightarrow
\cases
\Delta(xy)^-=\Delta(x)^-\sqcup x(\Delta(y)^-),\vspace{2mm}
\Delta(xy)^+=\Delta(x)^+\sqcup x(\Delta(y)^+).\endcases\leqno(9.1)$$
For $x,y\in W$
the $T$-module $\th\!\men_{\nu,xy,x}$ is the sum of the root
subspaces whose weight belong to the set
$x(\Delta(y)^-)\cap\th\!\Delta_\nu$,
and the $T$-module $\th\!\nen_{\nu,x}$ is the sum of the root
subspaces whose weight belong to the set
$x(\Delta^+)\cap\th\!\Delta_\nu$, see Sections 4.8, 4.9 for details.
Thus, by $(a)$, the first claim follows from the following equality
$$\Delta^+\sqcup\Delta(xy)^-\sqcup x(\Delta^+)=
\Delta^+\sqcup\Delta(x)^-\sqcup x(\Delta^+)\sqcup x(\Delta(y)^-).
$$
This equality is a consequence of the first identity in (9.1).
Now, let us concentrate on the second claim. 
Set
$$\gathered
S_{x,xy}=x(\Delta^+_s\cap y(\Delta^+_s)),
\quad
L_{x,xy}=x(\Delta^+_l\cap y(\Delta_l^+)).
\endgathered$$
There are integers $h_\a\geqslant 0$ such that the character of the
$T$-modules $\th\!E_\Vb$, $\th\!\een_{\nu,x}$ and $\th\!\een_{\nu,x,xy}$ are of the following form
$$\gathered
\ch(\th\!E_\Vb)=\sum_{\a}h_\a\,\a,\quad\a\in\Delta,
\vspace{2mm}
\ch(\th\!\een_{\nu,x})=\sum_\a h_\a\a+\sum_l\l_{i_l}\chi_l,\quad
\a\in x(\Delta^+_s),\quad
\chi_l\in x(\Delta^+_l),\vspace{2mm}
\ch(\th\!\een_{\nu,x,xy})=\sum_\a h_\a\a+\sum_l\l_{i_l}\chi_l,\quad
\a\in S_{x,xy},\quad
\chi_l\in L_{x,xy}.
\endgathered$$
See Section 4.9 for details. 
Let $$\gathered
S=S_{e,xy}\sqcup x(\Delta^+_s),\quad
S'=S_{e,x}\sqcup S_{x,xy},
\vspace{2mm}
L=L_{e,xy}\sqcup x(\Delta^+_l),\quad
L'=L_{e,x}\sqcup L_{x,xy}.
\endgathered$$ We obtain
$$\gathered
\ch(\th\!\een_{\nu,e,xy}\oplus\th\! \een_{\nu,x})=
\sum_\a h_\a\a+\sum_l\l_{i_l}\chi_l,\quad
\a\in S,\quad
\chi_l\in L,
\vspace{2mm}
\ch(\th\!\een_{\nu,e,x}\oplus
\th\!\een_{\nu,x,xy})=\sum_\a h_\a\a+\sum_l\l_{i_l}\chi_l,\quad
\a\in S',\quad
\chi_l\in L'.
\endgathered$$
Now, a short computation yields 
$$\aligned
S=S' \iff\Delta^+_s\cap\Delta(xy)^+=
\Delta^+_s\cap\bigl(\Delta(x)^+\sqcup x(\Delta(y)^+)\bigr),\vspace{2mm}
L=L' \iff\Delta^+_l\cap\Delta(xy)^+=
\Delta^+_l\cap\bigl(\Delta(x)^+\sqcup x(\Delta(y)^+)\bigr).
\endaligned$$
Thus the claim follows from (9.1).
\qed

\vskip2cm

\head 10.~Perverse sheaves on $\th\Eb_{\Lbb,\Vb}$ and the global bases of $\th\Vb(\l)$\endhead

Fix a quiver $\Gamma$ with set of vertices $I$ and set of arrows $H$.
Fix an involution $\theta$ on $\Gamma$. 
Assume that $\Gamma$ has no 1-loops and that
$\theta$ has no fixed points. 
Fix a dimension vector $\nu$ in $\th\NN I$ and a
dimension vector $\l$ in $\NN I$. Set $|\nu|=2m$.
Fix an object $(\Vb,\varpi)$ in $\th\Vcb_\nu$ 
and an object $\Lbb$ in $\Vcb_\l$.
To unburden the notation we'll abbreviate
\vskip1mm
$$\gathered
\th\!G_\nu=\th\!G_\Vb,\quad
\th\!\Rb_\nu=\th\!\Rb(\Gamma)_{\l,\nu}.
\endgathered$$

\vskip3mm

\subhead 10.1.~Perverse sheaves on $\th\! E_{\nu}$\endsubhead
First, we generalize the setting in Section 2. We define 
{\it an orientation $\Omega$
of $I$} to be a partition $I=\Omega\sqcup\bar\Omega$. Fix an orientation $\Omega$.
For each $I$-graded $\CC$-vector space $\Wb$ we write 
$\Wb_\Omega=\bigoplus_{i\in\Omega}\Wb_i$. Now, we define
$$L_{\Lbb,\Vb,\Omega}=\Hom_\Vcb(\Lbb_\Omega,\Vb_\Omega)\oplus
\Hom_\Vcb(\Vb_{\bar\Omega},\Lbb_{\bar\Omega}),\quad
\th\!E_{\Lbb,\Vb,\Omega}=\th\!E_\Vb\times L_{\Lbb,\Vb,\Omega}.
$$
An element of $\th\!E_{\Lbb,\Vb,\Omega}$ is a triplet
$(x,y,z)$ with 
$$x\in\th\!E_\Vb,\quad
y\in\Hom_\Vcb(\Lbb_\Omega,\Vb_\Omega),\quad
z\in\Hom_\Vcb(\Vb_{\bar\Omega},\Lbb_{\bar\Omega}).$$
For each $\yb$ in $\th Y^\nu$ we define also
$$\gathered
\th\!\widetilde F_{\Lbb,\Vb,\yb,\Omega}=
\{(x,y,z,\phi)\in\th\!E_{\Lbb,\Vb,\Omega}\times\th\!F_{\Vb,\yb}\,;\,
\phi=(\Vb^l),\,
x(\Vb^l)\subset\Vb^l,\,\vspace{2mm}
y(\Lbb)\subset\Vb^0,\,
z(\Vb^0)=0\}.\endgathered$$
To unburden the notation we'll abbreviate
\vskip1mm
$$\gathered
\vspace{2mm}
\th\!E_{\nu,\Omega}=\th\!E_{\Lbb,\Vb,\Omega},\quad
\th\!F_\yb=\th\!F_{\Vb,\yb},\quad
\th\! \widetilde F_{\yb,\Omega}=\th\! \widetilde F_{\Lbb,\Vb,\yb,\Omega}.
\endgathered$$
We define the semisimple complex
$\th\!\Lc_{\yb,\Omega}$ over $\th\! E_{\nu,\Omega}$ as the direct image
of the constant sheaf over $\th\! \widetilde F_{\yb,\Omega}$ by the obvious
projection.
We define $\th\Pcb_{\nu,\Omega}$ as the set of isomorphism classes of simple
perverse sheaves over $\th\! E_{\nu,\Omega}$ which
appear as a direct factor of
$\th\!\Lc_{\yb,\Omega}[d]$ for some $\yb\in \th\!Y^\nu$ and $d\in\ZZ$. 
Next, we define $\th\!\Qcb_{\nu,\Omega}$ as the full subcategory
of $\Dcb_{\th\! G_\Vb}(\th\! E_{\nu,\Omega})$ consisting of the objects 
which are isomorphic to finite direct sums of $\Lc[d]$ with
$\Lc\in\th\Pcb_{\nu,\Omega}$ and $d\in\ZZ$. 
When there is no risk of confusion we abbreviate
$$\th\Pcb=\th\Pcb_{\nu,\Omega},\quad\th\!\Qcb=\th\!\Qcb_{\nu,\Omega},
\quad\th\!\Lc_{\yb}=\th\!\Lc_{\yb,\Omega}.$$

\subhead 10.2.~Example\endsubhead
Let $\Gamma$, $\theta$, and $\l$ be as in Sections 6.2, 6.4.
Set $\bar\Omega=\emptyset$, and $\nu=i+\theta(i)$ for some $i\in I$.
We have $\th\!E_{\nu,\Omega}=L_i\times L_{\theta(i)}$
with $L_j=\Hom(\Lbb_j,\Vb_j)$,
$\th\! I^\nu=\{\ib,\theta(\ib)\}$ with $\ib=i\theta(i)$, and
$$\gathered\th\!\widetilde F_{\ib,\Omega}=
\{(\Vb\supset\Vb_{\theta(i)}\supset 0)\}\times L_{\theta(i)},\quad
\th\!\widetilde F_{\theta(\ib),\Omega}=\{
(\Vb\supset\Vb_i\supset 0)\}\times L_i.\endgathered$$
Therefore the following holds
\vskip1mm
\itemitem{$\bullet$}
if $\l_i+\l_{\theta(i)}\neq 0$ then 
$\th\Pcb_{\nu,\Omega}=
\{\kb_{L_i}[\l_i],\, \kb_{L_{\theta(i)}}[\l_{\theta(i)}]\}$,
$\th\!\Lc^\delta_{\theta(\ib)}=\kb_{L_i}[\l_i]$, and
$\th\!\Lc^\delta_{\ib}=\kb_{L_{\theta(i)}}[\l_{\theta(i)}],$
\vskip1mm
\itemitem{$\bullet$} if $\l_i+\l_{\theta(i)}=0$ then 
$\th\Pcb_{\nu,\Omega}=\{\kb_{\{0\}}\}$ and 
$\th\!\Lc^\delta_{\theta(\ib)}=\th\!\Lc^\delta_{\ib}=\kb_{\{0\}}.$

\vskip3mm

\subhead 10.3.~Multiplication of complexes\endsubhead
Set $\nu''=\nu+\nu'+\theta(\nu')$ with $\nu'\in\NN I$.
Fix $\Vb'\in\Vcb_{\nu'}$ and $\Vb''\in\th\Vcb_{\nu''}$.
Let $T$ be the set of triples $(V,\g,\g')$ where
\vskip1mm
\itemitem{$\bullet$}
$V$ is an $I$-graded subspace of $\Vb''$ such that
$\Vb''/V\in\Vcb_{\nu'}$ and
$V^\perp\subset V,$
\vskip1mm
\itemitem{$\bullet$}
$\g:\Vb\to V/V^\perp$ is an isomorphism 
in $\th\Vcb_\nu$, 
\vskip1mm
\itemitem{$\bullet$}
$\g':\Vb'\to\Vb''/V$ is an isomorphism 
in $\Vcb_{\nu'}$. 
\vskip1mm
\noindent
We consider the following diagram
$$\xymatrix{
\th\! E_{\nu,\Omega}\times E_{\nu'}&
\th\! E_{1,\Omega}\ar[l]_-{p_1}\ar[r]^{p_2}&
\th\! E_{2,\Omega}\ar[r]^-{p_3}&\th\! E_{\nu'',\Omega}.}$$
Here $\th\! E_{2,\Omega}$ is the variety of tuples $(x,y,z,V)$ where 
\vskip1mm
\itemitem{$\bullet$}
$V$ is an $I$-graded subspace of $\Vb''$ such that
$\Vb''/V\in\Vcb_{\nu'}$ and
$V^\perp\subset V,$
\vskip1mm
\itemitem{$\bullet$}
$(x,y,z)\in\th\! E_{\nu'',\Omega}$,
$y(\Lbb)\subset V$, $x(V)\subset V$, and $z(V^\perp)=0$,
\vskip1mm
\noindent
and $\th\! E_{1,\Omega}$ is the variety of tuples $(x,y,z,V,\g,\g')$ where 
\vskip1mm
\itemitem{$\bullet$}
$(V,\g,\g')\in T$,
\vskip1mm
\itemitem{$\bullet$}
$(x,y,z,V)\in\th\! E_{2,\Omega}$.
\vskip1mm
\noindent
Finally the maps are given by 
\vskip1mm
\itemitem{$\bullet$}
$p_1(x,y,z,V,\g,\g')=(x_\g,y_\g,z_\g,x_{\g'}),$
\vskip1mm
\itemitem{$\bullet$}
$p_2(x,y,z,V,\g,\g')=(x,y,z,V),$
\vskip1mm
\itemitem{$\bullet$}
$p_3(x,y,z,V)=(x,y,z),$
\vskip1mm
\noindent
where 
\vskip1mm
\itemitem{$\bullet$}
$x_\g=\g^{-1}\circ(x|_{V/V^\perp})\circ\g,$
\vskip1mm
\itemitem{$\bullet$}$x_{\g'}=(\g')^{-1}\circ(x|_{\Vb''/V})\circ\g',$
\vskip1mm
\itemitem{$\bullet$}$y_\g=\g^{-1}\circ y$.
\vskip1mm
\itemitem{$\bullet$}$z_\g=z\circ\g$.\vskip1mm
\noindent
The group $\th\!G_{\nu''}$ acts on $\th\!E_{1,\Omega}$, $\th\!E_{2,\Omega}$ and the maps
$p_2$, $p_3$ are $\th\!G_{\nu''}$-equivariant.
Note that $p_1$ is a smooth map with connected fibers, that $p_2$ is a 
principal bundle, and that $p_3$ is proper.
Therefore, for any complexes
$\Ec\in\Dcb_{\th\! G_\nu}(\th\! E_{\nu,\Omega})$ and 
$\Ec'\in\Dcb_{G_{\nu'}}(E_{\nu'})$
there is a unique complex $\Ec_2\in\Dcb_{\th\!G_{\nu''}}(\th\! E_{2,\Omega})$ such that
$$p_1^*(\Ec\boxtimes\Ec')=p_2^*(\Ec_2).$$
Then, we define a complex 
$\Ec''=\varphi_!(\Ec,\Ec')$ in $\Dcb_{\th\! G_{\nu''}}(\th\! E_{\nu'',\Omega})$ by
$$\Ec''=(p_3)_!(\Ec_2).$$
Now, let $\nu'=i$. Hence $E_{\nu'}=0$. 
Let $\Lc_i=\kb_{E_{\nu'}}$, the 
trivial complex over $E_{\nu'}$. 

\proclaim{10.4.~Definition}
Set $\nu'=i$.
For 
$\Ec$ in $\Dcb_{\th\!G_\nu}(\th\! E_{\nu,\Omega})$ 
we define the complex
$$\underline {f_i}(\Ec)=\varphi_!(\Ec,\Lc_i)[b_{\nu,i}],\quad
b_{\nu,i}=\nu_i+\sum_j\nu_j\,h_{i,j}+\l_{\Omega,\theta(i)}+\l_{\bar\Omega,i}.$$
\endproclaim

\proclaim{10.5.~Proposition} 
(a) $\underline{f_i}$ yields a functor 
$\th\!\Qcb\to\th\!\Qcb$.

\vskip1mm

(b)
$\underline{f_i}(\th\!\Lc_{\ib}^\delta)=
\th\!\Lc_{i\ib \theta(i)}^\delta$
for each $\ib$ in $\th\! I^\nu$. 
\endproclaim

\noindent{\sl Proof :} 
A standard computation yields
$$\varphi_!(\th\!\Lc_{\yb},\Lc_i)=\th\!\Lc_{i\yb\theta(i)},\quad
\yb\in\th Y^\nu.$$
See \cite{E, prop.~4.11}, \cite{L2, sec.~9.2.6-7}. 
This implies $(a)$.
The same computation as in Proposition 2.5 yields 
$$d_{\ib}=\ell_\nu/2+\sum_{k<l;\,k+l\neq 1}h_{i_k,i_l}/2+\sum_{1\leqslant l}
(\l_{\Omega,i_l}+\l_{\bar\Omega,\theta(i_l)}),
\leqno(10.1)$$
where $\l_{\Omega,i}=\l_i$ if $i\in\Omega$ and 0 else. Thus we have
$$d_{\l,i\ib\theta(i)}-d_{\l,\ib}=b_{\nu,i}.$$
Part $(b)$  follows from this equality.

\qed

\vskip3mm

\subhead 10.6.~Restriction of complexes\endsubhead
Set $\nu''=\nu+\nu'+\theta(\nu')$ with $\nu'\in\NN I$.
Fix $\Vb'\in\Vcb_{\nu'}$,
$\Vb''\in\th\Vcb_{\nu''}$,
and fix a triple
$(V,\g,\g')\in T$. 
Consider the diagram
$$\xymatrix{
\th\! E_{\nu,\Omega}\times E_{\nu'}&
\th\! E_{3,\Omega}\ar[l]_-{\kappa}\ar[r]^-{\iota}&
\th\! E_{\nu'',\Omega}.}$$
Here we have set
\vskip1mm
\itemitem{$\bullet$}
$\th\! E_{3,\Omega}=\{(x,y,z)\in\th\! E_{\nu'',\Omega}\,;
\,x(V)\subset V,\,y(\Lbb)\subset V,\,z(V^\perp)=0\},$
\vskip1mm
\itemitem{$\bullet$}
$\kappa(x,y,z)=(x_\g,y_\g,z_\g,x_{\g'}),$
\vskip1mm
\itemitem{$\bullet$}
$\iota(x,y,z)=(x,y,z).$
\vskip1mm
\noindent
For any 
$\Ec''$ in $\Dcb_{\th\! G_{\nu''}}(\th\! E_{\nu'',\Omega})$
we define a complex in 
$\Dcb_{\th\! G_\nu\times G_{\nu'}}(\th\! E_{\nu,\Omega}\times E_{\nu'})$ by
$$\varphi^*(\Ec'')=\kappa_!\iota^*(\Ec'').$$

\proclaim{10.7.~Definition}
Set $\nu'=i$. 
For any 
$\Ec''$ in $\Dcb_{\th\! G_{\nu''}}(\th\! E_{\nu'',\Omega})$
we define 
$$\underline{e_i}(\Ec'')=\varphi^*(\Ec'')[a_{\nu,i}],\quad
a_{\nu,i}=-\nu_i+\sum_j\nu_j\,h_{i,j}+
\l_{\Omega,\theta(i)}+\l_{\bar\Omega,i}.$$
\endproclaim

\proclaim{10.8.~Proposition}
(a) $\underline{e_i}$ yields a functor from $\th\!\Qcb\to\th\!\Qcb$.

\vskip1mm

(b)
$\underline{e_i}(\th\!\Lc_{\yb})=
\bigoplus_k\th\!\Lc_{\yb_k}[-2d_k]$ for some integer $d_k$.
The sum runs over all $k$ such that $i_k=i$, and
$$\yb_k=(\ib,\ab^{(k)}),\quad
\ab^{(k)}=(a^{(k)}_l),\quad a^{(k)}_l=a_l-\delta_{l,k}-\delta_{l,1-k}.$$

\vskip1mm

(c)
$\underline{e_i}(\th\!\Lc^\delta_{\ib})=
\bigoplus_{\ib'}\th\!\Lc^\delta_{\ib'}[\deg(\ib',\theta(i);\ib)]$,
where $\ib'$ runs over all sequences such that 
$\ib$ lies in $Sh(\ib',\theta(i))$. 
\endproclaim

\noindent{\sl Proof :}
Part $(a)$ follows from $(b)$.
Parts $(b)$ and $(c)$ are analogues of
\cite{E, prop.~4.11$(ii)$} (which itself is an analogue of
\cite{L2, sec.~9.2.6}), where the case 
$\l=0$ is considered.
Our proof is similar.
Assume that $\nu'=i$. Hence we have $E_{\nu'}=\{0\}$.
We define
$$\gathered
\th\! \widetilde F_{3,\Omega}=\{(x,y,z,\phi)\in\th\!\widetilde F_{\yb,\Omega}\,;\,
(x,y,z)\in\th\! E_{3,\Omega}\},\vspace{2mm}
\th\!F^{(k)}_3=
\{\phi=(\Vb^l)\in\th\!F_{\yb}\,;
\,\Vb^k\subset V,\,\Vb^{k-1}\not\subset V\},
\vspace{2mm}
\th\!\widetilde F^{(k)}_{3,\Omega}=
\{(x,y,z,\phi)\in\th\!\widetilde F_{3,\Omega};\,
\phi\in\th\!F^{(k)}_3\}.
\endgathered
$$
Note that
$\th\!\widetilde F_{3,\Omega}=\bigcup_k\th\!\widetilde F^{(k)}_{3,\Omega}$
is a partition into locally closed subsets.
Let $\yb_k$ be as above.
Consider the map
$$\gathered
f_k:\th\!\widetilde F^{(k)}_{3,\Omega}\to
\th\! \widetilde F_{\yb_k,\Omega},\quad
(x,y,z,\phi)\mapsto(x_\g,y_\g,z_\g,\phi_\g),
\endgathered$$
where $\phi_\g$ is the flag whose $l$-th term is equal to
$$\g^{-1}\bigl((V\cap\Vb^l+V^\perp)/V^\perp\bigr).$$
We get the following diagram, whose right square is Cartesian
$$\xymatrix{
\th\! \widetilde F_{\yb_k,\Omega}\ar[d]&
\th\! \widetilde F_{3,\Omega}^{(k)}\ar[r]\ar[l]_-{f_k}&
\th\! \widetilde F_{3,\Omega}\ar[d]\ar[r]&\th\!\widetilde F_{\yb,\Omega}\ar[d]\cr
\th\! E_{\nu,\Omega}&&
\th\! E_{3,\Omega}\ar[ll]_-\kappa\ar[r]^-\iota&\th\!E_{\nu'',\Omega}.
}$$
It is easy to prove that $f_k$ is an affine bundle.
Let $d_k=d_{f_k}$ be its relative dimension.
A standard argument using the diagram above yields
$$\varphi^*(\th\!\Lc_{\yb})=\kappa_!\iota^*(\th\!\Lc_{\yb})=
\bigoplus_{i_k=i}\th\!\Lc_{\yb_k}[-2d_k].$$
This proves $(b)$. Now, we concentrate on $(c)$.
Assume that $\yb=\ib$ lies in $\th\! I^{\nu''}$.
Therefore we have
$$\ib=(i_{-m},i_{1-m},\dots,i_{m+1}),\quad k=-m,1-m,\dots,m+1.$$
First, we compute explicitly the integer $d_k$.
The map $\th\! F^{(k)}_3\to
\th\!  F_{\yb_k}$, $\phi\mapsto\phi_\g$
is an affine bundle of relative dimension 
$$\sharp\{l\,;\,-m\leqslant l<k,\,i_l=i\}.$$
Further, for each tuple $(x,y,z,\phi)$ in $\th\! \widetilde F_{\ib_k,\Omega}$
and for each $\phi'\in\th\!F_3^{(k)}$such that $\phi'_\g=\phi$, the space 
of tuples  $(x',y',z')$ in $\th\! E_{3,\Omega}$ such that 
$(x',y',z',\phi')$ lies in $\th\!\widetilde F_{3,\Omega}$ and
$\kappa(x',y',z')=(x,y,z)$
has the dimension
$$\sum_{k<l\leqslant m+1}h_{i,i_l}+
\delta_{k\leqslant 0}\,(\l_{\Omega,\theta(i)}+\l_{\bar\Omega,i}-h_{i,\theta(i)}).$$
See \cite{E, prop.~4.11$(ii)$} for details. 
Therefore we have
$$d_k=
\sum_{k<l\leqslant m+1}h_{i,i_l}
+\delta_{k\leqslant 0}\,(\l_{\Omega,\theta(i)}+\l_{\bar\Omega,i}-h_{i,\theta(i)})+
\sharp\{l\,;\,-m\leqslant l<k,\,i_l=i\}.
$$
Next, (10.1) yields
$$\gathered
d_{\l,\ib}-d_{\l,\ib_k}
=\nu_i
+\sum_{-m\leqslant l<k}h_{i_l,i}+\sum_{k<l\leqslant m+1}h_{i,i_l}
+\delta_{k\geqslant 1}(\l_{\Omega,i}+\l_{\bar\Omega,\theta(i)}-h_{\theta(i),i})+\vspace{2mm}
+\delta_{k\leqslant 0}(\l_{\bar\Omega,i}+\l_{\Omega,\theta(i)}-h_{i,\theta(i)})
.\endgathered$$
Finally we have
$$a_{\nu,i}=-\nu_i+\sum_{-m\leqslant l\leqslant m+1}
\,h_{i,i_l}-h_{i,\theta(i)}+\l_{\Omega,\theta(i)}+\l_{\bar\Omega,i}.$$
Therefore we get
$$
a_{\nu,i}+d_{\l,\ib}-d_{\l,\ib_k}-2d_k
=-\sum_{-m\leqslant l<k}i\cdot i_l
+\delta_{k\geqslant 1}(i\cdot\theta(i)
+\l_i+\l_{\theta(i)}).
$$
On the other hand $\deg(\ib_k,\theta(i)\,;\ib)$ is the degree
of the homogeneous element
$\sigma_{\dot w}1_\ib$, where $\dot w$ is a reduced decomposition
of an element $w$ of $W_{m+1}$ such that
$w(\ib)=i\,\ib_k\theta(i)$.
If $k\leqslant 0$ then we can choose
$\dot w=s_ms_{m-1}\dots s_{1-k}$ and we get
$$\deg(\ib_k,\theta(i)\,;\ib)=-\sum_{-m\leqslant l<k}i\cdot i_l.$$
If $k\geqslant 1$ then we can choose
$\dot w=s_ms_{m-1}\dots s_1\eps_1s_1\dots s_{k-1}$ and we get
$$\aligned
\deg(\ib_k,\theta(i)\,;\ib)
&=-\sum_{-m\leqslant l\leqslant 0}i\cdot i_l+i\cdot\theta(i)+\l_i+\l_{\theta(i)}
-\sum_{1\leqslant l<k}i\cdot i_l,
\vspace{1mm}
&=-\sum_{-m\leqslant l<k}i\cdot i_l+i\cdot\theta(i)+\l_i+\l_{\theta(i)}
.\endaligned$$
The proposition is proved.

\qed

\vskip3mm

\subhead 10.9.~Example\endsubhead 
Let $\Gamma$, $\theta$, $\l$, $\nu$, and $\Omega$ be as in Example 10.2. 
Let $\kb$ denote the unique element of $\th\Pcb_{0,\Omega}$. 
We have
$$\xymatrix{
\{0\}&
L_{\theta(i)}
\ar[l]_-{\kappa}\ar[r]^-{\iota}&
L_i\times L_{\theta(i)}}.$$
We obtain 
$$\gathered
a_{\nu,i}=\l_{\theta(i)},\quad
\underline{e_i}(\kb_{L_i}[\l_i])=\kb[\l_i+\l_{\theta(i)}],\quad
\underline{e_i}(\kb_{L_{\theta(i)}}[\l_{\theta(i)}])=\kb,
\vspace{2mm}
b_{\nu,i}=\l_{\theta(i)},\quad
\underline{f_i}(\kb)=\kb_{L_{\theta(i)}}[\l_{\theta(i)}].
\endgathered$$

\vskip3mm

\subhead 10.10.~A key estimate\endsubhead
First, let us introduce the following notation.
For any complex of constructible sheaves $\Lc$ and any integer $d$ we'll write
$v^d\Lc$ for the shifted complex $\Lc[d]$.

\proclaim{10.11.~Proposition} For each $i\in I$ there are maps
$$\underline\eps_i:\th\!\Pcb\cup\{0\}\to\ZZ_{\geqslant 0},\quad
\tilde{\underline E}_i:\th\!\Pcb\to\th\!\Pcb\cup\{0\},\quad
\tilde{\underline F}_i:\th\!\Pcb\to\th\!\Pcb$$
such that $\underline\eps_i(0)=0$, and for each $\Lc$ in $\th\! \Pcb$ the following hold

\vskip2mm

(a)  we have
$\underline\eps_i(\tilde{\underline E}_i(\Lc))=\underline\eps_i(\Lc)-1$ and
$$\gathered
\underline{e_i}(\Lc)=v^{1-\underline\eps_i(\Lc)}\,\tilde{\underline E}_i(\Lc)+
\sum_{\Lc'}e_{\Lc,\Lc'}\Lc',\vspace{1mm}
\Lc'\in \th\!\Pcb,
\quad\underline\eps_i(\Lc')\geqslant \underline\eps_i(\Lc),\quad
e_{\Lc,\Lc'}\in v^{1-\underline\eps_i(\Lc')}\ZZ[v],
\endgathered$$

(b)  we have
$\underline\eps_i(\tilde{\underline F}_i(\Lc))=\underline\eps_i(\Lc)+1$
and
$$\gathered
\underline{f_i}(\Lc)=\la\underline\eps_i(\Lc)+1\ra\,\tilde{\underline F}_i(\Lc)+
\sum_{\Lc'}f_{\Lc,\Lc'}\Lc',\vspace{1mm}
\Lc'\in \th\!\Pcb,\quad
\underline\eps_i(\Lc')>\underline\eps_i(\Lc)+1,\quad
f_{\Lc,\Lc'}\in v^{2-\underline\eps_i(\Lc')}\ZZ[v],
\endgathered$$

(c) we have
$$\gathered
\underline\eps_i(0)=0,\quad
\tilde{\underline E}_i(\Lc)\neq 0\ \roman{if}\ \underline\eps_i(\Lc)>0,
\vspace{2mm}
\tilde{\underline E}_i\tilde{\underline F}_i(\Lc)=\Lc,\quad
\tilde{\underline F}_i\tilde{\underline E}_i(\Lc)=\Lc\ \roman{if}\ 
\tilde{\underline E}_i(\Lc)\neq 0, 
\endgathered
$$

\vskip1mm

(d) if $\Lc\in\th\Pcb$ is such that $\underline\eps_i(\Lc)=0$ 
for all $i$, then $\Lc\in\th\Pcb_{0,\Omega}$, 

\vskip2mm

(e) the elements of $\th\Pcb$ are selfdual.
\endproclaim

\noindent{\sl Proof :}
We'll prove the proposition for any quiver $\Gamma=(I,H)$ with an involution $\theta$ such that
$\Gamma$ has no 1-loops and 
$\theta$ has no fixed points. 
The  estimates in $(a)$, $(b)$ 
are analogue of \cite{E, thm.~5.3}, where
they are proved under the assumption $\l=0$.
Our proof is essentially the same. Fix a vertex $i$.
First, 
we can assume that 
\vskip1mm
\itemitem{$\bullet$}
$i$ is a sink of $\Gamma$, 
\vskip1mm
\itemitem{$\bullet$} $i\in\Omega$,
\vskip1mm
\itemitem{$\bullet$} $\theta(i)\in\bar\Omega$.
\vskip1mm
\noindent 
More precisely we have the following lemma.
Its proof is left to the reader. It is proved as in
\cite{E, thm.~4.19}, \cite{L2}, 
using Fourier transforms.

\proclaim{10.12.~Lemma} 
Let $(\Gamma^{(1)},\theta^{(1)})$, $(\Gamma^{(2)},\theta^{(2)})$ 
be two quivers with involutions without fixed points. 
Assume that $\Gamma^{(1)}$, $\Gamma^{(2)}$ 
have the set of vertices $I$ and that they have the same set of 
unoriented arrows. Let
$\Omega^{(1)}$, $\Omega^{(2)}$ be two orientations of $I$.
There is an equivalence of categories 
$\th\!\Qcb_{\nu,\Omega^{(1)}}\to\th\!\Qcb_{\nu,\Omega^{(2)}}$
which commutes with the functors $\underline{f_i}$, $\underline{e_i}$ and 
with the Verdier duality.
The categories $\th\!\Qcb_{\nu,\Omega^{(1)}}$ and 
$\th\!\Qcb_{\nu,\Omega^{(2)}}$
are relative to the quivers
$\Gamma^{(1)}$, $\Gamma^{(2)}$  respectively.
This equivalence yields a bijection 
$\th\Pcb_{\nu,\Omega^{(1)}}\to\th\Pcb_{\nu,\Omega^{(2)}}$.
\endproclaim

Then, for each integer
$r\geqslant 0$ let $\th\!E_{\nu,\Omega,\geqslant r}$ be the closed subset of
$\th\!E_{\nu,\Omega}$ consisting of the triples
$(x,y,z)$ such that there is a $I$-graded subspace
$W\subset\Vb$ of codimension vector $ri$ such that
$$\gathered
x(W)\subset W,\quad y(\Lbb)\subset W,\quad
z(W^\perp)=0.
\endgathered$$
Then, we set
$$\th\!E_{\nu,\Omega,r}=\th\!E_{\nu,\Omega,\geqslant r}\setminus\th\!E_{\nu,\Omega,\geqslant r+1},
\quad
\th\!E_{\nu,\Omega,\leqslant r}=\th\!E_{\nu,\Omega}\setminus\th\!E_{\nu,\Omega,\geqslant r+1}.$$
Finally we set
$\underline\eps_i(0)=0$ and for $\Lc\in\th\!\Pcb$ we define
$$\underline\eps_i(\Lc)=\max\{r\,;\,
\sup(\Lc)\subset\th\!E_{\nu,\Omega,\geqslant r}\}
.$$
Set $\nu'=i$ and consider the diagram
$$\xymatrix{
\th\! E_{\nu,\Omega}&
\th\! E_{1,\Omega}\ar[l]_-{p_1}\ar[r]^{p_2}&
\th\! E_{2,\Omega}\ar[r]^-{p_3}&\th\! E_{\nu'',\Omega}.}$$
Under restriction it yields the diagram
$$\xymatrix{
\th\! E_{\nu,\Omega,r}&
\th\! E_{1,\Omega,r}\ar[r]\ar[l]&
\th\! E_{2,\Omega,r+1}\ar[r]&\th\! E_{\nu'',\Omega,r+1},}$$
where
$$\gathered
\th\!E_{1,\Omega,r}=p_1^{-1}(\th\!E_{\nu,\Omega,r}),\quad
\th\!E_{2,\Omega,r+1}=p_3^{-1}(\th\!E_{\nu'',\Omega,r+1}).
\endgathered$$
Note that we have
$\th\!E_{1,r}=p_2^{-1}(\th\!E_{2,r+1})$
and that the map 
$\th\! E_{2,r+1}\to\th\! E_{\nu'',r+1}$
is a $\PP^{r}$-bundle.
Finally, we set $p=p_3p_2$
and we define
$\th\!E_{2,\Omega,\leqslant r}$ and 
$\th\!E_{1,\Omega,\leqslant r}$ in the obvious way.

\vskip2mm

Now, we concentrate on $(b)$.
Fix a simple perverse sheaf $\Lc\in\th\!\Pcb_{\nu,\Omega}$.
Set $\eps=\underline\eps_i(\Lc)$.
The maps $p_1$, $p_2$ are smooth with connected fibers of dimension
$d_{p_1}$, $d_{p_2}$ such that
$$b_{\nu,i}=d_{p_1}-d_{p_2}.\leqno(10.2)$$
Thus, there is a unique simple $\th\!G_{\nu''}$-equivariant perverse sheaf
$\Lc_2$ on $\th\!E_{2,\Omega}$ with 
$$p_1^*(\Lc)[b_{\nu,i}]=p_2^*(\Lc_2).$$ We have
$$\underline{f_i}(\Lc)=(p_3)_!(\Lc_2).$$
The complex $\underline{f_i}(\Lc)$ is supported on $\th\!E_{\nu'',\Omega,\eps+1}$.
Further, the restriction of $\Lc_2$
to $\th\!E_{2,\Omega,\leqslant\eps+1}$  
is supported on $\th\!E_{2,\Omega,\eps+1}$,
and it is constant along the fibers of $p_3$ by $\th\!G_{\nu''}$-equivariance.
Thus  $$\Lc_2|_{\th\!E_{2,\Omega,\leqslant\eps+1}}=p_3^*(\Lc'')[\eps]$$ 
for some simple perverse sheaf
$\Lc''$ on $\th\!E_{\nu'',\Omega,\leq\eps+1}$. 
Let $\Lc_0$ be the minimal perverse extension of $\Lc''$
to $\th\!E_{\nu'',\Omega}$. 
Since $\underline{f_i}(\Lc)$ is semi-simple, we get 
$$\gathered
\underline{f_i}(\Lc)=\la\eps+1\ra\,\Lc_0+
\sum_{\Lc'}f_{\Lc,\Lc'}\Lc',\vspace{1mm}
\Lc_0, \Lc'\in\th\Pcb_{\nu'',\Omega},\quad
\underline\eps_i(\Lc_0)=\eps+1,
\quad\underline\eps_i(\Lc')>\eps+1.\endgathered$$
Let us that $f_{\Lc,\Lc'}$ lies in 
$v^{ 2-\underline\eps_i(\Lc')}\ZZ[v]$.
Write
$$\underline{f_i}(\Lc)=\bigoplus_{\Lc'}\Lc'\otimes M_{\Lc'},$$
where $M_{\Lc'}$ is a complex of $\kb$-vector spaces.
Set $\eps'=\underline\eps_i(\Lc')$.
We have
$$\R\Hc om(\Lc',\Lc')|_{\th\!E_{\nu'',\leqslant\eps'}}\otimes M_{\Lc'}^*
\subset\R\Hc om((p_3)_!(\Lc_2),\Lc')|_{\th\!E_{\nu'',\leqslant\eps'}}.$$
On the other hand, since $p_3$ restricts to a $\PP^{\eps'-1}$-bundle
$\th\!E_{2,\Omega,\eps'}\to\th\!E_{\nu'',\Omega,\eps'}$, we have
$$\R\Hc om((p_3)_!(\Lc_2),\Lc')|_{\th\!E_{\nu'',\leqslant\eps'}}=
(p_3)_*\R\Hc om(\Lc_2,\Lc'[\eps'-1])|_{\th\!E_{2,\leqslant\eps'}}[\eps'-1].$$
Since
$\Lc'[\eps'-1]|_{\th\!E_{2,\leqslant\eps'}}$ is a perverse sheaf
the complex
$$\R\Hc om(\Lc_2,\Lc'[\eps'-1])|_{\th\!E_{2,\leqslant\eps'}}$$
is concentrated in degrees $\geqslant 0$.
Its 0-th cohomology group is zero because 
$\Lc_2$ and $\Lc'[\eps'-1]$
are simple and non isomorphic.
Thus the complex
$$\R\Hc om((p_3)_!(\Lc_2),\Lc')|_{\th\!E_{\nu'',\leqslant\eps'}}$$
is concentrated in degrees $>1-\eps'$.
This implies the estimate we want.

\vskip2mm

Next, we prove $(a)$. Fix a triple $(V,\g,\g')$ in $T$. 
Observe that the hypothesis on $\Gamma$, $\Omega$, $i$ implies that
for each $(x,y,z,W,\rho,\rho')$ in $\th\!E_{1,\Omega}$ we have
$x(W^\perp)=z(W^\perp)=0$, $x(\Vb),y(\Lbb)\subset W$, 
$z$ is completely determined by
its restriction to $W$,
and $y$ is completely determined by
its composition with the projection $\Vb\to\Vb/W^\perp$.
Hence $x$, $y$, $z$ are completely determined by $x_\rho$, $y_\rho$, $z_\rho$.
Therefore $\kappa$ is an isomorphism.
Consider the diagram
$$\xymatrix{
\th\! E_{\nu,\Omega}&
\th\! E_{3,\Omega}\ar@{=}[l]_-{\kappa}\ar[r]^-{\iota}\ar[d]^-{ s}&
\th\! E_{\nu'',\Omega}\cr
&
\th\!E_{1,\Omega},\ar[ur]_-{p}\ar[ul]-_{p_1}&}\leqno(10.3)$$
where 
$$\gathered
\kappa(x,y,z)=(x_\g,y_\g,z_\g),
\quad  s(x,y,z)=(x,y,z,V,\g,\g'),
\vspace{2mm}
p_1(x,y,z,W,\rho,\rho')=(x_\rho,y_\rho,z_\rho),\quad 
p(x,y,z,W,\rho,\rho')=(x,y,z).
\endgathered$$
The left square is Cartesian.
Fix a simple 
perverse sheaf $\Lc$ in $\th\!\Pcb_{\nu'',\Omega}$.
Set $\eps=\underline\eps_i(\Lc)$. 
We'll assume that $\eps>0$ (the case $\eps=0$ is left to the reader).
We have
$$\aligned
\underline{e_i}(\Lc)
=\kappa_! s^*p^*(\Lc)[a_{\nu,i}].
\endaligned$$
The restriction $\Lc|_{\th\!E_{\nu'',\Omega,\leqslant\eps}}$ 
is a simple 
$\th\!G_{\nu''}$-equivariant
perverse sheaf
supported on $\th\!E_{\nu'',\Omega,\eps}$.
Let $d_p$, $d_s$ be the relative dimension of the maps
$p$, $s$.
Since $p$ restricts to a smooth map 
$\th\!E_{1,\Omega,\eps-1}\to\th\!E_{\nu'',\Omega,\eps}$, the complex
$$\Lc_1=p^*(\Lc)[d_p]|_{\th\!E_{1,\Omega,\leqslant\eps-1}}$$ 
is again a simple 
$\th\!G_{\nu''}$-equivariant
perverse sheaf.
It is constant along the fibers of $p_1$ by 
$\th\!G_{\nu''}$-equivariance.
Thus $$\aligned
\Lc''
&=\kappa_!s^*p^*(\Lc)[d_p+d_s]|_{\th\!E_{\nu,\Omega,\leqslant\eps-1}}
\vspace{2mm}
&=\underline{e_i}(\Lc)
[d_p+d_s-a_{\nu,i}]|_{\th\!E_{\nu,\Omega,\leqslant\eps-1}}
\endaligned$$ 
is a simple perverse sheaf over $\th\! E_{\nu,\Omega,\leqslant\eps-1}$.
Using (10.2) we get
$$
\gathered
d_p+d_s=d_{p_2}+\eps-1-d_{p_1},\quad
d_{p_1}-d_{p_2}=b_{\nu,i},\quad
b_{\nu,i}=\nu_i,
\quad a_{\nu,i}=-\nu_i.
\endgathered$$
Therefore, we have
$$d_p+d_s-a_{\nu,i}=\eps-1.$$
Let $\Lc_0$ be the minimal perverse extension of
$\Lc''$ to $\th\! E_{\nu,\Omega,}$.
Since $\underline{e_i}(\Lc)$ is semi-simple we get
$$\gathered
\underline{e_i}(\Lc)=v^{1-\eps}\,\Lc_0+
\sum_{\Lc'}e_{\Lc,\Lc'}\Lc'\vspace{1mm}
\Lc_0,\Lc'\in\th\Pcb_{\nu,\Omega},\quad
\underline\eps_i(\Lc_0)=\eps-1,\quad
\underline\eps_i(\Lc')\geqslant\eps.\endgathered\leqno(10.4)$$
Now, one proves that $e_{\Lc,\Lc'}$ lies in 
$v^{1-\underline\eps_i(\Lc')}\ZZ[v]$ as in \cite{E, thm.~5.3}.
More precisely,
since $p_1^*\underline{e_i}(\Lc)$ and $p^*(\Lc)[a_{\nu,i}]$
are constant along the fibers of $p_1$ and 
since $$\underline{e_i}(\Lc)=\kappa_!s^*p^*(\Lc)[a_{\nu,i}],$$
we have
$$p_1^*\underline{e_i}(\Lc)=
p^*(\Lc)[-b_{\nu,i}].\leqno(10.5)$$
On the other hand, we have
$$p_1^*\underline{e_i}(\Lc)=
\bigoplus_{\Lc''}p_1^*(\Lc'')\otimes M_{\Lc''},$$
where the graded $\kb$-vector space $M_{\Lc''}$ is the multiplicity space of 
the simple perverse sheaf $\Lc''\in\th\Pcb_{\nu,\Omega}$ 
in $\underline{e_i}(\Lc)$.
Let $\Lc''_2$ be the perverse sheaf over $\th\!E_{2,\Omega}$ such that 
$$p_1^*(\Lc'')[b_{\nu,i}]=p_2^*(\Lc''_2).$$
We obtain
$$\bigoplus_{\Lc''}\Lc''_2\otimes M_{\Lc''}=
p_3^*(\Lc).$$
Now, let $\Lc'$ be as in (10.4).
Set $\eps'=\underline{\eps_i}(\Lc')$.
We have
$$\aligned
\bigoplus_{\Lc''}&\R\Hom(
\Lc''_2|_{\th\!E_{2,\Omega,\leqslant\eps'+1}},
\Lc'_2|_{\th\!E_{2,\Omega,\leqslant\eps'+1}})\otimes M_{\Lc''}^*=
\vspace{2mm}
&=\R\Hom(
p_3^*(\Lc)|_{\th\!E_{2,\Omega,\leqslant\eps'+1}},
\Lc'_2|_{\th\!E_{2,\Omega,\leqslant\eps'+1}})
\vspace{2mm}
&=\R\Hom(
\Lc|_{\th\!E_{\nu'',\Omega,\leqslant\eps'+1}},
(p_3)_!(\Lc'_2)|_{\th\!E_{\nu'',\Omega,\leqslant\eps'+1}})
\vspace{2mm}
&=\R\Hom(
\Lc|_{\th\!E_{\nu'',\Omega,\leqslant\eps'+1}},
\underline{f_i}(\Lc')|_{\th\!E_{\nu'',\Omega,\leqslant\eps'+1}})
\vspace{2mm}
&=\la\eps'+1\ra\,\R\Hom(
\Lc|_{\th\!E_{\nu'',\Omega,\leqslant\eps'+1}},
\tilde{\underline F}_i(\Lc')|_{\th\!E_{\nu'',\Omega,\leqslant\eps'+1}}),
\endaligned$$
where the last equality follows from part $(b)$. The complex 
$$\R\Hom(
\Lc|_{\th\!E_{\nu'',\Omega,\leqslant\eps'+1}},
\tilde{\underline F}_i(\Lc')|_{\th\!E_{\nu'',\Omega,\leqslant\eps'+1}}),$$
is concentrated in degrees $\geqslant 1$,
because the perverse sheaves
$\Lc'$ and $\tilde{\underline F}_i(\Lc')$ are simple and distincts.
Thus the complex
$$\R\Hom(
\Lc|_{\th\!E_{\nu'',\Omega,\leqslant\eps'+1}},
\underline{f_i}(\Lc')|_{\th\!E_{\nu'',\Omega,\leqslant\eps'+1}})$$
is concentrated in degrees $\geqslant 1-\eps'$.
Choosing $\Lc'=\Lc''$ we get that $M^*_{\Lc'}$ is also 
concentrated in degrees $\geqslant 1-\eps'$.
Therefore 
$$e'_i(\Lc)=
\bigoplus_{\Lc''}\Lc''\otimes M_{\Lc''}
=\bigoplus_{\Lc''}\bigoplus_{d\in\ZZ}v^{-d}\Lc''\otimes M_{\Lc'',d},$$
with $M_{\Lc'',d}=0$ unless $-d\geqslant 1-\eps'$. We are done.

\vskip2mm

Now, we concentrate on $(c)$. 
The second claim in $(c)$ is obvious.
Now, we prove that $\tilde{\underline E}_i\tilde{\underline F}_i(\Lc)=\Lc$
for $\Lc$ in $\th\Pcb_{\nu,\Omega}$.
Recall the diagram (10.3).
Set $\eps=\underline\eps_i(\Lc)$ and take a simple perverse sheaf
$\Lc_2$ on $\th\!E_{2,\Omega}$ such that
$$p_1^*(\Lc)[b_{\nu,i}]=p_2^*(\Lc_2),\quad
(p_3)_!(\Lc_2)=\underline{f_i}(\Lc).$$
We have
$$(p_3)_!(\Lc_2)|_{\th\!E_{\nu'',\Omega,\leqslant\eps+1}}=
\la\eps+1\ra\,\tilde{\underline F}_i(\Lc)|_{\th\!E_{\nu'',\Omega,\leqslant\eps+1}}.$$
On the other hand, since
$$\Lc_2|_{\th\!E_{2,\Omega,\leqslant\eps+1}}=
p_3^*(\tilde{\underline F}_i\Lc)[\eps]|_{\th\!E_{2,\Omega,\leqslant\eps+1}},$$
we have
$$p^*(\tilde{\underline F}_i\Lc)|_{\th\!E_{1,\Omega,\leqslant\eps+1}}=
p_2^*(\Lc_2)[-\eps]|_{\th\!E_{1,\Omega,\leqslant\eps+1}}=
p_1^*(\Lc)[b_{\nu,i}-\eps]|_{\th\!E_{1,\Omega,\leqslant\eps+1}}.$$
Therefore we have also
$$\aligned
\underline{e_i}
(\tilde{\underline F}_i\Lc)
|_{\th\!E_{\nu,\Omega,\leqslant\eps+1}}
&=
\kappa_!s^*p^*(\tilde{\underline F}_i\Lc)[a_{\nu,i}]
|_{\th\!E_{\nu,\Omega,\leqslant\eps+1}}
\vspace{2mm}
&=
\Lc[-\eps]|_{\th\!E_{\nu,\Omega,\leqslant\eps+1}}.
\endaligned$$
Therefore $\tilde{\underline E}_i\tilde{\underline F}_i(\Lc)=\Lc$.
Finally, fix $\Lc\in\th\Pcb_{\nu'',\Omega}$ such that
$\underline\eps_i(\Lc)>0$ and let us
prove that
$\tilde{\underline F}_i\tilde{\underline E}_i(\Lc)=\Lc$.
Write $\eps=\underline\eps_i(\Lc)$.
By (10.5) we have
$$p_1^*\underline{e_i}(\Lc)=p^*(\Lc)[a_{\nu,i}].$$
Hence we have also
$$p_1^*(\tilde{\underline E_i}\Lc)[-a_{\nu,i}]
|_{\th\! E_{1,\Omega,\leqslant\eps-1}}=
p^*(\Lc)[\eps-1]
|_{\th\! E_{1,\Omega,\leqslant\eps-1}}
.$$
Since $p_3^*(\Lc)[\eps-1]|_{\th\! E_{2,\Omega,\eps}}$
is a simple perverse sheaf, we have
$$\underline{f_i}(\tilde{\underline E}_i\Lc)
|_{\th\! E_{\nu'',\Omega,\leqslant \eps}}
=
(p_3)_!p_3^*(\Lc)[\eps-1]
|_{\th\! E_{\nu'',\Omega,\leqslant \eps}}
=\la\eps\ra\,\Lc|_{\th\! E_{\nu'',\Omega,\leqslant \eps}}.
$$
This implies that $\tilde{\underline F}_i\tilde{\underline E}_i(\Lc)=\Lc$.

\vskip2mm

Next, $(d)$ is obvious. 
If $\nu\neq 0$ we choose $\yb$, $d$ such that $\Lc[d]$ is a direct summand of
$\th\Lc_\yb$. We may assume that $\yb=(\ib,\ab)$ with $a_1>0$. Then
$\underline\eps_{i_1}(\Lc)>0$ by $(b)$ and Proposition 10.5$(b)$.

\vskip2mm

Finally, we prove $(e)$ by descending induction on $\nu$.
Any element in $\th\Pcb_{0,\Omega}$ is selfdual. 
Assume that $\nu>0$. By part $(d)$ there is $i$ such that
$\underline\eps_i(\Lc)>0$.
Set $\eps=\underline\eps_i(\Lc)$.
We prove that $\Lc$ is selfdual by descending induction on 
$\eps$.
By parts $(b)$, $(c)$ we have
$$\gathered
\underline{f_i}(\tilde{\underline E}_i\Lc)=\la\eps\ra\,\Lc+
\sum_{\Lc'}f_{\Lc,\Lc'}\Lc',\quad
\underline\eps_i(\Lc')>\eps.
\endgathered$$
The perverse sheaf $\tilde{\underline E}_i\Lc$ is selfdual
by the induction hypothesis on $\nu$.
It is easy to see that $\underline{f_i}$ 
commutes with the Verdier duality.
Hence the laft hand side is also selfdual. 
We have
$$\underline{f_i}(\tilde{\underline E}_i\Lc)|_
{\th\!E_{\nu,\Omega,\leqslant\eps}}=\la\eps\ra\,\Lc
|_{\th\!E_{\nu,\Omega,\leqslant\eps}}.$$
Since $\Lc$ is the minimal extension of its restriction to 
$\th\!E_{\nu,\Omega,\leqslant\eps}$, it is selfdual.

\qed

\vskip3mm

Let $K(\th\!\Qcb_{\nu,\Omega})$ be the Abelian group with one generator $[\Lc]$ 
for each isomorphism class of objects of
$\th\!\Qcb_{\nu,\Omega}$ and with relations $[\Lc]+[\Lc']=[\Lc'']$ 
whenever $\Lc''$ is isomorphic to
$\Lc\oplus\Lc'$. To unburden the notation we'll abbreviate
$$K(\th\!\Qcb)=\bigoplus_\nu K(\th\!\Qcb_{\nu,\Omega}),\quad\Lc=[\Lc].$$
Note that $K(\th\!\Qcb)$ is a free $\Ac$-module such that $v\Lc=\Lc[1]$
and $v^{-1}\Lc=\Lc[-1]$. Further
the Verdier duality yields an $\Ac$-antilinear map
$K(\th\!\Qcb)\to K(\th\!\Qcb)$.

\proclaim{10.13.~Corollary} 
The $\Ac$-module $K(\th\!\Qcb_{\nu,\Omega})$ is spanned by 
$\{\th\!\Lc_\yb^\delta;\,\yb\in\th Y^\nu\}$.
\endproclaim

\noindent{\sl Proof :}
The corollary is proved as in Lemma 8.34, using Proposition 10.11 instead of 
Propositions 8.22, 8.23.

\qed

\vskip3mm

\subhead 10.14.~Example\endsubhead 
Let $\Gamma$, $\theta$, $\l$, $\nu$ 
be as in Example 10.2, 
and set $\Omega=\{i\}$.
We have $\th\!E_{\nu,\Omega}=L_i\times L_{\theta(i)}^*$, 
$\th\!E_{\nu,\Omega,0}=\th\!E_{\nu,\Omega}\setminus\{0\}$,
and $\th\!E_{\nu,\Omega,1}=\{0\}$.
We have also
\vskip1mm
\itemitem{$\bullet$}
if $\l_i+\l_{\theta(i)}\neq 0$ then 
$\th\Pcb_{\nu,\Omega}=\{\kb_{\th\!E_{\nu,\Omega}}[\l_i+\l_{\theta(i)}],
\, \kb_{\{0\}}\}$, and
$$\gathered
\underline\eps_i(\kb_{\th\!E_{\nu,\Omega}}[\l_i+\l_{\theta(i)}])=0,\quad
\underline\eps_i(\kb_{\{0\}})=1,
\quad
\tilde{\underline F}_i(\kb)=\kb_{\{0\}},
\endgathered$$
\vskip1mm
\itemitem{$\bullet$} if $\l_i+\l_{\theta(i)}=0$ then 
$\th\Pcb_{\nu,\Omega}=\{\kb_{\{0\}}\}$, and
$$\gathered
\underline\eps_i(\kb_{\{0\}})=1,
\quad
\tilde{\underline F}_i(\kb)=\kb_{\{0\}}.
\endgathered$$

\vskip3mm

\subhead 10.15.~Comparison of the crystals\endsubhead
We choose $\Gamma$, $\theta$ and $\l$ as in Sections 6.2, 6.4,
and we set $\Omega=I$.
We define a functor 
$$\Yb:\th\!\Qcb_{\nu,\Omega}\to\th\Rb_\nu\text{-}\modb,\quad \Yb(\Lc)=
\bigoplus_{\ib\in\th\! I^\nu}\Ext^*_{\th\! G_\nu}(\th\!\Lc^\delta_\ib,\Lc).$$
The functor $\Yb$ is additive and it commutes with the shift 
(the shift of complexes in the left hand side 
and the shift of the grading in the 
right hand side).

\proclaim{10.16.~Proposition} 
(a) $\Yb$ takes $\th\!\Qcb$ to $\th\Rb\text{-}\proj$, and
$\th\!\Lc_\yb^\delta$ to $\th\Rb_\yb$. 
It maps $\th\Pcb$  bijectively to the set of
$\sharp$-selfdual indecomposable projective graded modules. 

\vskip1mm

(b) $\Yb$ yields an $\Ac$-module isomorphism $K(\th\!\Qcb)\to\th\Kb_I$ 
which maps $\th\Pcb$ bijectively onto $\th\!\Gb^\low(\l)$. 
It commutes with the duality.
We have 
$$e_i\circ\Yb=\Yb\circ\underline{e_{\theta(i)}},\quad 
f_i\circ\Yb=\Yb\circ\underline{f_{\theta(i)}}.$$
\endproclaim

\noindent{\sl Proof :}
Theorem 5.8, proved in Section 9, yields a
graded $\kb$-algebra isomorphism
$$\th\Rb_\nu=\th\Zb^\delta_\nu.$$
Recall that the right hans side is the graded $\kb$-algebra 
$$\th\Zb^\delta_\nu=\bigoplus_{\ib,\ib'\in\th\! I^\nu}\Ext^*_{\th\! G_\nu}(\th\!\Lc^\delta_\ib,\th\!\Lc^\delta_{\ib'}),$$
equipped with the Yoneda composition, see Sections 2.6, 2.8.
Therefore the first claim of $(a)$ is obvious. 
If the sequence $\ib$ of $\th\! I^\nu$ is the expansion of the pair $\yb$ in $\th Y^\nu$ then we have
$$\th\Rb_{\ib}=\la\bb\ra!\,\th\Rb_\yb,\quad
\th\!\Lc_{\ib}^\delta=\la\bb\ra!\,\th\!\Lc^\delta_\yb,$$ 
where $\bb$ is a sequence such that the multiplicity of $\yb$
is $\theta(\bb)\bb$.
See Remark 2.7 and (8.5).
Therefore to prove the second claim of $(a)$ 
it is enough to observe that we have
$\Yb(\th\!\Lc^\delta_{\ib})=\th\Rb^\delta_{\ib}$. 
Next, the same proof as in \cite{VV, sec.~4.7} implies
that $\Yb$ takes any element in
$\th\Pcb$ to an indecomposable projective graded module. 
Indeed, since $\Yb(\th\!\Lc^\delta_{\yb})=\th\Rb^\delta_{\yb}$ and both sides are 
selfdual, the functor $\Yb$ takes
the elements of $\th\Pcb$ to $\sharp$-selfdual indecomposable projective graded modules,
see Sections 2.6 and 8.10.
Further, any  $\sharp$-selfdual indecomposable projective graded module is a direct summand of
$\Yb(\th\!\Lc^\delta_{\yb})=\th\Rb^\delta_{\yb}$ for some $\yb$, hence is the image by $\Yb$
of an element of $\th\Pcb$.
Part $(a)$ is proved.
Next, the first claim of $(b)$ follows from Definition 8.3, Proposition 8.4$(c)$ and 
Corollary 10.13.
Finally the last claim of $(b)$ follows from
Propositions 10.5$(b)$, 10.8$(c)$ and Proposition 8.14.

\qed

\vskip3mm

Recall the set $\th\!\Gb^\low(\l)$ introduced in Definition 8.3.
For $b\in\th\!B(\l)$ let $\th\!\Lc(b)$ 
denote the unique element in $\th\Pcb$ such that
$$\Yb(\th\!\Lc(b))=\th\!G^\low(b).\leqno(10.6)$$ Hence we have 
$\th\Pcb=\{\th\!\Lc(b);\,b\in \th\!B(\l)\}$.
We'll set also $\th\!\Lc(0)=0$.
Combining Propositions 8.23 and 10.11 we can now compare the crystal
$(\th\!B(\l),\tilde{E}_i,\tilde{F}_i,\eps_i)$ from Proposition 8.22
with the crystal
$(\th\Pcb,\tilde{\underline E}_i,\tilde{\underline F}_i,\underline\eps_i)$
from Proposition 10.11.

\proclaim{10.17.~Proposition} 
For $i\in I$ and $b\in \th\!B(\l)$ we have 
$$\tilde{\underline E}_i(\th\!\Lc(b))=\th\Lc(\tilde E_{\theta(i)}b), 
\quad
\tilde{\underline F}_i(\th\!\Lc(b))=\th\Lc(\tilde F_{\theta(i)}b),
\quad
\underline\eps_i(\th\!\Lc(b))=\eps_{\theta(i)}(b).$$\endproclaim

\noindent{\sl Proof :}
We can regard
$\underline\eps_i$,
$\tilde{\underline E}_i$, and
$\tilde{\underline F}_i$ as maps
$$\underline\eps_i:\th\!B(\l)\cup\{0\}\to\ZZ_{\geqslant 0},\quad
\tilde{\underline E}_i:\th\!B(\l)\to\th\!B(\l)\cup\{0\},\quad
\tilde{\underline F}_i:\th\!B(\l)\to\th\!B(\l).$$
Propositions 10.11$(b)$, 10.16$(b)$ yield
$$f_{\theta(i)}\th\!G^\low(b)=\la\underline\eps_i(b)+1\ra\,
\th\!G^\low(\tilde{\underline F}_ib)+
\sum_{b'}f_{b,b'}\th\!G^\low(b'),\quad
\underline\eps_i(b')>\underline\eps_i(b)+1.$$ 
Taking the transpose, Definition 8.8 and Proposition 8.4$(a)$ yield 
$$e_{\theta(i)}\th\!G^\up(b)=\la\underline\eps_i(b)\ra\,
\th\!G^\up(\tilde{\underline E}_ib)+
\sum_{b'}f_{b',b}\th\!G^\up(b'),\quad
\underline\eps_i(b')<\underline\eps_i(b)-1.$$
Now, recall that
$$\eps_{\theta(i)}(b)=
\max\{n\geqslant 0;e_{\theta(i)}^n\th\!G^\up(b)\neq 0\},
\quad
\underline\eps_i(\tilde{\underline E}_ib)=\underline\eps_i(b)-1.$$
Thus, using Proposition 8.17 and (8.8) we get
$\underline\eps_i=\eps_{\theta(i)}$. 
Then, comparing the formulas above with Proposition 8.23 we get
$\tilde{\underline F}_i=\tilde F_{\theta(i)}$. 
Finally, Proposition 8.22$(c)$ and 10.11$(c)$ yield
$\tilde{\underline E}_i=\tilde E_{\theta(i)}$.

\qed

\subhead 10.18.~The global bases of $\th\Vb(\l)$\endsubhead
Since the operators $e_i$, $f_i$ on $\th\Vb(\l)$ satisfy the relations
$e_if_i=v^{-2}f_ie_i+1$, we can define the modified root operators
$\tilde{\eb}_i$, $\tilde{\fb}_i$ on the $\th\Bcb$-module 
$\th\Vb(\l)$ as follows.
For $u\in\th\Vb(\l)$ we write 
$$\gathered
u=\sum_{n\geqslant 0}f_i^{(n)}u_n\ \roman{with}\ e_iu_n=0,
\vspace{1mm}
\tilde{\eb}_i(u)=\sum_{n\geqslant 1}f_i^{(n-1)}u_n,\quad
\tilde{\fb}_i(u)=\sum_{n\geqslant 0}f_i^{(n+1)}u_n.
\endgathered$$
Let $\Rc\subset\Kc$ be the set of functions which are regular at $v=0$.
Let $\th\Lb(\l)$ be the $\Rc$-submodule of $\th\Vb(\l)$ spanned by the elements
$\tilde{\fb}_{i_1}\dots\tilde{\fb}_{i_l}(\phi_\l)$ with
$l\geqslant 0$, $i_1,\dots,i_l\in I$.
We can now apply the results in  \cite{EK3}. 
Together with Propositions 10.16 and 10.17 this yields the following,
which is the main result of the paper.

\proclaim{10.19.~Theorem} (a) We have
$$\gathered
\th\Lb(\l)=\bigoplus_{b\in \th\! B(\l)}\Rc\, \th\!G^\low(b),\quad
\tilde{\eb}_i(\th\Lb(\l))\subset\th\Lb(\l),
\quad\tilde{\fb}_i(\th\Lb(\l))\subset\th\Lb(\l),\vspace{2mm}
\tilde{\eb}_i\th\!G^\low(b)=\th\!G^\low(\tilde E_ib)\ \mod\ v\,\th\Lb(\l),\quad
\tilde{\fb}_i\th\!G^\low(b)=\th\!G^\low(\tilde F_ib)\ \mod\ v\,\th\Lb(\l).
\endgathered$$

(b) The assignment $b\mapsto \th\!G^\low(b)\ \mod\ v\,\th\Lb(\l)$ 
yields a bijection from $\th\! B(\l)$
to the subset of $\th\Lb(\l)/v\th\Lb(\l)$ consisting of the 
$\tilde{\fb}_{i_1}\dots\tilde{\fb}_{i_l}(\phi_\l)$'s.
Further $\th\!G^\low(b)$ is the unique element $x$ of $\th\Vb(\l)$
satisfying the following conditions
$$x^\sharp=x,\quad x=\th\!G^\low(b)\ \mod\ v\,\th\Lb(\l).$$

(c) For $b,b'\in\th\!B(\l)$ let $E_{i,b,b'}, F_{i,b,b'}\in\Ac$ 
be the coefficients of
$\th\!G^\low(b')$ in $e_{\theta(i)}\th\!G^\low(b)$, 
$f_i\th\!G^\low(b)$ respectively. Then we have
$$\gathered
E_{i,b,b'}|_{v=1}=
[F_i\Psi\forb(\th\!G^\up(b')):\Psi\forb(\th\!G^\up(b))],\vspace{2mm}
F_{i,b,b'}|_{v=1}=[E_i\Psi\forb(\th\!G^\up(b')):\Psi\forb(\th\!G^\up(b))].
\endgathered
$$
\endproclaim

\noindent{\sl Proof :}
Proposition 10.17 implies that $\Yb$ intertwines
the crystal operators $\tilde E_{\theta(i)}$, $\tilde F_{\theta(i)}$ 
on $\th B(\l)$ and the crystal operators $\tilde{\underline E}_i$,
$\tilde{\underline F}_i$ on $\th\pmb\Pc$.
Proposition 10.16 implies that $\Yb$ intertwines
the operators $e_{\theta(i)}$, $f_{\theta(i)}$ and $\underline e_i$,
$\underline f_i$.
Therefore, formula (10.6) and Proposition 10.11 yield estimate for the
action of $e_i$, $f_i$ on $\th\Gb^\low(\l)$ which were not available
in Proposition 8.23.
Using these estimates, part $(a)$ follows from
\cite{EK3, thm.~4.1, cor.~4.4}, \cite{E, Section 2.3}.
The first claim in $(b)$ follows from $(a)$ and Proposition 8.22.
The second one is obvious.
Part $(c)$ follows from Proposition 8.17.
More precisely, by Cartan duality
we can regard the elements $E_{i,b,b'}$, $F_{i,b,b'}$ of $\Ac$ 
as the coefficients of
$\th\!G^\up(b)$ in the expansion of $f_{\theta(i)}\th\!G^\up(b')$, 
$e_i\th\!G^\up(b')$ with respect to the basis $\th\Gb^\up(\l)$. 
Therefore, by Proposition 8.17 
we can regard the integers $E_{i,b,b'}|_{v=1}$, $F_{i,b,b'}|_{v=1}$ 
as the coefficients of
$\Psi\forb(\th\!G^\up(b))$ in $F_i\Psi\forb(\th\!G^\up(b'))$, 
$E_i\Psi\forb(\th\!G^\up(b'))$ respectively.

\qed

\vskip2cm

\head A.~Appendix\endhead

The statements above generalizes to affine Hecke algebras
of type C. 

\subhead A.1.~Affine Hecke algebras of type C\endsubhead 
Fix $p$, $q_0$, $q_1$ in $\kb^\times$. 
For any integer $m\geqslant 0$ we define
the affine Hecke algebra $\Hb_m$ of type $\C_m$ 
to be the affine Hecke algebra of $Sp(2m)$.
It admits the following presentation.
If $m>0$ then $\Hb_m$ is the $\kb$-algebra generated by
$$T_k,\quad X_l^{\pm 1},\quad 
k=0,1,\dots,m-1,\quad l=1,2,\dots,m$$ satisfying the
following defining relations :

\vskip2mm
\itemitem{$(a)$} $X_lX_{l'}=X_{l'}X_l$,

\vskip2mm

\itemitem{$(b)$} $(T_0T_1)^2=(T_1T_0)^2$,
$T_kT_{k-1}T_k=T_{k-1}T_kT_{k-1}$ if $k\neq 0,1$, and
$T_kT_{k'}=T_{k'}T_k$ if $|k-k'|\neq 1$,

\vskip2mm

\itemitem{$(c)$} $T_0X_1^{-1}-X_1T_0=(q_1^{-1}-q_0)X_1+(q_0q_1^{-1}-1)$,
$T_kX_kT_k=X_{k+1}$ if $k\neq 0$, and $T_kX_l=X_lT_k$ if $l\neq
k,k+1$,

\vskip2mm

\itemitem{$(d)$}
$(T_k-p)(T_k+p^{-1})=0$ if $k\neq 0$, and $(T_0-q_0)(T_0+q_1^{-1})=0$.

\vskip3mm

\noindent If $m=0$ then $\Hb_0=\kb$, the trivial $\kb$-algebra.

\vskip3mm

\subhead A.2.~Remark\endsubhead
The affine Hecke algebra of type $\B_m$ is equal to  $\Hb_m/(q_0-q,q_1-q)$.

\vskip3mm

\subhead A.3.~Intertwiners and blocks of $\Hb_m$\endsubhead 
We define
$$\gathered
\Ab'=\Ab[\Sigma^{-1}],\quad
\Hb'_m=\Ab'\otimes_\Ab\Hb_m,
\endgathered$$ 
where $\Sigma$ is the multiplicative system generated by
$$X_{l'}^{\pm 1}-X_l,\quad X_{l'}^{\pm 1}-p^2X_l,\quad 
1-X_l^{2},\quad 1+q_0X_l^{\pm 1},\quad 1-q_1X_l^{\pm 1},\quad l\neq l'.$$
For $k=0,\dots,m-1$ the intertwiner
$\varphi_k$ in $\Hb'_m$ is given by the
following formulas
$$
\matrix
&\displaystyle{\varphi_k-1={X_k-X_{k+1}\over pX_k-p^{-1}X_{k+1}}\,(T_k-p)}
\hfill
&\roman{if}\ k\neq 0, \vspace{2mm}
&\displaystyle{\varphi_0-1=q_1{X_{1}^{2}-1\over
(X_1+q_0)(X_1-q_1)}\,(T_0-q_0).}\hfill&
\endmatrix\leqno(A.1)$$
There is an isomorphism of $\Ab'$-algebras
$$\Ab'\rtimes W_m\to\Hb'_m,\quad
s_k\mapsto\varphi_k,\quad \eps_1\mapsto\varphi_0,\quad k\neq 0.$$
The semi-direct product
group $\ZZ\rtimes\ZZ_2$ acts on $\kb^\times$ as in Section 6.2.
Given a $\ZZ\rtimes\ZZ_2$-invariant subset $I$ of $\kb^\times$ we denote by
$\Hb_m\text{-}\Modb_I$ the category of all $\Hb_m$-modules
such that the action of $X_1,X_2,\dots,X_m$  is locally finite and
all the eigenvalues belong to $I$. 
We associate to the set $I$ a quiver $\Gamma$  with an involution $\theta$ 
as in loc.~cit. Finally, we assume that  
$$1,-1\notin I,\quad p,q_0,q_1\neq 1,-1. $$
Next, we define the element $\l$ of $\NN I$ as
$$\l=\cases
\sum_{i}i,\quad i\in I\cap\{-q_0,q_1\},&\roman{if}\ -q_0\neq q_1,\vspace{2mm}
2\sum_{i}i,\quad i\in I\cap\{q_1\},&\roman{if}\ -q_0=q_1.
\endcases
\leqno(A.2)$$ 
Finally we define 
$\th\Rb_{m}$
and  $\th\Rb_{m}\text{-}\Modb_0$
as in Section 6.4.
Note that if $q_0=q_1=q$ then $\l$ is the same as in (6.2).
Fix any formal series $f(\varkappa)$ in $\kb[[\varkappa]]$ such that
$f(\varkappa)=1+\varkappa$ modulo $(\varkappa^2)$.

\proclaim{A.4.~Theorem} There is an equivalence of categories 
$$\th\Rb_{m}\text{-}\Modb_0\to\Hb_m\text{-}\Modb_I,\quad M\mapsto M$$
which is given by
\vskip1mm

\itemitem{$(a)$}
$X_l$ acts on $1_\ib M$ by $i_l^{-1}f(\varkappa_l)$ for $l=1,2,\dots,m$, 
\vskip1mm
\itemitem{$(b)$}$T_k$ acts on $1_\ib M$ 
as follows for $k=1,2,\dots, m-1,$
$$\matrix
&\displaystyle
{{(pf(\varkappa_k)-p^{-1}f(\varkappa_{k+1}))(\varkappa_k-\varkappa_{k+1})
\over f(\varkappa_k)-f(\varkappa_{k+1})}\sigma_k+p}
\hfill&
\text{if}\ i_{k+1}=i_{k},\hfill\vspace{2mm}
&\displaystyle{
{f(\varkappa_k)-f(\varkappa_{k+1})\over
(p^{-1}f(\varkappa_k)-pf(\varkappa_{k+1}))
({\varkappa_k}-{\varkappa_{k+1})}}\sigma_k+{(p^{-2}-1)f(\varkappa_{k+1})
\over pf(\varkappa_k)-p^{-1}f(\varkappa_{k+1})}}
&\text{if}\  i_{k+1}=p^2i_{k},\hfill\vspace{2mm}
&\displaystyle{
{pi_{k}f(\varkappa_k)-p^{-1}i_{k+1}f(\varkappa_{k+1})\over
i_{k}f(\varkappa_k)-i_{k+1}f(\varkappa_{k+1})}\sigma_k+{(p^{-1}-p)i_{k}f(\varkappa_{k+1})\over
i_{k+1}f(\varkappa_k)-i_{k}f(\varkappa_{k+1})}}
\hfill&\text{if}\  i_{k+1}\neq i_{k},p^2i_{k} ,\hfill
\endmatrix$$
\vskip1mm
\itemitem{$(c)$}$T_0$ acts on $1_\ib M$ 
as follows 
$$\matrix
\displaystyle{
{(f(\varkappa_1)-1)^2\over
(q_1f(\varkappa_1)^2-q_1^{-1})\varkappa_1^2}\pi_1+
{(q_0^{-1}q_1^{-1}-1)f(\varkappa_1)^2+2f(\varkappa_1)\over
q_1^{-1}f(\varkappa_1)^2-q_1}}
\hfill&\roman{if}\ i_1=-q_0=q_1,
\vspace{2mm}
\displaystyle{{(q_1f(\varkappa_1)+q_0)(f(\varkappa_1)-1)\over
(1-q_1^2f(\varkappa_1)^2)\varkappa_1}\pi_1+
{(q_0-q_1^{-1})f(\varkappa_1)^2+(q_1-q_0)f(\varkappa_1)\over
f(\varkappa_1)^2-q_1^2}}
\hfill&\roman{if}\ i_1=q_1\neq -q_0,
\vspace{2mm}
\displaystyle{{(q_0f(\varkappa_1)+q_1)(1-f(\varkappa_1))\over
q_1(q_0f(\varkappa_1)^2-q_0^{-1})\varkappa_1}\pi_1+
{(q_1-q_0^{-1})f(\varkappa_1)^2-(q_1-q_0)f(\varkappa_1)\over
q_1(q_0^{-1}f(\varkappa_1)^2-q_0)}}
\hfill&\roman{if}\ i_1=-q_0\neq q_1,
\vspace{2mm}
\displaystyle{{(i_1f(\varkappa_1)+q_0)(i_1f(\varkappa_1)-q_1)\over
q_1(i_1^2f(\varkappa_1)^2-1)}\pi_1+{(q_1q_0-1)f(\varkappa_1)^2+(q_1-q_0)i_1f(\varkappa_1)\over
q_1(f(\varkappa_1)^2-i_1^2)}}\hfill&\roman{if}\   i_1\neq-q_0,q_1.
\hfill
\endmatrix$$
\endproclaim

\noindent{\sl Proof :}
Formula  $(A.2)$ yields $$
\cases
\lambda_{i_1}=2&\roman{if}\ i_1=-q_0=q_1,\vspace{2mm} 
\lambda_{i_1}=0&\roman{if}\  i_1\neq-q_0, q_1,
\vspace{2mm} 
\lambda_{i_1}=1&\roman{else}.
\endcases$$
The proof is the same as the proof  Theorem 6.5, using (5.4) and (A.1).

\qed

\vskip3mm

Now, all the statements in Section 8 generalizes. The proof is straightforward and is left to the reader.
In particular, if $\th\Kb_I$ denotes the Grothendieck group of the category 
$\th\Rb$-$\proj$, then
we have a canonical isomorphism $$\th\Vb(\l)=\Kc\otimes_\Ac\th\Kb_I,$$
where $\th\Vb(\l)$ the same $\th\Bcb$-module as in Theorem 8.31, 
with $\l$ given by $(A.2)$ instead of (6.2).
Theorem 10.19 generalizes as well.

\vskip3mm

\head Index of notation\endhead

\vskip2mm
\itemitem{0.1 :} 
$\mb$, $\Sen_m$, $\Sen_\mb$, $\ell_m$, $\ell_\mb$, $\la m\ra$,

\vskip2mm
\itemitem{0.2 :} 
$K(\Rb)$, $G(\Rb)$, $\hom_\Rb$, $\Hom_\Rb$, $\Ac$, $\gdim$, $\kb$,

\vskip2mm
\itemitem{0.3 :} 
$\Sb_G$,

\vskip2mm
\itemitem{1.1 :} 
$\Gamma$, $H_{i,j}$, $h_{i,j}$, $i\to j$, $i\not\to j$, 
$i\cdot j$, $E_\Vb$, $G_\Vb$, $\Vcb$, $Y^\nu$, $I^\nu$, $Y^m$,

\vskip2mm
\itemitem{1.2 :}
$F_{\Vb,\yb}$, $\tilde F_{\Vb,\yb}$, $\pi_\yb$, $d_\yb$,

\vskip2mm
\itemitem{1.3 :}
$\Sb_\Vb$, $\Lc_\yb$, $\Lc^\delta_\yb$, $\Zb_\Vb$, $\Zb^\delta_\Vb$, $\Fb_\Vb$, 
$\Rb(\Gamma)_\nu$, 

\vskip2mm
\itemitem{2.1 :} 
$\theta$, $\varpi$, $\th\NN I$,
$\th\!E_\Vb$, $\th\!E_{\Lbb,\Vb}$, $L_{\Lbb,\Vb}$,
$\th\!G_\Vb$, $\th\!\Vcb$, 

\vskip2mm
\itemitem{2.2 :} 
$F(\Wb)$, $F(\Wb,\varpi)$, 

\vskip2mm
\itemitem{2.3 :} 
$\th\!I^\nu$, $\th\!Y^\nu$, $\th\!Y^m$,

\vskip2mm
\itemitem{2.4 :} 
$\th\!F_{\Vb,\yb}$,
$\th\!\tilde F_{\Lbb,\Vb,\yb}$,
$\th\!\pi_{\Lbb,\yb}$,
$d_{\l,\yb}$,

\vskip2mm
\itemitem{2.6 :} 
$\th\!\Lc_\yb$, 
$\th\!\Lc^\delta_\yb$, 
$\th\Sb_\Vb$,
$\th\Zb_{\Lbb,\Vb}$,
$\th\Fb_{\Lbb,\Vb}$,
$1_{\Lbb,\Vb,\ib}$,

\vskip2mm
\itemitem{2.7 :} 
$\theta(\bb)\bb$,

\vskip2mm
\itemitem{2.8 :} 
$\th\Zb^\delta_{\Lbb,\Vb}$,

\vskip2mm
\itemitem{3 :} 
$\th\!\Fc_{\Lbb,\Vb}$,
$\th\!\Zc_{\Lbb,\Vb}$,

\vskip2mm
\itemitem{4.1 :} 
$G=O(\Vb,\varpi)$,
$F=F(\Vb,\varpi)$,
$T$,
$W$,
$W_\Vb$,

\vskip2mm
\itemitem{4.2 :} 
$\phi_\Vb$,
$\Db_l$,
$\Delta$,
$\Delta^+$,
$\Pi$,
$\th\!B_\Vb$,
$\th\!\Delta_\Vb$,

\vskip2mm
\itemitem{4.3 :} 
$\Sen_m$,
$W_m$,
$w(\ib)$,

\vskip2mm
\itemitem{4.4 :} 
$\phi_{\Vb,w}$,
$\ib_w$,
$W_\nu$,
$\th\!B_{\Vb,w}$,
$\th\!N_{\Vb,w}$,

\vskip2mm
\itemitem{4.5 :} 
$\th\!O^w_\Vb$,
$\th\!O^w_{\Vb,x,y}$,
$\th\!P_{\Vb,x,y}$,
$\th\! Z_{\Lbb,\Vb}^{x}$,

\vskip2mm
\itemitem{4.7 :} 
$\Sb$,
$\chi_l$,
$M[\l]$,
$\Eu(M)$,
$\Lambda_w$,
$\Lambda^x_{w,w'}$,

\vskip2mm
\itemitem{4.8 :} 
$\th\!\een_{\Lbb,\Vb,w}$,
$\th\!\nen_{\Lbb,\Vb,w}$,

\vskip2mm
\itemitem{4.9 :} 
$\th\!\een_{\Lbb,\Vb,w,w'}$,
$\th\!\den_{\Lbb,\Vb,w,w'}$,
$\th\!\nen_{\Lbb,\Vb,w,w'}$,
$\th\!\men_{\Lbb,\Vb,w,w'}$,

\vskip2mm
\itemitem{4.11 :} 
$x_\ib(l)$,

\vskip2mm
\itemitem{4.12 :} 
$\Qb$,
$\psi_w$,
$\psi_{w,w'}$,

\vskip2mm
\itemitem{4.14 :} 
$\l_\ib(l)$,
$h_\ib(k)$,
$\sigma_{\Lbb,\Vb}(k)$,
$\varkappa_{\Lbb,\Vb}(l)$,
$\pi_{\Lbb,\Vb}(1)$,
$\sigma_{\Lbb,\Vb,\ib,\ib'}(k)$,
$\varkappa_{\Lbb,\Vb,\ib,\ib'}(l)$,
$\pi_{\Lbb,\Vb,\ib,\ib'}(1)$,

\vskip2mm
\itemitem{4.16 :} 
$\sigma_{\Lbb,\Vb,\ib}(k)$,
$\varkappa_{\Lbb,\Vb,\ib}(l)$,
$\pi_{\Lbb,\Vb,\ib}(1)$,

\vskip2mm
\itemitem{5.1 :} 
$\th\Rb(\Gamma)_{\l,\nu}$,
$1_\ib$,
$\sigma_k$,
$\varkappa_l$,
$\pi_1$,
$Q_{i,j}(u,v)$,
$\omega$,

\vskip2mm
\itemitem{5.3 :} 
$\dot w$,
$\sigma_{\dot w}$,

\vskip2mm
\itemitem{6.1 :} 
$\Hb_m$,
$T_k$,
$X_l$,

\vskip2mm
\itemitem{6.2 :} 
$\varphi_k$,
$\Hb_{m}\text{-}\Modb_I$,

\vskip2mm
\itemitem{6.4 :} 
$\th\Rb_m$,
$\th\Rb_\nu$,
$1_{\nu,\nu'}$,
$1_{m,\nu'}$,
$\th\Rb_{m}\text{-}\Modb_0$,
$\th\Rb_{m}\text{-}\fModb_0$,
$\Hb_{m}\text{-}\fModb_I$,
$\Psi$,

\vskip2mm
\itemitem{6.9 :} 
$E_i$,
$F_i$,
$\kb_i$,

\vskip2mm
\itemitem{7.1 :} 
$\Rb_m$,
$\omega$,
$\tau$,
$\iota$,
$\kappa$,
$w_m$,

\vskip2mm
\itemitem{7.2 :} 
$\Rb_{m,m'}$,
$\phi_!$,
$\phi^*$,
$\phi_*$,
$P^\sharp$,
$P^\omega$,
$\Bc$,
$(\bullet :\bullet )$,
$\Kb_I$,
$\Gb_I$,
$\la\bullet :\bullet \ra$,
$M^\flat$,
$\Bc I^m$,
$\ch(M)$,

\vskip2mm
\itemitem{7.4 :} 
$\Rb_\yb$,

\vskip2mm
\itemitem{7.5 :} 
$\Lb_{\ib}$,
$\Lb_\yb$,
$\Lb_{mi}$,

\vskip2mm
\itemitem{7.6 :} 
$\Kc$,
$\theta_i$,
$\theta_i^{(a)}$,
$\theta_\yb$,
$r$,
$\fb$,
${}_\Ac\fb$,
$\fb_\nu$,
$\Gb^\up$,
$\Gb^\low$,
$B(\infty)$,
$(\bullet :\bullet )$,

\vskip2mm
\itemitem{8.1 :} 
$\th\Kb_I$,
$\th\Gb_I$,
$P^\sharp$,
$M^\flat$,
$\bar f$,
$\th\Gb^\low(\l)$,
$\th\Gb^\up(\l)$,
$\th\!B(\l)$,
$\la\bullet:\bullet\ra$,
$(\bullet:\bullet)$,

\vskip2mm
\itemitem{8.5 :} 
$\th\Lb_\ib$,

\vskip2mm
\itemitem{8.6 :} 
$D_{m,m'}$,
$W_{m,m'}$,
$D_{m,m';n,n'}$,
$W(w)$,
$\th\Rb_{m,m'}$,
$\psi_!$,
$\psi^*$,
$\psi_*$,
$e_i$,
$e'_i$,
$f_i$,

\vskip2mm
\itemitem{8.10 :} 
$\th\Rb_\ib$,
$\th\Rb_\yb$,
$\th\!\Bc I^m$,
$\ch(M)$,
$\deg(\ib,\ib';\,\ib'')$,
$Sh(\ib,\ib')$,

\vskip2mm
\itemitem{8.16 :} 
$N^\kappa$,
$\forb$,
$E_i$,
$F_i$,

\vskip2mm
\itemitem{8.20 :} 
$\tilde e_i$,
$\tilde f_i$,
$\eps_i$,
$\Delta_{ni}$,
$\tilde E_i$,
$\tilde F_i$,

\vskip2mm
\itemitem{8.29 :} 
$\th\Bcb$,
$\th\Vb(\l)$,

\vskip2mm
\itemitem{9.1 :} 
$\Delta^+_s$,
$\Delta^+_l$,
$\Delta^+$,

\vskip2mm
\itemitem{9.3 :} 
$\th\Rb_\nu^{\leqslant x}$,

\vskip2mm
\itemitem{10.1 :} 
$\Omega$,
$L_{\Lbb,\Vb,\Omega}$,
$\th\!E_{\Lbb,\Vb,\Omega}$,
$\th\!E_{\nu,\Omega}$,
$\th\!F_{\yb}$,
$\th\!\widetilde F_{\yb,\Omega}$,
$\th\!\Pcb$,
$\th\!\Qcb$,
$\th\!\Lc_{\yb}$,

\vskip2mm
\itemitem{10.3 :} 
$\th\!E_{1,\Omega}$,
$\th\!E_{2,\Omega}$,
$p_1$, $p_2$, $p_3$,
$\underline{f_i}$,

\vskip2mm
\itemitem{10.6 :} 
$\th\!E_{3,\Omega}$,
$\kappa$, $\iota$,
$\underline{e_i}$,

\vskip2mm
\itemitem{10.10 :} 
$\underline{\widetilde E_i}$,
$\underline{\widetilde F_i}$,
$\underline\eps_i$,

\vskip2mm
\itemitem{10.15 :} 
$\Yb$,

\vskip2mm
\itemitem{10.18 :} 
$\tilde\eb_i$,
$\tilde\fb_i$,
$\th\Lb(\l)$,
$\Rc$,

\enddocument